\documentclass{cmslatex}
\usepackage{amsfonts}
\usepackage{amsmath}
\usepackage{tikz}
\usetikzlibrary{trees}
\usetikzlibrary{arrows,shapes,snakes,automata,backgrounds,petri}

\newtheorem{thm}{Theorem}[section]

\newtheorem{lem}[thm]{Lemma}
\newtheorem{cor}[thm]{Corollary}

\begin{document}

\title{Graph-theoretic approaches to injectivity and multiple equilibria in systems of interacting elements}

\author{Murad Banaji
\thanks{Corresponding author. Department of Mathematics, University College London, Gower Street, London WC1E 6BT, UK.}
\and Gheorghe Craciun
\thanks{Department of Mathematics and Department of Biomolecular Chemistry, University of Wisconsin, Madison, WI 53706-1388.}
}

\maketitle

\begin{abstract}
We extend previous work on injectivity in chemical reaction networks to general interaction networks. Matrix- and graph-theoretic conditions for injectivity of these systems are presented. A particular signed, directed, labelled, bipartite multigraph, termed the ``DSR graph'', is shown to be a useful representation of an interaction network when discussing questions of injectivity. A graph-theoretic condition, developed previously in the context of chemical reaction networks, is shown to be sufficient to guarantee injectivity for a large class of systems. The graph-theoretic condition is simple to state and often easy to check. Examples are presented to illustrate the wide applicability of the theory developed. 
\end{abstract}


\section{Introduction}

\nopagebreak
Dynamical systems involving networks of interacting elements arise in many fields, including chemistry, systems biology, ecosystem modelling, and even beyond the natural sciences. The  description of these systems requires enumeration of the {\bf species} involved and the {\bf interactions} between them. Associated with a species is a ``state'', usually a population or a concentration, and associated with an interaction is a ``rate'' which describes the frequency of the interaction. An important theoretical question is: ``what claims can be made on the basis only of qualitative knowledge of the interactions, i.e. without detailed knowledge of functional forms and parameter values?''. In this paper we are concerned with questions of {\em injectivity} and hence the ability of systems to have multiple equilibria. In particular: 
\begin{enumerate}
\item Given a wide class of systems (to be defined below) we show how to associate with each such system an object termed a DSR graph;
\item We then show how a simple computation, or even observation, on this DSR graph can suffice to tell us that a system cannot have multiple equilibria (because the associated vector field is injective).
\end{enumerate}

The results presented here largely build on those in \cite{banajiSIAM,banajicraciun}. However the graph-theoretic work in this paper is closely connected to two previous strands of work: discussions of injectivity/multiple equilibria via {\bf interaction graphs} \cite{gouze98,soule,cinquin,kaufman}, and discussions of injectivity in chemical reaction networks (CRNs) with reference to so-called {\bf SR graphs} \cite{craciun1, craciun3}. The classes of systems treated in these references are special cases of interaction networks, as defined generally here.

There is also other important work in this tradition. For example, graph-theoretic approaches to questions of monotonicity appear in \cite{angelileenheersontag,kunzesiegelpositivity} and several results in \cite{banajidynsys} have easy graph-theoretic interpretations. Important questions of persistence in CRNs are treated graph-theoretically in \cite{angelipetrinet, anderson08}.

The main result is that given any interaction network it is possible to construct a multigraph called the {\bf directed SR graph} (DSR graph), and check a condition on this graph which will rule out the possibility of multiple equilibria. At a practical level, the DSR graph is an intuitively meaningful object, closely related to the diagrams drawn by researchers in biochemistry and chemical engineering. The condition sufficient for injectivity, here termed {\bf Condition~($*$)}, is simple to state and has been previously presented \cite{craciun1,banajicraciun}. What is remarkable is that a condition originally developed for CRNs with mass action kinetics \cite{craciun1}, and then generalised to CRNs with looser kinetics \cite{banajicraciun}, can actually be applied to arbitrary networks of interacting elements.

Before proving any results, we will discuss the notion of an ``interaction network'', both informally and formally. We will see that essentially arbitrary dynamical systems on appropriate subsets of $\mathbb{R}^n$ can be cast as interaction networks, and thus the results in this paper have wide applicability.

\subsection{Species and interactions: an informal discussion}

\nopagebreak
A species $S$ may participate in an interaction $R$ in the following three ways: 
\begin{enumerate}
\item  {\bf Two-way} ($S \leftrightarrow R$): the species influences the interaction and is itself influenced by the interaction. 
\item {\bf Species-to-interaction} ($S \rightarrow R$): the species influences the interaction, but is unaffected by the interaction. 
\item {\bf Interaction-to-species} ($S \leftarrow R$): the species is influenced by the interaction, but does not affect the interaction. 
\end{enumerate}
Case (1) describes perhaps the most common situation. During predation, for example, we expect the population of prey both to affect the rate of predation and to be affected by it. In case (2) the species is a {\bf modulator} of the interaction, and this modulation may be in a defined direction (``activation'' or ``inhibition'') or an undefined direction. For example, a simplification common in CRN modelling is to treat enzymes as modulators of reactions. In case (3), we say that the species participates {\bf irreversibly} in the interaction. Note that this is slightly different to the usual definition of irreversibility for a chemical reaction, because by this definition, if a species $S$ is the product of an irreversible reaction and it affects the rate of that reaction, then we say that $S$ participates in that reaction {\em reversibly}, even though the reaction could never run backwards (see for example the model of the TCA cycle presented later). 

A key feature of the DSR graph associated with an interaction network is that it has two types of vertices, one representing the species, and another representing the interactions. A species vertex $S$ and interaction vertex $R$ may be linked by up to two edges, which may be directed or undirected. These encode information about how species $S$ affects the rate of interaction $R$, and conversely how the dynamics of $S$ are affected by the rate of interaction. A formal treatment will be developed later, but the reader might like to glance at Appendix~\ref{appDSR} to get a feel for the meanings of various {\em motifs} in DSR graphs. 

\subsection{Summary of the results}

\nopagebreak
{\bf Rectangular domains.} A subset $X$ of $\mathbb{R}^n$ is a rectangular domain iff for $i = 1, \ldots, n$ there exist quantities $-\infty \leq x_{i, min} < \infty$ and $-\infty < x_{i, max} \leq \infty$ satisfying $x_{i, min} \leq x_{i, max}$, and for each $x \equiv [x_1, \ldots, x_n]^T \in X$
\[
x_{i, min} \prec x_i \prec x_{i, max} \quad i =1, \ldots, n\,.
\]
Here the symbol $\prec$ can mean $<$ or $\leq$. For example $X$ could be all of $\mathbb{R}^n$, an orthant, etc. Note that $X$ may be of lower dimension than $\mathbb{R}^n$. 

{\bf Injectivity and multiple equilibria.} A function $f:X \to Y$ between arbitrary spaces $X$ and $Y$ is {\bf injective} iff whenever $a, b \in X$, $a \not = b$, then $f(a) \not = f(b)$. Let $X \subset \mathbb{R}^n$, and $f:X \to \mathbb{R}^n$ an arbitrary $C^1$ (i.e. continuously differentiable) function on $X$. Supposing some model of the natural world gives rise to a vector field $\dot x = f(x)$, or the map $x_{n+1} = x_n + f(x_n)$, where $x \in X$. Suppose further that we have incomplete data, so that all we know is that $f$ belongs to some set of functions $\mathcal{F}$, but by some means we can show that all $f \in \mathcal{F}$ are injective. This immediately implies that each $f \in \mathcal{F}$ can have no more than one zero, i.e. the vector field can have no more than one equilibrium, or the map no more than one fixed point.

{\bf The key idea: decomposing functions.} Our central theme in this paper is the following: we start with a set of functions $\mathcal{F}$ from some rectangular domain $X \subset \mathbb{R}^n$ to $\mathbb{R}^n$. Now for some integer $m \geq 1$, and for each $f \in \mathcal{F}$ we write $f = f_2 \circ f_1$, where $f_1:X \to \mathbb{R}^m$ and $f_2: \mathbb{R}^m \to \mathbb{R}^n$. Such decompositions always exist (for example we can choose $m=n$, $f_1 = f$, and $f_2$ to be the identity, or vice versa), but often there are natural decompositions implied by the way that the models are constructed. Intuitively, we can think of $X$ as ``species space'', the set of allowed species values, and $\mathbb{R}^m$ as ``interaction space'', the space of allowed interaction rates. Differentiating $f = f_2 \circ f_1$ gives $Df(x) = Df_2(f_1(x))Df_1(x)$, i.e. the Jacobian $Df$ has a natural product structure. Specifically, $(Df_2)_{ij}$ encodes information about how interaction $j$ affects affects quantity $i$, while $(Df_1)_{ji}$ encodes information about how quantity $i$ affects the rate of interaction $j$. 

{\bf $P_0^{(-)}$ matrices.} The product structure in the Jacobian forms the theoretical starting point for the subsequent treatment. In particular, the class of $P_0^{(-)}$ matrices (defined in Appendix~\ref{app1}) is shown to be a useful class to consider when asking questions about injectivity. Two very general lemmas are developed: Lemma~\ref{genlem} gives a sufficient condition for any set of matrices with product structure to be  $P_0^{(-)}$ matrices; Lemma~\ref{genlemconv} provides restrictions under which the sufficient conditions in Lemma~\ref{genlem} are also necessary. 

{\bf DSR graph results.} While the matrix-theoretic approaches give sharp results, there are very elegant and intuitive graph-theoretic results which arise as corollaries of these. Roughly speaking, the DSR graph is constructed as follows: returning to the decomposition of the Jacobian $Df(x) = Df_2(f_1(x))Df_1(x)$, the elements $(Df_2)_{ij}$ and $(Df_1)_{ji}$ together contain information on how species $i$ is affected by, and affects, interaction $j$, and unless both are identically zero, together they map to one or more edges between species vertex $i$ and interaction vertex $j$ in the DSR graph. Via a series of lemmas we find that a previously developed condition on cycles in the SR graph, restated for the DSR graph, guarantees injectivity of the interaction network with some outflow/degradation of each species. The implications for multiple equilibria under different assumptions are summarised in Corollary~\ref{maincor}, which can be seen as the main result in this paper.

\subsection{An example: the TCA cycle}

\nopagebreak
Before presenting the theoretical development, or even rigorously defining interaction networks and DSR graphs, it is helpful to present a nontrivial example. Further discussion of this and other examples will be provided after development of the theory. 

We consider the model of the TCA cycle discussed in \cite{cortassa}. The model is a fairly complicated biological model, and is presented in stages. At each stage the relevant question is whether it can admit ``multiple positive nondegenerate equilibria'', a notion to be made precise later. 

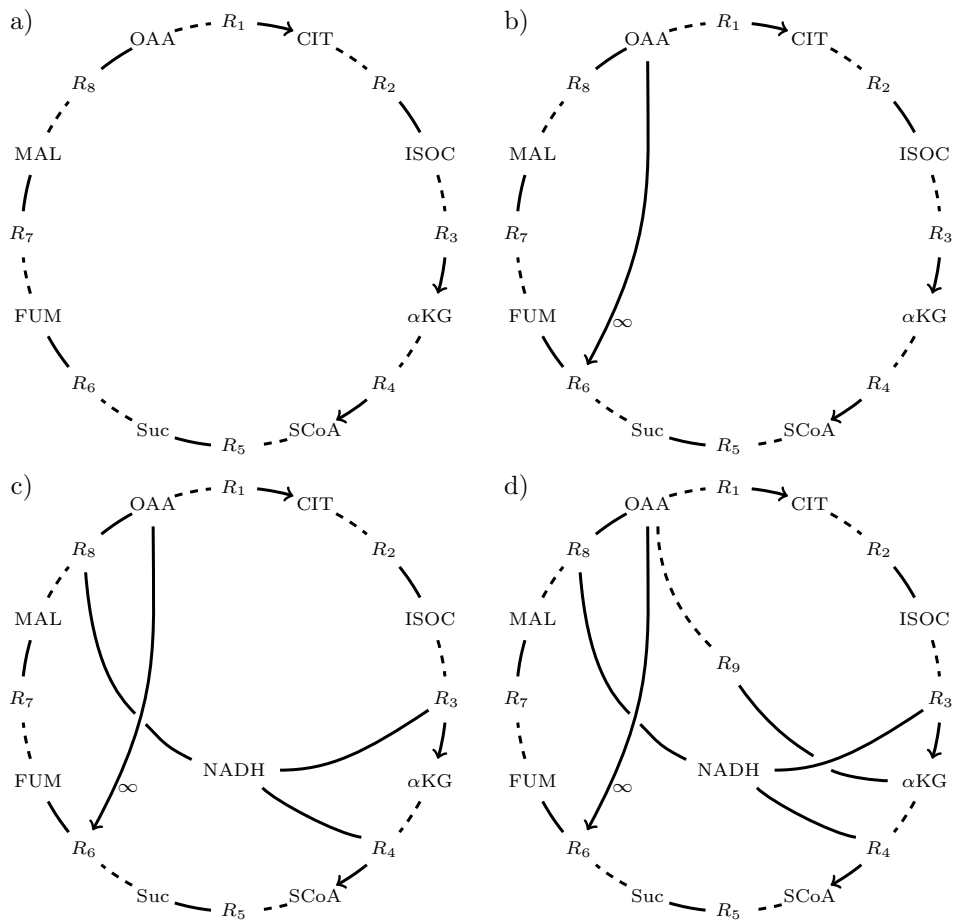
\begin{figure}[h]
\begin{minipage}{0.99\textwidth}
\begin{minipage}{0.49\textwidth}
\begin{tikzpicture}[domain=-6:6,scale=0.47]

\node at (-6,6) {a)};

\path (90: 6cm) coordinate (R1);
\node at (R1) {$\scriptstyle{R_1}$};
\path (67.5: 6cm) coordinate (CIT);
\node at (CIT) {$\scriptstyle{\mathrm{CIT}}$};
\path (45: 6cm) coordinate (R2);
\node at (R2) {$\scriptstyle{R_2}$};
\path (22.5: 6cm) coordinate (ISOC);
\node at (ISOC) {$\scriptstyle{\mathrm{ISOC}}$};
\path (0: 6cm) coordinate (R3);
\node at (R3) {$\scriptstyle{R_3}$};
\path (-22.5: 6cm) coordinate (AKG);
\node at (AKG) {$\scriptstyle{\alpha\mathrm{KG}}$};
\path (-45: 6cm) coordinate (R4);
\node at (R4) {$\scriptstyle{R_4}$};
\path (-67.5: 6cm) coordinate (SCoA);
\node at (SCoA) {$\scriptstyle{\mathrm{SCoA}}$};
\path (-90: 6cm) coordinate (R5);
\node at (R5) {$\scriptstyle{R_5}$};
\path (-112.5: 6cm) coordinate (Suc);
\node at (Suc) {$\scriptstyle{\mathrm{Suc}}$};
\path (-135: 6cm) coordinate (R6);
\node at (R6) {$\scriptstyle{R_6}$};
\path (-157.5: 6cm) coordinate (FUM);
\node at (FUM) {$\scriptstyle{\mathrm{FUM}}$};
\path (180: 6cm) coordinate (R7);
\node at (R7) {$\scriptstyle{R_7}$};
\path (157.5: 6cm) coordinate (MAL);
\node at (MAL) {$\scriptstyle{\mathrm{MAL}}$};
\path (135: 6cm) coordinate (R8);
\node at (R8) {$\scriptstyle{R_8}$};
\path (112.5: 6cm) coordinate (OAA);
\node at (OAA) {$\scriptstyle{\mathrm{OAA}}$};

\path (90-6.25: 6cm) coordinate (R1end);
\path (67.5-6.25: 6cm) coordinate (CITend);
\path (45-6.25: 6cm) coordinate (R2end);
\path (22.5-6.25: 6cm) coordinate (ISOCend);
\path (0-6.25: 6cm) coordinate (R3end);
\path (-22.5-6.25: 6cm) coordinate (AKGend);
\path (-45-6.25: 6cm) coordinate (R4end);
\path (-67.5-8: 6cm) coordinate (SCoAend);
\path (-90-6.25: 6cm) coordinate (R5end);
\path (-112.5-6.25: 6cm) coordinate (Sucend);
\path (-135-6.25: 6cm) coordinate (R6end);
\path (-157.5-6.25: 6cm) coordinate (FUMend);
\path (180-6.25: 6cm) coordinate (R7end);
\path (157.5-6.25: 6cm) coordinate (MALend);
\path (135-6.25: 6cm) coordinate (R8end);
\path (112.5-6.25: 6cm) coordinate (OAAend);

\draw[->, line width=0.04cm] (R1end) arc (90-6.25:90-16.25:6cm);
\draw[-, dashed, line width=0.04cm] (CITend)  arc (67.5-6.25:67.5-16.25:6cm);
\draw[-, line width=0.04cm] (R2end)  arc (45-6.25:45-16.25:6cm);
\draw[-, dashed, line width=0.04cm] (ISOCend)  arc (22.5-6.25:22.5-16.25:6cm);
\draw[->, line width=0.04cm] (R3end)  arc (-6.25:-16.25:6cm);
\draw[-, dashed, line width=0.04cm] (AKGend)  arc (-22.5-6.25:-22.5-16.25:6cm);
\draw[->, line width=0.04cm] (R4end)  arc (-45-6.25:-45-16.25:6cm);

\draw[-, dashed, line width=0.04cm] (SCoAend)  arc (-67.5-8:-67.5-16.25:6cm);
\draw[-, line width=0.04cm] (R5end)  arc (-90-6.25:-90-16.25:6cm);
\draw[-, dashed, line width=0.04cm] (Sucend)  arc (-112.5-6.25:-112.5-16.25:6cm);
\draw[-, line width=0.04cm] (R6end)  arc (-135-6.25:-135-16.25:6cm);

\draw[-, dashed, line width=0.04cm] (FUMend)  arc (-157.5-6.25:-157.5-16.25:6cm);
\draw[-, line width=0.04cm] (R7end)  arc (180-6.25:180-16.25:6cm);
\draw[-, dashed, line width=0.04cm] (MALend)  arc (157.5-6.25:157.5-16.25:6cm);
\draw[-, line width=0.04cm] (R8end)  arc (135-6.25:135-16.25:6cm);
\draw[-, dashed, line width=0.04cm] (OAAend)  arc (112.5-6.25:112.5-16.25:6cm);

\end{tikzpicture}
\end{minipage}
\hfill
\begin{minipage}{0.49\textwidth}
\begin{tikzpicture}[domain=-6:6,scale=0.47]

\node at (-6,6) {b)};

\path (90: 6cm) coordinate (R1);
\node at (R1) {$\scriptstyle{R_1}$};
\path (67.5: 6cm) coordinate (CIT);
\node at (CIT) {$\scriptstyle{\mathrm{CIT}}$};
\path (45: 6cm) coordinate (R2);
\node at (R2) {$\scriptstyle{R_2}$};
\path (22.5: 6cm) coordinate (ISOC);
\node at (ISOC) {$\scriptstyle{\mathrm{ISOC}}$};
\path (0: 6cm) coordinate (R3);
\node at (R3) {$\scriptstyle{R_3}$};
\path (-22.5: 6cm) coordinate (AKG);
\node at (AKG) {$\scriptstyle{\alpha\mathrm{KG}}$};
\path (-45: 6cm) coordinate (R4);
\node at (R4) {$\scriptstyle{R_4}$};
\path (-67.5: 6cm) coordinate (SCoA);
\node at (SCoA) {$\scriptstyle{\mathrm{SCoA}}$};
\path (-90: 6cm) coordinate (R5);
\node at (R5) {$\scriptstyle{R_5}$};
\path (-112.5: 6cm) coordinate (Suc);
\node at (Suc) {$\scriptstyle{\mathrm{Suc}}$};
\path (-135: 6cm) coordinate (R6);
\node at (R6) {$\scriptstyle{R_6}$};
\path (-157.5: 6cm) coordinate (FUM);
\node at (FUM) {$\scriptstyle{\mathrm{FUM}}$};
\path (180: 6cm) coordinate (R7);
\node at (R7) {$\scriptstyle{R_7}$};
\path (157.5: 6cm) coordinate (MAL);
\node at (MAL) {$\scriptstyle{\mathrm{MAL}}$};
\path (135: 6cm) coordinate (R8);
\node at (R8) {$\scriptstyle{R_8}$};
\path (112.5: 6cm) coordinate (OAA);
\node at (OAA) {$\scriptstyle{\mathrm{OAA}}$};

\path (90-6.25: 6cm) coordinate (R1end);
\path (67.5-6.25: 6cm) coordinate (CITend);
\path (45-6.25: 6cm) coordinate (R2end);
\path (22.5-6.25: 6cm) coordinate (ISOCend);
\path (0-6.25: 6cm) coordinate (R3end);
\path (-22.5-6.25: 6cm) coordinate (AKGend);
\path (-45-6.25: 6cm) coordinate (R4end);
\path (-67.5-8: 6cm) coordinate (SCoAend);
\path (-90-6.25: 6cm) coordinate (R5end);
\path (-112.5-6.25: 6cm) coordinate (Sucend);
\path (-135-6.25: 6cm) coordinate (R6end);
\path (-157.5-6.25: 6cm) coordinate (FUMend);
\path (180-6.25: 6cm) coordinate (R7end);
\path (157.5-6.25: 6cm) coordinate (MALend);
\path (135-6.25: 6cm) coordinate (R8end);
\path (112.5-6.25: 6cm) coordinate (OAAend);

\draw[->, line width=0.04cm] (R1end) arc (90-6.25:90-16.25:6cm);
\draw[-, dashed, line width=0.04cm] (CITend)  arc (67.5-6.25:67.5-16.25:6cm);
\draw[-, line width=0.04cm] (R2end)  arc (45-6.25:45-16.25:6cm);
\draw[-, dashed, line width=0.04cm] (ISOCend)  arc (22.5-6.25:22.5-16.25:6cm);
\draw[->, line width=0.04cm] (R3end)  arc (-6.25:-16.25:6cm);
\draw[-, dashed, line width=0.04cm] (AKGend)  arc (-22.5-6.25:-22.5-16.25:6cm);
\draw[->, line width=0.04cm] (R4end)  arc (-45-6.25:-45-16.25:6cm);

\draw[-, dashed, line width=0.04cm] (SCoAend)  arc (-67.5-8:-67.5-16.25:6cm);
\draw[-, line width=0.04cm] (R5end)  arc (-90-6.25:-90-16.25:6cm);
\draw[-, dashed, line width=0.04cm] (Sucend)  arc (-112.5-6.25:-112.5-16.25:6cm);
\draw[-, line width=0.04cm] (R6end)  arc (-135-6.25:-135-16.25:6cm);

\draw[-, dashed, line width=0.04cm] (FUMend)  arc (-157.5-6.25:-157.5-16.25:6cm);
\draw[-, line width=0.04cm] (R7end)  arc (180-6.25:180-16.25:6cm);
\draw[-, dashed, line width=0.04cm] (MALend)  arc (157.5-6.25:157.5-16.25:6cm);
\draw[-, line width=0.04cm] (R8end)  arc (135-6.25:135-16.25:6cm);
\draw[-, dashed, line width=0.04cm] (OAAend)  arc (112.5-6.25:112.5-16.25:6cm);

\draw[->, line width=0.04cm] (3.71-6, 11.2-6.3) .. controls (3.71-6,7.54-6) and (4-6,0) .. (2.0-6,2.0-5.7);
\node at (3-6, 3.5-6) {$\scriptstyle{\infty}$};

\end{tikzpicture}
\end{minipage}
\end{minipage}
\begin{minipage}{0.99\textwidth}
\begin{minipage}{0.49\textwidth}
\begin{tikzpicture}[domain=-6:6,scale=0.47]

\node at (-6,6) {c)};

\path (90: 6cm) coordinate (R1);
\node at (R1) {$\scriptstyle{R_1}$};
\path (67.5: 6cm) coordinate (CIT);
\node at (CIT) {$\scriptstyle{\mathrm{CIT}}$};
\path (45: 6cm) coordinate (R2);
\node at (R2) {$\scriptstyle{R_2}$};
\path (22.5: 6cm) coordinate (ISOC);
\node at (ISOC) {$\scriptstyle{\mathrm{ISOC}}$};
\path (0: 6cm) coordinate (R3);
\node at (R3) {$\scriptstyle{R_3}$};
\path (-22.5: 6cm) coordinate (AKG);
\node at (AKG) {$\scriptstyle{\alpha\mathrm{KG}}$};
\path (-45: 6cm) coordinate (R4);
\node at (R4) {$\scriptstyle{R_4}$};
\path (-67.5: 6cm) coordinate (SCoA);
\node at (SCoA) {$\scriptstyle{\mathrm{SCoA}}$};
\path (-90: 6cm) coordinate (R5);
\node at (R5) {$\scriptstyle{R_5}$};
\path (-112.5: 6cm) coordinate (Suc);
\node at (Suc) {$\scriptstyle{\mathrm{Suc}}$};
\path (-135: 6cm) coordinate (R6);
\node at (R6) {$\scriptstyle{R_6}$};
\path (-157.5: 6cm) coordinate (FUM);
\node at (FUM) {$\scriptstyle{\mathrm{FUM}}$};
\path (180: 6cm) coordinate (R7);
\node at (R7) {$\scriptstyle{R_7}$};
\path (157.5: 6cm) coordinate (MAL);
\node at (MAL) {$\scriptstyle{\mathrm{MAL}}$};
\path (135: 6cm) coordinate (R8);
\node at (R8) {$\scriptstyle{R_8}$};
\path (112.5: 6cm) coordinate (OAA);
\node at (OAA) {$\scriptstyle{\mathrm{OAA}}$};

\path (90-6.25: 6cm) coordinate (R1end);
\path (67.5-6.25: 6cm) coordinate (CITend);
\path (45-6.25: 6cm) coordinate (R2end);
\path (22.5-6.25: 6cm) coordinate (ISOCend);
\path (0-6.25: 6cm) coordinate (R3end);
\path (-22.5-6.25: 6cm) coordinate (AKGend);
\path (-45-6.25: 6cm) coordinate (R4end);
\path (-67.5-8: 6cm) coordinate (SCoAend);
\path (-90-6.25: 6cm) coordinate (R5end);
\path (-112.5-6.25: 6cm) coordinate (Sucend);
\path (-135-6.25: 6cm) coordinate (R6end);
\path (-157.5-6.25: 6cm) coordinate (FUMend);
\path (180-6.25: 6cm) coordinate (R7end);
\path (157.5-6.25: 6cm) coordinate (MALend);
\path (135-6.25: 6cm) coordinate (R8end);
\path (112.5-6.25: 6cm) coordinate (OAAend);

\draw[->, line width=0.04cm] (R1end) arc (90-6.25:90-16.25:6cm);
\draw[-, dashed, line width=0.04cm] (CITend)  arc (67.5-6.25:67.5-16.25:6cm);
\draw[-, line width=0.04cm] (R2end)  arc (45-6.25:45-16.25:6cm);
\draw[-, dashed, line width=0.04cm] (ISOCend)  arc (22.5-6.25:22.5-16.25:6cm);
\draw[->, line width=0.04cm] (R3end)  arc (-6.25:-16.25:6cm);
\draw[-, dashed, line width=0.04cm] (AKGend)  arc (-22.5-6.25:-22.5-16.25:6cm);
\draw[->, line width=0.04cm] (R4end)  arc (-45-6.25:-45-16.25:6cm);

\draw[-, dashed, line width=0.04cm] (SCoAend)  arc (-67.5-8:-67.5-16.25:6cm);
\draw[-, line width=0.04cm] (R5end)  arc (-90-6.25:-90-16.25:6cm);
\draw[-, dashed, line width=0.04cm] (Sucend)  arc (-112.5-6.25:-112.5-16.25:6cm);
\draw[-, line width=0.04cm] (R6end)  arc (-135-6.25:-135-16.25:6cm);

\draw[-, dashed, line width=0.04cm] (FUMend)  arc (-157.5-6.25:-157.5-16.25:6cm);
\draw[-, line width=0.04cm] (R7end)  arc (180-6.25:180-16.25:6cm);
\draw[-, dashed, line width=0.04cm] (MALend)  arc (157.5-6.25:157.5-16.25:6cm);
\draw[-, line width=0.04cm] (R8end)  arc (135-6.25:135-16.25:6cm);
\draw[-, dashed, line width=0.04cm] (OAAend)  arc (112.5-6.25:112.5-16.25:6cm);

\draw[->, line width=0.04cm] (3.71-6, 11.2-6.3) .. controls (3.71-6,7.54-6) and (4-6,0) .. (2.0-6,2.0-5.7);
\node at (3-6, 3.5-6) {$\scriptstyle{\infty}$};

\path (0,-2) coordinate (NADH);
\node at (NADH) {$\scriptstyle{\mathrm{NADH}}$};
\draw[-, line width=0.04cm] (1.3,-2) .. controls (8.7-6,4-6) and (9.7-6,4.5-6) .. (11.7-6.2, 5.8-6.1);
\draw[-, line width=0.04cm] (0.8,3.7-6.2) .. controls (7.2-6,3.1-6) and (8.9-6,2.2-6) .. (9.9-6.3,1.8-5.7);
\draw[-, line width=0.04cm] (5.1-6.3,4.1-5.8) .. controls (4-6,4.7-6) and (4-6,4.8-6) .. (3.5-6,5.3-6);
\draw[-, line width=0.04cm] (3.2-6,5.6-6) .. controls (2.7-6,6.1-6) and (2-6,7-6) .. (1.8-6,9.9-6.3);

\end{tikzpicture}
\end{minipage}
\hfill
\begin{minipage}{0.49\textwidth}
\begin{tikzpicture}[domain=-6:6,scale=0.47]

\node at (-6,6) {d)};

\path (90: 6cm) coordinate (R1);
\node at (R1) {$\scriptstyle{R_1}$};
\path (67.5: 6cm) coordinate (CIT);
\node at (CIT) {$\scriptstyle{\mathrm{CIT}}$};
\path (45: 6cm) coordinate (R2);
\node at (R2) {$\scriptstyle{R_2}$};
\path (22.5: 6cm) coordinate (ISOC);
\node at (ISOC) {$\scriptstyle{\mathrm{ISOC}}$};
\path (0: 6cm) coordinate (R3);
\node at (R3) {$\scriptstyle{R_3}$};
\path (-22.5: 6cm) coordinate (AKG);
\node at (AKG) {$\scriptstyle{\alpha\mathrm{KG}}$};
\path (-45: 6cm) coordinate (R4);
\node at (R4) {$\scriptstyle{R_4}$};
\path (-67.5: 6cm) coordinate (SCoA);
\node at (SCoA) {$\scriptstyle{\mathrm{SCoA}}$};
\path (-90: 6cm) coordinate (R5);
\node at (R5) {$\scriptstyle{R_5}$};
\path (-112.5: 6cm) coordinate (Suc);
\node at (Suc) {$\scriptstyle{\mathrm{Suc}}$};
\path (-135: 6cm) coordinate (R6);
\node at (R6) {$\scriptstyle{R_6}$};
\path (-157.5: 6cm) coordinate (FUM);
\node at (FUM) {$\scriptstyle{\mathrm{FUM}}$};
\path (180: 6cm) coordinate (R7);
\node at (R7) {$\scriptstyle{R_7}$};
\path (157.5: 6cm) coordinate (MAL);
\node at (MAL) {$\scriptstyle{\mathrm{MAL}}$};
\path (135: 6cm) coordinate (R8);
\node at (R8) {$\scriptstyle{R_8}$};
\path (112.5: 6cm) coordinate (OAA);
\node at (OAA) {$\scriptstyle{\mathrm{OAA}}$};

\path (90-6.25: 6cm) coordinate (R1end);
\path (67.5-6.25: 6cm) coordinate (CITend);
\path (45-6.25: 6cm) coordinate (R2end);
\path (22.5-6.25: 6cm) coordinate (ISOCend);
\path (0-6.25: 6cm) coordinate (R3end);
\path (-22.5-6.25: 6cm) coordinate (AKGend);
\path (-45-6.25: 6cm) coordinate (R4end);
\path (-67.5-8: 6cm) coordinate (SCoAend);
\path (-90-6.25: 6cm) coordinate (R5end);
\path (-112.5-6.25: 6cm) coordinate (Sucend);
\path (-135-6.25: 6cm) coordinate (R6end);
\path (-157.5-6.25: 6cm) coordinate (FUMend);
\path (180-6.25: 6cm) coordinate (R7end);
\path (157.5-6.25: 6cm) coordinate (MALend);
\path (135-6.25: 6cm) coordinate (R8end);
\path (112.5-6.25: 6cm) coordinate (OAAend);

\draw[->, line width=0.04cm] (R1end) arc (90-6.25:90-16.25:6cm);
\draw[-, dashed, line width=0.04cm] (CITend)  arc (67.5-6.25:67.5-16.25:6cm);
\draw[-, line width=0.04cm] (R2end)  arc (45-6.25:45-16.25:6cm);
\draw[-, dashed, line width=0.04cm] (ISOCend)  arc (22.5-6.25:22.5-16.25:6cm);
\draw[->, line width=0.04cm] (R3end)  arc (-6.25:-16.25:6cm);
\draw[-, dashed, line width=0.04cm] (AKGend)  arc (-22.5-6.25:-22.5-16.25:6cm);
\draw[->, line width=0.04cm] (R4end)  arc (-45-6.25:-45-16.25:6cm);

\draw[-, dashed, line width=0.04cm] (SCoAend)  arc (-67.5-8:-67.5-16.25:6cm);
\draw[-, line width=0.04cm] (R5end)  arc (-90-6.25:-90-16.25:6cm);
\draw[-, dashed, line width=0.04cm] (Sucend)  arc (-112.5-6.25:-112.5-16.25:6cm);
\draw[-, line width=0.04cm] (R6end)  arc (-135-6.25:-135-16.25:6cm);

\draw[-, dashed, line width=0.04cm] (FUMend)  arc (-157.5-6.25:-157.5-16.25:6cm);
\draw[-, line width=0.04cm] (R7end)  arc (180-6.25:180-16.25:6cm);
\draw[-, dashed, line width=0.04cm] (MALend)  arc (157.5-6.25:157.5-16.25:6cm);
\draw[-, line width=0.04cm] (R8end)  arc (135-6.25:135-16.25:6cm);
\draw[-, dashed, line width=0.04cm] (OAAend)  arc (112.5-6.25:112.5-16.25:6cm);

\draw[->, line width=0.04cm] (3.71-6, 11.2-6.3) .. controls (3.71-6,7.54-6) and (4-6,0) .. (2.0-6,2.0-5.7);
\node at (3-6, 3.5-6) {$\scriptstyle{\infty}$};

\path (0,-2) coordinate (NADH);
\node at (NADH) {$\scriptstyle{\mathrm{NADH}}$};
\draw[-, line width=0.04cm] (1.3,-2) .. controls (8.7-6,4-6) and (9.7-6,4.5-6) .. (11.7-6.2, 5.8-6.1);
\draw[-, line width=0.04cm] (0.8,3.7-6.2) .. controls (7.2-6,3.1-6) and (8.9-6,2.2-6) .. (9.9-6.3,1.8-5.7);
\draw[-, line width=0.04cm] (5.1-6.3,4.1-5.8) .. controls (4-6,4.7-6) and (4-6,4.8-6) .. (3.5-6,5.3-6);
\draw[-, line width=0.04cm] (3.2-6,5.6-6) .. controls (2.7-6,6.1-6) and (2-6,7-6) .. (1.8-6,9.9-6.3);

\path (0,1) coordinate (R9);
\node at (R9) {$\scriptstyle{R_9}$};
\draw[dashed, line width=0.04cm] (-2, 11.2-6.3) .. controls (4-6,9-6) and (5.2-6,7.8-6) .. (5.7-6.2,7.3-5.8);
\draw[line width=0.04cm] (6.2-5.9, 6.6-6.2) .. controls (7.2-6,5.1-6) and (8.1-6,4.6-6) .. (8.4-6,4.4-6);
\draw[line width=0.04cm] (8.9-6, 4-6) .. controls (9.3-6,3.8-6) and (10.2-6,3.7-6) .. (10.8-6.3,3.7-6);

\end{tikzpicture}
\end{minipage}
\end{minipage}
\caption{\label{SRTCA} The DSR graph for the model of the TCA cycle presented in \cite{cortassa} at four different stages of construction. Details of the biochemistry can be found in \cite{cortassa}. Negative edges are represented with dashed lines, while positive edges are represented with bold lines. Apart from the one edge labelled $\infty$ all edges have edge-label $1$. Various quantities which do not affect the cycle structure have been omitted from the DSR graphs. Implications are described in the text. }
\end{figure}

The backbone of the model is a cycle of eight interconversions, that is, very basic chemical reactions in which one substrate is simply converted into another\footnote{The analysis remains the same if we treat it as a cycle of nine interconversions, including cis-aconitate.}. Importantly, some of the reactions are assumed to be irreversible while others are reversible. The DSR graph for the model at four different levels of complexity is shown in Figure~\ref{SRTCA}.

\begin{enumerate}
\item The basic structure of the model gives the DSR graph shown in Figure~\ref{SRTCA}a. Models with this simple structure cannot display multiple positive nondegenerate equilibria. In fact previous theory \cite{banajimathchem} indicates that with mild additional assumptions, there must be a unique, globally stable, equilibrium. 
\item Adding the inhibition by oxaloacetate of succinate dehydrogenase gives rise to the DSR graph shown in Figure~\ref{SRTCA}b. Analysis of this graph tells us that the conclusion about multiple nondegenerate equilibria remains true. 
\item Several reactions in the TCA cycle cause reduction of NAD. As $[\mathrm{NAD}] + [\mathrm{NADH}]$ is a conserved quantity, one of the pair can be eliminated. Including NADH gives the DSR graph shown in Figure~\ref{SRTCA}c. Again, the model does not allow multiple positive nondegenerate equilibria, despite the large number of cycles in the DSR graph. 
\item Finally, adding the aspartate amino transferase (AAT)-catalysed reaction effectively adds an extra interconversion between $\alpha$-ketoglutarate and oxaloacetate giving rise to the DSR graph shown in Figure~\ref{SRTCA}d. Now the DSR graph contains ``bad cycles'' and it is no longer possible to rule out multiple positive nondegenerate equilibria from the graph. 
\end{enumerate}

\section{Interaction networks}

\nopagebreak
Assume that there are $n$ species with ``amounts'' (concentrations, populations, etc.) $x_1, \ldots, x_n$, and define $x = [x_1, \ldots, x_n]^T$. Let there be $m$ interactions which occur with rates $v_1(x), \ldots, v_m(x)$, each involving any subset of the species, and define $v(x) = [v_1(x), \ldots, v_m(x)]^T$. Finally, define the $i$th {\bf interaction function}, $f_{i}(v(x))$, to be the total rate of production/consumption of species $i$. The evolution of the system is then given by:
\begin{equation}
\label{genprodform}
\dot x_i = f_{i}(v(x)) \qquad i = 1, \ldots, n.
\end{equation}
It is assumed that all functions are $C^1$ and that each $x_i$ takes values in some interval (bounded or unbounded), so that the state space is a rectangular domain in $\mathbb{R}^n$. Two absolutely general features of (\ref{genprodform}) are: 
\begin{itemize}
\item By the chain rule, the Jacobian (at each point in space) allows a decomposition $J = SV$ where $S_{ij} = \frac{\partial f_{i}}{\partial v_j}$, and $V_{ji} = \frac{\partial v_j}{\partial x_i}$. 
\item The interaction structure can be represented as a bipartite multigraph (discussed informally above, and to be defined formally later) with $n$ species vertices and $m$ interaction vertices.  
\end{itemize}

A slight, but important, variant on Equation~(\ref{genprodform}) involves assuming nonzero outflow or degradation rates for each species. We introduce a set of scalar $C^1$ functions $q_i(x_i)$, $i = 1, \ldots, n$, each satisfying $\frac{\partial q_i}{\partial x_i} > 0$ throughout its domain of definition. The system becomes
\begin{equation}
\label{genprodformout}
\dot x_i = f_{i}(v(x)) -q_i(x_i)\qquad i = 1, \ldots, n\,.
\end{equation}
The Jacobian is now $J = SV - D$ where $D$ is the positive diagonal matrix defined by $D_{ii} \equiv \frac{\partial q_i}{\partial x_i}$. Defining 
\[
f(v(x)) = [f_1(v(x)), \ldots, f_n(v(x))]^T \qquad \mbox{and} \qquad Q(x) = [q_1(x_1), \ldots, q_n(x_n)]^T\,,
\] 
we can write (\ref{genprodform})~and~(\ref{genprodformout}) in abbreviated forms which we term (\ref{genprodform1}) and (\ref{genprodformout1}):
\begin{equation*}
\tag{$N_0$}
\label{genprodform1}
\dot x = f(v(x))
\end{equation*}
and
\begin{equation*}
\tag{$N_{+}$}
\label{genprodformout1}
\dot x = f(v(x)) - Q(x).
\end{equation*}
The situation where some components of $Q$ are identically zero is also important. Given any $\theta \subset \{1, \ldots, n\}$ we can define a class of systems which we term (\ref{genprodformconserved}):
\begin{equation*}
\tag{$N_\theta$}
\label{genprodformconserved}
\dot x = f(v(x)) - Q_\theta(x),
\end{equation*}
with $Q_\theta(x) = [q_1(x_1), \ldots, q_n(x_n)]^T$, $q_i(x_i) = 0$ for $i \not \in \theta$ and $\frac{\partial q_i}{\partial x_i} > 0$ for $i \in \theta$. Note that both (\ref{genprodform1}) and (\ref{genprodformout1}) are now special cases of (\ref{genprodformconserved}), with $\theta = \emptyset$ and $\theta = \{1, \ldots, n\}$ respectively.

\subsection{Arbitrary dynamical systems as interaction networks}

\nopagebreak
Consider a rectangular domain $X \subset \mathbb{R}^n$, and some $C^1$ function $F:X \to \mathbb{R}^n$. With $Q(x)$  and $Q_\theta(x)$ defined as above, consider the dynamical systems: 
\begin{equation}
\label{verygeneqn}
\dot x = F(x), \quad \dot x = F(x) - Q(x) \quad \mbox{and} \quad \dot x = F(x) - Q_\theta(x). 
\end{equation}
Defining $\mathrm{id}(\cdot)$ to be the identity on $\mathbb{R}^n$, we can certainly write $F = F \circ \mathrm{id}$, and the three systems can be rewritten
\begin{equation}
\label{verygeneqn1}
\dot x = F(\mathrm{id}(x)), \quad \dot x = F(\mathrm{id}(x)) - Q(x) \quad \mbox{and} \quad \dot x = F(\mathrm{id}(x)) - Q_\theta(x).
\end{equation}
Thus the three systems in (\ref{verygeneqn}) can be cast as instances of (\ref{genprodform1}), (\ref{genprodformout1}) and (\ref{genprodformconserved}), and all theory developed in this paper can be applied. Of course the choice $F = F \circ \mathrm{id}$ is not unique. We will see by example later that this formal approach may be useful, but does not necessarily give the strongest results. This theme will be more fully developed in future work.

\subsection{Linear interaction functions}

\nopagebreak
If the functions $f_i$ take the form
\[
f_i(v(x)) = \sum_{j=1}^m S_{ij}v_j(x),
\]
for some constants $S_{ij}$, then (\ref{genprodform1}), (\ref{genprodformout1}) and (\ref{genprodformconserved}) reduce, respectively, to 
\[
\dot x = Sv(x), \quad \dot x = Sv(x) - Q(x) \quad \mbox{and} \quad \dot x = Sv(x) - Q_\theta(x).
\]
These equations are familiar: they are the equations (with or without outflow) for a general CRN. In this context, $x$ becomes the vector of reactant concentrations, $v$ is the vector of reaction rates, and $S$ is the ``stoichiometric matrix''. However the assumption of linearity is also common beyond CRNs, for example in models of regulatory networks (such as that in \cite{elowitz}, discussed in the examples later).

\section{Matrix-theoretic results}

\nopagebreak
\subsection{Sets of real numbers and matrices}
\label{matsets}

\nopagebreak
This paper is concerned with what can be said about the dynamics of (\ref{genprodform1}), (\ref{genprodformout1}) and (\ref{genprodformconserved}), knowing only that $S$ belongs to some matrix-set $\mathcal{S}$, and $V$ belongs to some matrix-set $\mathcal{V}$, possibly related to $\mathcal{S}$. For this reason we start by developing notation and ideas connected with sets of matrices. Definitions and notation which are well known are summarised in Appendix~\ref{app1}. 

Define $\mathbb{R}_{>0} \equiv (0, \infty)$, $\mathbb{R}_{<0} \equiv (-\infty, 0)$, $\mathbb{R}_{\geq 0} \equiv [0, \infty)$ and $\mathbb{R}_{\leq 0} \equiv (-\infty, 0]$. A set of real numbers $\mathcal{R}$ is {\bf signed} if $\mathcal{R} \subset \mathbb{R}_{>0}$ or $\mathcal{R} \subset \mathbb{R}_{<0}$ or $\mathcal{R}  = \{0\}$, and is {\bf weakly signed} if $\mathcal{R} \subset \mathbb{R}_{\geq 0}$ or $\mathcal{R} \subset \mathbb{R}_{\leq 0}$. A set of real numbers which fails to be weakly signed (i.e. intersects both $\mathbb{R}_{>0}$ and $\mathbb{R}_{<0}$) is {\bf unsigned}.

When $\mathcal{M}$ is some set of matrices or real numbers, the statement ``$M$ has property $P$ for all $M \in \mathcal{M}$'' will be abbreviated to ``$\mathcal{M}$ has property $P$''. So, for example, the statement ``$\mathrm{det}(\mathcal{M}) = 0$'' should be read as ``$\mathrm{det}(M) = 0$ for each $M \in \mathcal{M}$''.

{\bf Sums and products of matrix-sets.} Given two sets of matrices, $\mathcal{A}$ and $\mathcal{B}$, we will {\em always} assume that $\mathcal{A}$ and $\mathcal{B}$ are defined as the ranges of matrix-valued functions $\widehat{\mathcal{A}}$, $\widehat{\mathcal{B}}$ over some underlying space $X$, so that $\mathcal{A} = \{\widehat{\mathcal{A}}(x)|x \in X\}$, $\mathcal{B} = \{\widehat{\mathcal{B}}(x)|x \in X\}$. This entails no loss of generality, because even where the matrix-sets are defined purely set-theoretically, or via functions over different spaces, it is an easy matter to redefine the sets as the ranges of matrix-valued functions over some common space. Suppose, for example, $\mathcal{A}$ is defined as the set of all $n \times n$ matrices, and $\mathcal{D}$ as the set of all $n \times n$ positive diagonal matrices. In this case we could define $X = \mathbb{R}^{n \times n} \times \mathbb{R}^n_{>0}$, with $\mathcal{A}$ now the range of a function $\widehat{\mathcal{A}}$ taking the first $n^2$ coordinates of $X$ to elements of $\mathcal{A}$, and $\mathcal{D}$ the range of a function $\widehat{\mathcal{D}}$ taking the final $n$ coordinates to the diagonal elements of $\mathcal{D}$.

Now as a convenient abbreviation we define $\mathcal{A}\mathcal{B} = \{\widehat{\mathcal{A}}(x)\widehat{\mathcal{B}}(x)|x \in X\}$ and $\mathcal{A} + \mathcal{B} = \{\widehat{\mathcal{A}}(x) + \widehat{\mathcal{B}}(x)|x \in X\}$. If, for each $A \in \mathcal{A}, B \in \mathcal{B}$, there exists $x \in X$ such that $A = \widehat{\mathcal{A}}(x), B = \widehat{\mathcal{B}}(x)$, then we will say that $\mathcal{A}$ and $\mathcal{B}$ are {\bf independent}. For independent sets of matrices, it follows that $\mathcal{A}\mathcal{B} = \{AB|A \in \mathcal{A}, B \in \mathcal{B}\}$, and similarly $\mathcal{A} + \mathcal{B} = \{A+B|A \in \mathcal{A}, B \in \mathcal{B}\}$. Thus if two sets are independent then their sum/product is their pointwise sum/product.

As an example, let $x = [x_1, x_2, x_3, x_4]^T \in \mathbb{R}^4$ and for each $x$ define
\[
\widehat{\mathcal{A}}(x) = \left[
\begin{array}{rr}
x_1 & x_2\\x_3 & x_4
\end{array}
\right] \quad \mbox{and} \quad
\widehat{\mathcal{B}}(x) = \left[
\begin{array}{rr}
1 - x_1 & 2 - x_2\\-x_3 & 1- x_4
\end{array}
\right]\,.
\]
Here $\mathcal{A} = \{\widehat{\mathcal{A}}(x)|x \in \mathbb{R}^4\}$ and $\mathcal{B}= \{\widehat{\mathcal{B}}(x)|x \in \mathbb{R}^4\}$ are clearly not independent. In fact, $\mathcal{A} + \mathcal{B}= \{\widehat{\mathcal{A}}(x)+\widehat{\mathcal{B}}(x)|x \in \mathbb{R}^4\}$ consists of the single matrix
\[
\left[
\begin{array}{cc}
1& 2\\0 & 1
\end{array}
\right]\,.
\]
On the other hand, suppose
\[
\widehat{\mathcal{A}}(x) = \left[
\begin{array}{c}
x_1\\x_2
\end{array}
\right] \quad \mbox{and} \quad
\widehat{\mathcal{B}}(x) = \left[
\begin{array}{cc}
x_3-x_1 & x_4-x_2
\end{array}
\right]\,.
\]
In this case $\mathcal{A} = \{\widehat{\mathcal{A}}(x)|x \in \mathbb{R}^4\}$ and $\mathcal{B}= \{\widehat{\mathcal{B}}(x)|x \in \mathbb{R}^4\}$ are independent, and it is easy to confirm that $\mathcal{A}\mathcal{B} = \{AB|A \in \mathcal{A}, B \in \mathcal{B}\}$.

{\bf Notation.} Let $M$ be an $n \times m$ matrix, with $\delta \subset \{1, \ldots, n\}$ and $\gamma \subset \{1,\ldots, m\}$. The following notation will be used:
\begin{itemize}
\item $M(\delta|\gamma)$ is the submatrix of $M$ with rows indexed by $\delta$ and columns indexed by $\gamma$. 
\item $M[\delta|\gamma] \equiv \mathrm{det}(M(\delta|\gamma))$. We write $M[\delta]$ as shorthand for $M[\delta|\delta]$.
\item $M_{\delta\gamma}$ is an $n \times m$ matrix defined by $(M_{\delta\gamma})_{ij} = M_{ij}$ if $i \in \delta$ and $j \in \gamma$ and $(M_{\delta\gamma})_{ij} = 0$ otherwise. Note that all $|\delta| \times |\gamma|$ submatrices of $M_{\delta\gamma}$, apart possibly from $M(\delta|\gamma)$, must contain a row or column of zeros (and hence, if they are square, must be identically singular).
\end{itemize}

A set of real $n \times m$ matrices can be seen as a subset of $\mathbb{R}^{n \times m}$ and thus inherits topological properties such as openness, connectedness, etc. The closure of a matrix-set $\mathcal{M}$ is denoted by $\mathrm{cl}(\mathcal{M})$. $0$ and $I$ denote the zero and identity matrices with dimensions being clear from the context. 

{\bf Entries in matrix-sets and minors of matrix-sets.} Given a matrix-set $\mathcal{A} = \{\widehat{\mathcal{A}}(x)|x \in X\}$, we define $\widehat{\mathcal{A}}_{ij}:X \to \mathbb{R}$ by $\widehat{\mathcal{A}}_{ij}(x) = (\widehat{\mathcal{A}}(x))_{ij}$, and define $\mathcal{A}_{ij} = \{\widehat{\mathcal{A}}_{ij}(x)\,|\,x \in X\}$. A product of entries such as $\mathcal{A}_{ij}\mathcal{A}_{kl}$ means $\{\widehat{\mathcal{A}}_{ij}(x)\widehat{\mathcal{A}}_{kl}(x)\,|\,x \in X\}$, and this notation extends to arbitrary products and sums of products of entries. Given matrix-sets $\mathcal{A} = \{\widehat{\mathcal{A}}(x)|x \in X\}$ and $\mathcal{B} = \{\widehat{\mathcal{B}}(x)|x \in X\}$, $\mathcal{A}[\delta|\gamma]\mathcal{B}[\gamma|\delta]$ should be interpreted as $\{\widehat{\mathcal{A}}(x)[\delta|\gamma]\widehat{\mathcal{B}}(x)[\gamma|\delta]\,|\,x \in X\}$. Expressions involving sums/products of entries/minors always follow these conventions, which considerably simplify notation. 

{\bf Sign-classes.} A set of $n \times m$ matrices $\mathcal{M}$ will be defined as a {\bf sign-class} if given any $M \in \mathcal{M}$, and any $\delta \subset \{1, \ldots, n\}$, $\gamma \subset \{1, \ldots, m\}$, $M_{\delta\gamma} \in \mathrm{cl}(\mathcal{M})$. Since $\gamma$ and $\delta$ may both be empty, $\mathrm{cl}(\mathcal{M})$ contains the zero matrix. 

Sign-classes may have related elements. For example, let $x = [x_1, x_2, x_3, x_4]^T \in \mathbb{R}^4$ and define
\[
\widehat{\mathcal{M}}(x) = \left[
\begin{array}{cc}
x_1x_3 & x_1x_4\\x_2x_3 & x_2x_4
\end{array}
\right], \quad \mathcal{M} = \{\widehat{\mathcal{M}}(x)\,|\,x \in \mathbb{R}^4\}\,.
\]
It is easy to check that $\mathcal{M}$ is a sign-class, even though the entries are not independent, and there is no $M \in \mathrm{cl}(\mathcal{M})$ with $M_{11} = 0$ but $M_{12}, M_{21} \not = 0$.

{\bf Determinant expansions.} Given a set of $n \times n$ matrices $\mathcal{M} = \{\widehat{\mathcal{M}}(x)|x \in X\}$, and a permutation $\alpha$ of the ordered set $[1, \ldots, n]$, define $P(\alpha)$ to be the parity of the permutation, i.e. $P(\alpha) = 1$ for an even permutation and $P(\alpha) = -1$ for an odd permutation. Corresponding to $\alpha$ there is a {\bf term in the determinant expansion} of $\mathcal{M}$,  $T_\alpha$. Define $\widehat{T}_\alpha:X \to \mathbb{R}$, by $\widehat{T}_\alpha(x) = P(\alpha)\prod_{i=1}^n (\widehat{\mathcal{M}}(x))_{i\alpha_i}$, and $T_\alpha = \{\widehat{T}_\alpha(x)|x \in X\}$. The conventions adopted for sums and products of matrix-sets apply to terms in determinant expansions as well, so it is possible for two terms to sum to zero, even though neither consists of a singleton. As an example, let $x = [x_1, x_2]^T \in \mathbb{R}^2$ and define
\[
\widehat{\mathcal{M}}(x) = \left[
\begin{array}{cc}
x_1 & x_2\\x_1 & x_2
\end{array}
\right],  \qquad \mathcal{M} = \{\widehat{\mathcal{M}}(x)|x \in \mathbb{R}^2\}.
\]
Then for $\alpha$ the identity and $\beta = \left(\substack{1\,\,\,2\\ 2\,\,\,1}\right)$, $T_\alpha = \{x_1x_2\,|\,x_1, x_2 \in \mathbb{R}\} = \mathbb{R}$, $T_\beta = \{-x_1x_2\,|\,x_1, x_2 \in \mathbb{R}\} = \mathbb{R}$; however $T_\alpha + T_\beta = \{x_1x_2-x_1x_2\,|\,x_1, x_2 \in \mathbb{R}\} = \{0\}$.

\subsection{$P_0^{(-)}$ systems}
\label{secP0}

\nopagebreak
Our concern is to find conditions which tell us when an interaction network forbids multiple equilibria. This is the case if the system is injective, i.e., the right hand side of the dynamical system is injective. We cannot in general expect (\ref{genprodformconserved}) to be injective: quite generally we may get conserved quantities and hence invariant manifolds foliating the state space (termed stoichiometric classes for CRNs). Often these manifolds are affine subspaces, that is, they are cosets of some linear subspace. Define the conditions:\\
{\bf C0.} (\ref{genprodform1}) have $P_0^{(-)}$ Jacobians.\\
{\bf C1.} (\ref{genprodformout1}) have $P^{(-)}$ Jacobians.\\
{\bf C2.} (\ref{genprodformout1}) have nonsingular Jacobians.\\
{\bf C3.} (\ref{genprodformout1}) is injective on any rectangular domain.\\
{\bf C4.} (\ref{genprodformconserved}) has no more than one nondegenerate equilibrium in the relative interior of any invariant affine subset of the phase space. \\

By basic results on $P_0$ matrices,  {\bf C0} implies {\bf C1} and hence, trivially, {\bf C2}. Below, in Corollary~\ref{P0systems}, we show that {\bf C2}  implies {\bf C0}, and so {\bf C0, C1 and C2 are equivalent}. Techniques such as those employed in \cite{cinquin,CraciunHeltonWilliams,helton} to prove absence of multiple equilibria using degree theory require only a nonsingular Jacobian rather than $P$ or $P^{(-)}$ Jacobians. The equivalence of {\bf C0}, {\bf C1} and {\bf C2} shows that this requirement is only apparently weaker for systems which may have arbitrary outflows. It is well known that if a function $f$ has a $P$ matrix (or $P^{(-)}$ matrix) Jacobian on some rectangular domain in $\mathbb{R}^n$, then $f$ is injective on this domain (Theorem~4 and subsequent remarks in \cite{gale}), a fact which was applied in \cite{soule,banajiSIAM}. Thus {\bf C1 implies C3}. Finally, {\bf C3 implies C4}. The meaning and proof of this fact are presented in Appendix~\ref{apppersistence}. Thus {\bf C0} implies all the other statements. Defining a ``$P_0^{(-)}$ system'' to be a system with $P_0^{(-)}$ Jacobian, we have:
\begin{quote}
If (\ref{genprodform1}) is a $P_0^(-)$ system, then (\ref{genprodformout1}) is injective, and (\ref{genprodformconserved}) admits no multiple nondegenerate equilibria on the relative interior of any invariant affine subset of the phase space. Note that (\ref{genprodformout1}) was obtained from (\ref{genprodform1}) by assuming that each $x_i$ is subject to outflow or degradation, and (\ref{genprodformconserved}) was obtained from $(N_0)$ by assuming that some $x_i$'s may be subject to outflow or degradation.
\end{quote}

{\bf Remark.} Note that in the case of some chemical reaction systems, inectivity of (\ref{genprodformout1}) implies injectivity of (\ref{genprodformconserved}) for any $\theta$ \cite{craciun4}.

\begin{lem}
\label{P0sys}
Consider any $n \times n$ matrix $A$ and let $\mathcal{D}$ be the set of $n \times n$ positive diagonal matrices. If $A$ fails to be a $P_0^{(-)}$ matrix, then $A- \mathcal{D} = \{A - D\,|\,D \in \mathcal{D}\}$ contains matrices with determinants of all signs. 
\end{lem}
\begin{proof}
Given any $D \in \mathcal{D}$, write $A-D = -D(I + D^{-1}(-A))$, so that $\mathrm{sign}(\mathrm{det}(A-D)) = -\mathrm{sign}(\mathrm{det}(I + D^{-1}(-A)))$. Expanding, we have
\[
\mathrm{det}(I + D^{-1}(-A)) = 1 + \sum_{\delta \subset\{1, \ldots, n\}} (D^{-1}(-A))[\delta]
\]
where the sum is taken over all nonempty subsets of $\{1, \ldots, n\}$. Using the Cauchy-Binet formula \cite{gantmacher}, and noting that the only nonzero minors of a diagonal matrix are the principal minors, gives $(D^{-1}(-A))[\delta] = D^{-1}[\delta](-A)[\delta]$. 

First, for any $\epsilon > 0$ define $D(\epsilon) \in \mathcal{D}$ by $D(\epsilon)_{ii} = 1/\epsilon$. Then $[D(\epsilon)^{-1}]_{ii} = \epsilon$, and $(D(\epsilon)^{-1})[\delta] = \epsilon^{|\delta|}$. So $\mathrm{det}(I + D(\epsilon)^{-1}(-A)) = 1 + \sum_{\delta \subset\{1, \ldots, n\}} \epsilon^{|\delta|}(-A)[\delta]$. For small enough $\epsilon>0$, clearly $\mathrm{det}(I + D(\epsilon)^{-1}(-A)) > 0$. Thus there exists a positive diagonal matrix $D_1$ such that $\mathrm{sign}(\mathrm{det}(A-D_1)) = -1$. 

Suppose that $A$ is not a $P_0^{(-)}$ matrix. This means that there is some nonempty $\delta_0 \subset \{1, \ldots, n\}$ such that $(-A)[\delta_0] < 0$. Now for any $\epsilon > 0$ define $D(\epsilon) \in \mathcal{D}$ by $D(\epsilon)_{ii} = \epsilon$ if $i \in \delta_0$ and $D(\epsilon)_{ii} = 1/\epsilon$ if $i \not \in \delta_0$. In this case $(D(\epsilon)^{-1})[\delta] = \epsilon^{|\delta\backslash \delta_0| - |\delta\cap \delta_0|}$. From this we get that 
\begin{eqnarray}
\mathrm{det}(I + D(\epsilon)^{-1}(-A)) & = & 1 + \sum_{\delta \subset\{1, \ldots, n\}} \epsilon^{|\delta\backslash \delta_0| - |\delta\cap \delta_0|}(-A)[\delta] \nonumber \\
& = & \epsilon^{-|\delta_0|}(-A)[\delta_0] + \mbox{ higher order terms in $\epsilon$}.\nonumber
\end{eqnarray}
Clearly, for small enough $\epsilon > 0$, $\mathrm{det}(I + D(\epsilon)^{-1}(-A)) < 0$. Thus there exists a positive diagonal matrix $D_2$ such that $\mathrm{sign}(\mathrm{det}(A-D_2)) = 1$. By continuity of the determinant, there exists $\lambda \in (0,1)$ such that $\mathrm{det}(A-\lambda D_1 - (1-\lambda)D_2) = 0$. As $\lambda D_1 + (1-\lambda)D_2 \in \mathcal{D}$, the result is proved. 
\end{proof}

\begin{cor}
\label{P0systems}
Consider any set of  $n \times n$ matrices $\mathcal{A}$ and let $\mathcal{D}$ be the set of $n \times n$ positive diagonal matrices. If any $A \in \mathcal{A}$ fails to be a $P_0^{(-)}$ matrix, then $\mathcal{A}- \mathcal{D} = \{A - D\,|\,A \in \mathcal{A}, D \in \mathcal{D}\}$ contains matrices with determinants of all signs. 
\end{cor}
\begin{proof}
This is an immediate corollary of Lemma~\ref{P0sys}.
\end{proof}

{\bf Remark.} Together with the basic fact that a $P_0$ matrix plus a positive diagonal matrix is a $P$ matrix, Lemma~\ref{P0sys} implies the following characterisation of $P_0$ matrices: a matrix $A$ is a $P_0$ matrix iff $\mathrm{det}(A + D) > 0$ for every positive diagonal matrix $D$.

\subsection{Identifying matrix products as $P_0^{(-)}$ matrices}

\nopagebreak
Having seen that whether a set of matrices are $P_0^{(-)}$ matrices is central to questions concerning multiple equilibria, the next stage is to present a general rule for deciding when the product of two matrix-sets $\mathcal{S}$ and $\mathcal{V}$ (independent or otherwise), consists entirely of $P_0^{(-)}$ matrices. A very general sufficient condition is provided in Lemma~\ref{genlem}. This condition is also proved to be necessary in Lemma~\ref{genlemconv}, provided $\mathcal{S}$ and $\mathcal{V}$ are independent and at least one of $\mathcal{S}$ or $\mathcal{V}$ is a sign-class.

\begin{lem}
\label{genlem}
Consider any set of $n \times m$ matrices $\mathcal{S}$ and any set of $m \times n$ matrices $\mathcal{V}$. Assume that for any $\delta \subset \{1, \ldots, n\}$, $\gamma\subset \{1, \ldots, m\}$ satisfying $|\gamma| = |\delta|$, $\mathcal{S}[\delta|\gamma]\mathcal{V}[\gamma|\delta] \subset (-1)^{|\delta|}\mathbb{R}_{\geq 0}$. Then $\mathcal{S}\mathcal{V}$ consists of $P_0^{(-)}$ matrices.
\end{lem}
\begin{proof}
The result follows from the Cauchy-Binet formula. In the usual way, let $\mathcal{S} = \{\widehat{\mathcal{S}}(x)\,|\,x \in X\}$, and $\mathcal{V} = \{\widehat{\mathcal{V}}(x)\,|\,x \in X\}$. Given some $x \in X$ and some $S = \widehat{\mathcal{S}}(x)$ and $V = \widehat{\mathcal{V}}(x)$ we have
\[
(SV)[\delta] = \sum_{\substack{\gamma \subset \{1, \ldots, m\}\\ |\gamma| = |\delta|}}S[\delta|\gamma]V[\gamma|\delta]\,.
\]
Supposing the conditions of the lemma are fulfilled, then $S[\delta|\gamma]V[\gamma|\delta] \in (-1)^{|\delta|}\mathbb{R}_{\geq 0}$, and thus $(SV)[\delta] \in (-1)^{|\delta|}\mathbb{R}_{\geq 0}$. This proves that $SV$ is a $P_0^{(-)}$ matrix. 
\end{proof}

{\bf Remark.} If $\mathcal{S}$ and $\mathcal{V}$ are independent, the condition $\mathcal{S}[\delta|\gamma]\mathcal{V}[\gamma|\delta] \subset (-1)^{|\delta|}\mathbb{R}_{\geq 0}$ means that either, i) $\mathcal{S}[\delta|\gamma] = 0$, or ii) $\mathcal{V}[\gamma|\delta] = 0$, or iii) $\mathcal{S}[\delta|\gamma]$ and $\mathcal{V}[\gamma|\delta]$ are both weakly signed with either $\mathcal{S}[\delta|\gamma],(-\mathcal{V})[\gamma|\delta] \subset \mathbb{R}_{\geq 0}$ or $\mathcal{S}[\delta|\gamma],(-\mathcal{V})[\gamma|\delta] \subset \mathbb{R}_{\leq 0}$.

Depending on $\mathcal{S}$ and $\mathcal{V}$, the conditions in Lemma~\ref{genlem} may be necessary as well as sufficient to guarantee that $\mathcal{S}\mathcal{V}$ consists of $P_0^{(-)}$ matrices. In particular the following case where either $\mathcal{S}$ or $\mathcal{V}$ (or both) are sign-classes of matrices often arises:

\begin{lem}
\label{genlemconv}
Consider any set of $n \times m$ matrices $\mathcal{S}$ and any set of $m \times n$ matrices $\mathcal{V}$. Assume that $\mathcal{S}$ and $\mathcal{V}$ are independent, one of $\mathcal{S}$ or $\mathcal{V}$ is a sign-class, and $\mathcal{S}\mathcal{V}$ consists of $P_0^{(-)}$ matrices. Then $\mathcal{S}[\delta|\gamma]\mathcal{V}[\gamma|\delta] \subset (-1)^{|\delta|}\mathbb{R}_{\geq 0}$ for each $\delta \subset \{1, \ldots, n\}$, $\gamma\subset \{1, \ldots, m\}$ satisfying $|\gamma| = |\delta|$.
\end{lem}
\begin{proof}
Suppose there are some sets $\delta$ and $\gamma$ such that $\mathcal{S}[\delta|\gamma]\mathcal{V}[\gamma|\delta] \not \subset (-1)^{|\delta|}\mathbb{R}_{\geq 0}$. This means that there exist $S \in \mathcal{S}$ and $V \in \mathcal{V}$ such that $S[\delta|\gamma]V[\gamma|\delta] \in (-1)^{|\delta|+1}\mathbb{R}_{> 0}$. For definiteness assume that $\mathcal{V}$ is a sign-class of matrices, so that $V_{\gamma\delta} \in \mathrm{cl}(\mathcal{V})$, and so $(SV_{\gamma\delta})[\delta] = S[\delta|\gamma]V[\gamma|\delta] \in (-1)^{|\delta|+1}\mathbb{R}_{> 0}$. Thus $SV_{\gamma\delta}$ fails to be a $P_0^{(-)}$ matrix. Since the class of $P_0^{(-)}$ matrices is closed, its complement is open. Thus any matrix sufficiently near to $SV_{\gamma\delta}$ fails to be a $P_0^{(-)}$ matrix. Since, by independence of $\mathcal{S}$ and $\mathcal{V}$, $SV_{\gamma\delta} \in \mathrm{cl}(\mathcal{S}\mathcal{V})$, there are matrices in $\mathcal{S}\mathcal{V}$ which fail to be $P_0^{(-)}$. The argument works equally well if $\mathcal{S}$ is a sign-class. 
\end{proof}

There are a couple of immediate observations to be made from the above lemmas. Lemma~\ref{genlemconv} places the following basic restriction on systems with unsigned entries and $P_0^{(-)}$ Jacobians:
\begin{cor}
\label{onereac}
Consider a system defined by matrix-sets $\mathcal{S}$ and $\mathcal{V}$ fulfilling the conditions in Lemma~\ref{genlemconv}. If for some $i\in \{1, \ldots, n\}$, $j \in \{1, \ldots, m\}$,  $\mathcal{V}_{ji}$ is unsigned, then $\mathcal{S}_{ij} = 0$. Similarly if $\mathcal{S}_{ij}$ is unsigned, then $\mathcal{V}_{ji} = 0$. 
\end{cor}
\begin{proof}
Setting $\delta = \{i\}$ and $\gamma = \{j\}$ in Lemma~\ref{genlemconv} gives us the result. 
\end{proof}

Loosely speaking, this corollary tells us that in a $P_0^{(-)}$ system, if a quantity has influence of unknown sign on an interaction, then it should not itself be affected by the interaction. In particular, in CRNs, if a substrate occurs on both sides of a reaction with unknown influence on the rate, then the substrate must occur with the same stoichiometry on both sides. As an example of how more complicated ``forbidden scenarios'' can also be formulated, we have the following result:
\begin{cor}
\label{tworeac}
Consider a system defined by matrix-sets $\mathcal{S}$ and $\mathcal{V}$ fulfilling the conditions in Lemma~\ref{genlemconv}, with $\mathcal{V}$ a sign-class. Suppose for some indices $i, j, k$, $\mathcal{S}_{ji} \not = 0$ and $\mathcal{V}_{kj}$ is unsigned. Then there is no index $l$ such that $\mathcal{S}_{lk} \not = 0$ and $\mathcal{V}_{il} \not = 0$.  
\end{cor}

\begin{proof}
Assume the contrary and consider the $2 \times 2$ submatrices
\[
\mathcal{S}(\{j, l\}|\{i, k\}) = \left[
\begin{array}{cc}
\mathcal{S}_{ji} & \mathcal{S}_{jk}\\
\mathcal{S}_{li} & \mathcal{S}_{lk}
\end{array}
\right] \quad \mbox{and} \quad
\mathcal{V}(\{i, k\}|\{j, l\}) = \left[
\begin{array}{cc}
\mathcal{V}_{ij} & \mathcal{V}_{il}\\
\mathcal{V}_{kj} & \mathcal{V}_{kl}
\end{array}
\right]
\]
From Corollary~\ref{onereac}, $\mathcal{S}_{jk}$ must be zero, and by assumption $\mathcal{S}_{ji}, \mathcal{S}_{lk} \not = 0$, and so $\mathcal{S}[\{j, l\}|\{i, k\}] \not = 0$. Further, $\mathcal{V}_{il} \not = 0$, and $\mathcal{V}_{kj}$ is unsigned, so, since $\mathcal{V}$ is a sign-class, $\mathcal{V}[\{i, k\}|\{j, l\}]$ is unsigned. By Lemma~\ref{genlemconv}, $\mathcal{S}\mathcal{V}$ are not all $P_0^{(-)}$ matrices. 
 \end{proof}

The graphical presentation of this result in Section~\ref{tworeacexample} provides some intution as to its meaning.

{\bf Remark on previous matrix-theoretic results.} The special case described in \cite{banajiSIAM,banajicraciun} of chemical reactions with no substrate on both sides of any reaction (termed NAC systems in the first reference and N1C system in the second) corresponds to $\mathcal{S}$ consisting of a single matrix, i.e. $\mathcal{S} = \{S\}$, and $\mathcal{V} = \mathcal{Q}_0(-S^T)$ (see Appendix~\ref{app1} for a definition). In that case the condition $\mathcal{S}[\delta|\gamma]\mathcal{V}[\gamma|\delta] \subset (-1)^{|\delta|}\mathbb{R}_{\geq 0}$ in Lemma~\ref{genlem} implies that all square submatrices of $S$ are either singular or sign nonsingular (see Appendix~\ref{app1}). Since $\{S\}$ and $\mathcal{Q}_0(-S^T)$ are trivially independent, and $\mathcal{Q}_0(-S^T)$ is a sign-class of matrices, this condition on $S$ is both necessary and sufficient for all matrices $S\mathcal{V}$ to be $P_0^{(-)}$ matrices. In \cite{banajiSIAM} there was also some discussion of a model where $\mathcal{S}$ was not a single matrix. 

Lemmas~\ref{genlem} and \ref{genlemconv} directly give a number of easy generalisations of matrix-theoretic results in \cite{banajiSIAM} which will be developed in future work. Here we concentrate on the graph-theoretic corollaries.

\section{SR graphs and DSR graphs}

\nopagebreak
\subsection{SR graphs}

\nopagebreak
{\bf SR graphs for single (rectangular) matrices.} SR graphs (or species-reaction graphs) are bipartite multigraphs with two vertex-sets termed S-vertices and R-vertices, and signed, labelled edges. Although originally defined for CRNs \cite{craciun1}, they can also directly be associated with matrices as in \cite{banajicraciun}: given any rectangular matrix $M$, we can construct the associated SR graph $G_{M}$ where an edge exists between S-vertex $i$ and R-vertex $j$ iff $M_{ij} \not = 0$. The edge takes the sign of $M_{ij}$ and the edge-label $|M_{ij}|$. When the matrix $M$ is the stoichiometric matrix of a CRN with SR graph $G$, then $G_{M} = G$ for any reasonable kinetics iff no substrates occur on both sides of any reaction in the system. This special case is treated in detail in \cite{banajicraciun}.

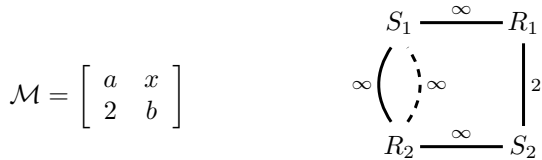
\begin{figure}[h]
\begin{minipage}{0.4\textwidth}
\[
\mathcal{M} = \left[\begin{array}{cc}a & x\\2&b\end{array}\right]
\]

\end{minipage}
\hfill
\begin{minipage}{0.55\textwidth}
\begin{tikzpicture}[domain=0:4,scale=0.55]

\node at (1,4) {$S_1$};
\node at (4,4) {$R_1$};
\node at (4,1) {$S_2$};
\node at (1,1) {$R_2$};

\draw[-, line width=0.04cm] (0.8,1.6) .. controls (0.4,2.2) and (0.4,2.8) .. (0.8,3.4);
\draw[-, dashed, line width=0.04cm] (1.2,1.6) .. controls (1.6,2.2) and (1.6,2.8) .. (1.2,3.4);

\draw[-, line width=0.04cm] (1.5,4.0) -- (3.5, 4.0);
\draw[-, line width=0.04cm] (1.5,1.0) -- (3.5, 1.0);
\draw[-, line width=0.04cm] (4,1.5) -- (4, 3.5);

\node at (2.5, 4.3) {$\scriptstyle{\infty}$};
\node at (2.5, 1.3) {$\scriptstyle{\infty}$};

\node at (4.3, 2.5) {$\scriptstyle{2}$};

\node at (1.9, 2.5) {$\scriptstyle{\infty}$};
\node at (0.1, 2.5) {$\scriptstyle{\infty}$};

\end{tikzpicture}

\end{minipage}
\caption{\label{SRexample} Assume that $a, b$ can take all values in $[0, 1]$ and $x$ can take all values in $[-1, 1]$. The SR graph associated with $\mathcal{M}$ contains five edges. Positive edges are represented with bold lines, negative edges are represented with dashed lines. These conventions will be followed throughout.}
\end{figure}

{\bf SR graphs for sets of (rectangular) matrices.} Given any matrix-set $\mathcal{M}$, we can construct an SR graph $G_{\mathcal{M}}$ which is, roughly speaking, the amalgamation of the SR graphs $G_M$ associated with each $M \in \mathcal{M}$. If $\mathcal{M}_{ij} = \{0\}$, then there is no edge between S-vertex $i$ and R-vertex $j$; if $\mathcal{M}_{ij} \not = \{0\}$ and $\mathcal{M}_{ij} \subset [0, \infty)$, then there is an edge between S-vertex $i$ and R-vertex $j$ with sign $+1$; if $\mathcal{M}_{ij} \not = \{0\}$ and $\mathcal{M}_{ij}  \subset (-\infty, 0]$, then there is an edge between S-vertex $i$ and R-vertex $j$ with sign $-1$; finally if $\mathcal{M}_{ij}$ is unsigned then we introduce a pair of oppositely signed edges between S-vertex $i$ and R-vertex $j$. If $\mathcal{M}_{ij} = \{k\}$ where $k \not = 0$, then we give the unique edge between S-vertex $i$ and R-vertex $j$ an edge-label of $|k|$, otherwise we give all edges between these vertices edge-labels of $\infty$. An example of a matrix-set and the associated SR graph is shown in Figure~\ref{SRexample}.

{\bf Cycles in SR graphs.} Cycles in SR graphs are minimal undirected paths from some vertex to itself. As seen in Figure~\ref{SRexample}, an SR graph may have a pair of oppositely signed edges between S-vertex $i$ and R-vertex $j$ forming a cycle of length $2$. Such cycles will be termed {\bf short cycles}. Cycles of length greater than $2$ are then {\bf long cycles}. Since all edges in an SR graph are signed, all cycles in an SR graph have a sign, defined as the product of signs of edges in the cycle. A cycle $C$ in an SR graph also has a {\bf parity} $P(C)$: it is either an {\bf o-cycle} or an {\bf e-cycle} according to whether
\[
P(C) = (-1)^{|C|/2}\mathrm{sign}(C)
\]
is negative or positive. Note that short cycles in an SR graph are always e-cycles. Given an edge $e$, define $\mathrm{val}(e)$ to be the edge-label of $e$. When $C$ is a cycle containing edges $e_1, e_2, \ldots, e_{2r}$ such that $e_i$ and $e_{(i \mod 2r) +1}$ are adjacent for each $i =1,\ldots,2r$, we can define its stoichiometry as follows: if any edge in $C$ has edge-label $\infty$, then we set $\mathrm{stoich}(C) = \infty$; otherwise
\[
\mathrm{stoich}(C) = \left|\prod_{i = 1}^{r}\mathrm{val}(e_{2i-1}) - \prod_{i = 1}^{r}\mathrm{val}(e_{2i})\right|\,.
\]
Note that this definition is independent of the starting point chosen on the cycle. A cycle with $\mathrm{stoich}(C) = 0$ is termed an {\bf s-cycle}. An e-cycle which is also an s-cycle is an {\bf es-cycle}. A {\bf disconnecting partition} of a cycle $C$ is the (unique) partition of the cycle into two edge-sets, such that no two edges in either set share a vertex.

{\bf S-to-R intersection in SR graphs}. The intersection of two cycles in an SR graph can be divided into a set of vertex-disjoint paths. We say that two cycles have S-to-R intersection if each component of their intersection is an S-to-R path, i.e. a path between an S-vertex and an R-vertex.

\subsection{DSR graphs}

\nopagebreak
Where SR graphs are associated with matrix-sets, DSR graphs are associated with {\em pairs} of matrix-sets. Given any two sets of $n \times m$ matrices, $\mathcal{A}$ and $\mathcal{B}$, the DSR graph $\mathcal{G}_{\mathcal{A},\mathcal{B}}$ is a signed, labelled, bipartite multigraph with edges which may be directed or undirected. Here the definition is presented, while a more intuitive discussion is presented in Appendix~\ref{appDSR}. 

First, we take the SR graph $G_{\mathcal{A}}$ and create a directed version $\overleftarrow{G}_{\mathcal{A}}$ by insisting that all edges are directed from R- to S-vertices. Similarly, we define $\overrightarrow{G}_{\mathcal{B}}$, a directed version of $G_{\mathcal{B}}$ with all edges pointing from S- to R-vertices, and all edge-labels set to be $\infty$. 
\begin{enumerate}
\item Since S-vertices in $\overleftarrow{G}_{\mathcal{A}}$ and $\overrightarrow{G}_{\mathcal{B}}$ can be identified, and similarly for R-vertices, $\overleftarrow{G}_{\mathcal{A}}$ and $\overrightarrow{G}_{\mathcal{B}}$ can be amalgamated into a single signed, labelled, directed multigraph $\tilde{G}$ with all edges from both $\overleftarrow{G}_{\mathcal{A}}$ and $\overrightarrow{G}_{\mathcal{B}}$. (Any S-vertex and R-vertex in $\tilde{G}$ may be connected by $0, 1, 2, 3,$ or $4$ edges.)
\item If two vertices in $\tilde{G}$ are connected by a pair of edges with the same sign but opposite orientation, then we replace these with a single undirected edge with sign and edge-label imported from $\overleftarrow{G}_{\mathcal{A}}$. The resulting graph $G_{\mathcal{A},\mathcal{B}}$ is now termed the DSR graph. Any S-vertex and R-vertex in $G_{\mathcal{A},\mathcal{B}}$ may be connected by $0, 1$ or $2$ edges, some of which may be directed.
\end{enumerate}

Every edge-set in $G_{\mathcal{A}}$ has a corresponding edge-set in $G_{\mathcal{A},\mathcal{B}}$ with the same edge-signs and edge-labels. Thus ignoring direction of edges, $G_{\mathcal{A}}$ is a true subgraph of $G_{\mathcal{A},\mathcal{B}}$. Similarly, ignoring direction and labels on edges, $G_{\mathcal{B}}$ is a true subgraph of $G_{\mathcal{A},\mathcal{B}}$. Given matrix-sets $\mathcal{A}$ and $\mathcal{B}$, $G_{\mathcal{A},\mathcal{B}}$ is not in general the same as $G_{\mathcal{B},\mathcal{A}}$, although if we ignore edge-labels, reversing the orientation of all directed edges takes $G_{\mathcal{A},\mathcal{B}}$ to $G_{\mathcal{B},\mathcal{A}}$. 

{\bf Subgraphs of DSR graphs}. If a DSR graph $G$ is associated with a pair $(\mathcal{A}, \mathcal{B})$ of sets of $n \times m$ matrices, then given $\delta \subset \{1, \ldots, n\}$ and $\gamma \subset \{1, \ldots, m\}$, $G(\delta|\gamma)$ will refer to the subgraph of $G$ associated with the pair $(\mathcal{A}(\delta|\gamma), \mathcal{B}(\delta|\gamma))$. DSR subgraphs with an equal number of S- and R-vertices will be referred to as {\bf square} DSR graphs. 

{\bf Directed paths in a DSR graph}. An edge in a DSR graph has {\bf R-to-S direction} if it is undirected or directed from an R-vertex to an S-vertex. It has {\bf S-to-R direction} if it is undirected or directed from an S-vertex to an R-vertex. So undirected edges have both R-to-S direction and S-to-R direction. A subset of the edge-set of a DSR graph has R-to-S (S-to-R) direction if all edges in the set have R-to-S (S-to-R) direction. Two edge-sets in a DSR graph are {\bf oppositely directed} if all edges in one have R-to-S direction and all edges in the other have S-to-R direction (one or both sets may contain undirected edges). A directed path in a DSR graph consists of alternating S-to-R and R-to-S edges.

{\bf Cycles in DSR graphs.} A cycle $C$ in a DSR graph is a minimal directed path from some vertex to itself. Of course, some edges in a cycle may be undirected. Cycles are always either e-cycles or o-cycles, since all edges are signed. They may also be s-cycles. The following two remarks follow immediately from the definition: all short cycles in a DSR graph are e-cycles but not es-cycles because at least one edge has edge-label $\infty$; any cycle has a disconnecting partition consisting of an S-to-R edge-set and an R-to-S edge-set of equal size. 

{\bf Formal cycles.} Given a DSR graph $G$ we can construct an undirected version of it, $G^{'}$ by making all directed edges in $G$ undirected in $G^{'}$. Cycles in $G^{'}$ will be termed formal cycles in $G$ and may or may not be cycles. Sometimes we will refer to a cycle in $G$ as a ``genuine'' cycle to emphasise that the cycle is directed. 

{\bf Orientation of cycles.} If a (genuine) cycle $C$ in a DSR graph contains only undirected edges, then it has two natural orientations. On the other hand, if $C$ contains some edge which fails to have both S-to-R and R-to-S direction, then $C$ has one natural orientation. Thus we can always choose either one or two orientations for any genuine cycle. If we choose an orientation for a cycle $C$, then each edge (including undirected edges) in $C$ inherits an orientation, which we can call that edge's ``$C$-orientation''. We say that two cycles $C$ and $D$ have {\bf compatible orientation} if we can choose an orientation for $C$ and an orientation for $D$ such that each edge in their intersection has the same $C$-orientation and $D$-orientation. Note that two cycles with empty edge-intersection trivially have compatible orientation, and that two undirected cycles may be incompatibly oriented. 

{\bf S-to-R intersection between cycles}. Two cycles in a DSR graph are said to have S-to-R intersection if they have nonempty intersection, compatible orientation, and moreover each component of their intersection has odd length, i.e. it is either an S-to-R path or an R-to-S path. (Note that there is no implied direction in the term ``S-to-R intersection'').

\subsection{The main results}

\nopagebreak
We have now developed sufficient terminology to state the main results of this paper. Define:

{\bf Condition ($*$)}: All e-cycles in a DSR graph are s-cycles, and no two e-cycles have S-to-R intersection. 

\begin{thm}
\label{mainthm}
Consider a set of $n \times m$ matrices $\mathcal{S}$, a set of $m \times n$ matrices $\mathcal{V}$, and the corresponding DSR graph $G \equiv G_{\mathcal{S}, -\mathcal{V}^T}$. If $G$ satisfies Condition~($*$), then $\mathcal{S}$ and $\mathcal{V}$ satisfy the conditions in Lemma~\ref{genlem}.
\end{thm}

The proof of this theorem will be presented after considerable machinery has been built up. The following corollary, which is the main result in this paper, tells us what Theorem~\ref{mainthm} implies.

\begin{cor}
\label{maincor}
Consider an interaction network (\ref{genprodform1}) defined on a rectangular domain $X \subset \mathbb{R}^n$, with DSR graph $G$ satisfying Condition~($*$). Then the Jacobians of (\ref{genprodform1}) are all $P_0^{(-)}$ matrices. Associated systems of the form (\ref{genprodformout1}) have Jacobians which are all $P^{(-)}$ matrices, and hence cannot have multiple equilibria in $X$. Associated systems of the form (\ref{genprodformconserved}), including (\ref{genprodform1}) itself, cannot have multiple nondegenerate equilibria in the relative interior of any invariant affine subset of $X$. 
\end{cor}
\begin{proof}
These claims follow immediately from Theorem~\ref{mainthm} coupled with the discussion in Section~\ref{secP0}. 
\end{proof}

{\bf Remark.} Note that there is an alternative ``pointwise'' phrasing of Corollary~\ref{maincor}. Rather than associating a single DSR graph $G$ with the interaction network (\ref{genprodform1}), we can associate a DSR graph $G(x)$ with each point $x \in X$. Provided each $G(x)$ satisfies Condition~($*$), then the Jacobians of (\ref{genprodform1}) are all $P_0^{(-)}$ matrices, and all the other conclusions follow. It can occur that $G$ fails Condition~($*$), but at each $x \in X$, the DSR graph $G(x)$ satisfies Condition~($*$). Although generating DSR graphs for each $x \in X$ gives stronger results, application of the result is simplified if a single graph is generated. 

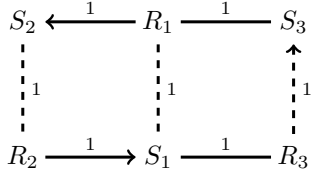
\begin{figure}[h]
\begin{minipage}{\textwidth}
\begin{center}
\begin{tikzpicture}[domain=0:4,scale=0.6]

\node at (1,4) {$S_2$};
\node at (4,4) {$R_1$};
\node at (7,4) {$S_3$};
\node at (1,1) {$R_2$};
\node at (4,1) {$S_1$};
\node at (7,1) {$R_3$};

\draw[<-, line width=0.04cm] (1.5, 4) -- (3.5,4);
\draw[->, line width=0.04cm] (1.5, 1) -- (3.5,1);
\draw[-, line width=0.04cm] (4.5, 4) -- (6.5,4);
\draw[-, line width=0.04cm] (4.5, 1) -- (6.5,1);

\draw[-, dashed, line width=0.04cm] (1, 3.5) -- (1,1.5);
\draw[-, dashed, line width=0.04cm] (4, 3.5) -- (4,1.5);
\draw[<-, dashed, line width=0.04cm] (7, 3.5) -- (7,1.5);

\node at (2.5, 4.3) {$\scriptstyle{1}$};
\node at (2.5, 1.3) {$\scriptstyle{1}$};

\node at (1.3, 2.5) {$\scriptstyle{1}$};
\node at (4.3, 2.5) {$\scriptstyle{1}$};

\node at (5.5, 4.3) {$\scriptstyle{1}$};
\node at (5.5, 1.3) {$\scriptstyle{1}$};

\node at (7.3, 2.5) {$\scriptstyle{1}$};
\end{tikzpicture}
\end{center}

\end{minipage}
\caption{\label{Cstar} The DSR graph for the pair $A, \mathcal{B}$ in (\ref{orientationeq}) above. The graph satisfies Condition~($*$).}
\end{figure}

{\bf Example.} As an illustration of Condition~($*$) consider the following pair of matrix-sets (the first consists of a single matrix):
\begin{equation}
\label{orientationeq}
A = \left[\begin{array}{rrr}-1 & 1 & 1\\1 & -1 & 0\\1 & 0 & -1\end{array}\right],\quad \mathcal{B} = \left[\begin{array}{rrr}-a & 0 & b\\0 & -c & \,\,\,0\\d & 0 & 0\end{array}\right]\,\,.
\end{equation}
Assume that $a, b, c, d > 0$. The associated DSR graph $G_{A, \mathcal{B}}$ is shown in Figure~\ref{Cstar}. Although the es-cycles $S_1\!-\!R_1\!-\!S_2\!-\!R_2$ and $R_1\!-\!S_1\!-\!R_3\!-\!S_3$ appear to have S-to-R intersection, the are not compatibly oriented, and in fact Condition~($*$) holds. 

\subsection{Determinant expansions and structures in SR and DSR graphs}

\nopagebreak
Here we develop some methodology relating terms in determinant expansions and objects in SR and DSR graphs. The most important notions are {\bf signed subterms} in determinant expansions which are in one-to-one correspondence with {\bf signed term subgraphs} in SR or DSR graphs.

Below, $\mathcal{A}$ and $\mathcal{B}$ are sets of $n \times m$ matrices, $G_{\mathcal{A}}$ and $G_{\mathcal{B}}$ are the associated SR graphs, and $G \equiv G_{\mathcal{A},\mathcal{B}}$ is the associated DSR graph. $\gamma =[\gamma_1, \gamma_2, \ldots, \gamma_k] \subset \{1, \ldots, m\}$ and $\delta = [ \delta_1, \delta_2, \ldots, \delta_k] \subset \{1, \ldots, n\}$ are fixed, ordered sets of equal size.

{\bf Signed subentries in matrix-sets}. Let $\mathcal{A} = \{\widehat{\mathcal{A}}(x)|x \in X\}$, with $\widehat{\mathcal{A}}_{ij}:X \to \mathbb{R}$ defined by 
$\widehat{\mathcal{A}}_{ij}(x) = (\widehat{\mathcal{A}}(x))_{ij}$, and $\mathcal{A}_{ij} = \{\widehat{\mathcal{A}}_{ij}(x)\,|\,x \in X\}$. Since $\widehat{\mathcal{A}}_{ij}$ is a scalar function on $X$, we can split it into $\widehat{\mathcal{A}}_{ij} = \widehat{\mathcal{A}}_{ij}^{+} + \widehat{\mathcal{A}}_{ij}^{-}$ where $\widehat{\mathcal{A}}_{ij}^{+}(x)  = \max\{ \widehat{\mathcal{A}}_{ij}(x), 0\}$ and $\widehat{\mathcal{A}}_{ij}^{-}(x)  = \min\{ \widehat{\mathcal{A}}_{ij}(x), 0\}$. We call $\mathcal{A}_{ij}^{+} = \{\widehat{\mathcal{A}}_{ij}^{+}(x)\,|\,x \in X\}$ and $\mathcal{A}_{ij}^{-} = \{\widehat{\mathcal{A}}_{ij}^{-}(x)\,|\,x \in X\}$ signed subentries in $\mathcal{A}$ with the usual rules for addition and multiplication. Corresponding to each of $\mathcal{A}_{ij}^{+}$ and $\mathcal{A}_{ij}^{-}$ is an edge in $G_{\mathcal{A}}$. The same methodology can be applied to unsigned entries in $\mathcal{B}$.

{\bf Terms in determinant expansions:} Given a permutation $\alpha$ of $\gamma$, define:
\[
\overleftarrow{T}_\alpha = P(\alpha)\prod_{i=1}^{|\delta|}\mathcal{A}_{\delta_i\alpha_i},\qquad \overrightarrow{T}_\alpha = P(\alpha)\prod_{i=1}^{|\delta|}\mathcal{B}_{\delta_i\alpha_i}\,.
\]
For an interpretation of these expressions, the reader is referred back to Section~\ref{matsets}. $\overleftarrow{T}_\alpha$ is a term in the expansion of $\mathcal{A}[\delta|\gamma]$, while $\overrightarrow{T}_\alpha$ is the corresponding term in the expansion of $\mathcal{B}[\delta|\gamma]$. From now on, where an ordered pair of matrix-sets $(\mathcal{A}, \mathcal{B})$ is concerned, objects with left-pointing arrows above them are derived from the first matrix-set, and objects with right-pointing arrows above them are derived from the second.

{\bf Signed subterms in determinant expansions}. Given a permutation $\alpha$ of $\gamma$, consider the term (possibly zero) $\overrightarrow{T}_\alpha$ in the expansion of $\mathcal{B}[\delta|\gamma]$. We can think of $\overrightarrow{T}_\alpha$ as a sum of signed subterms which can be enumerated in various ways. For example, for $j = 1, \ldots, 2^{|\delta|}$, define $r^{(j)}$ to be the integer $j-1$ written out as a binary string with $|\delta|$ digits. Then define the subentries $\overrightarrow{T}_{\delta_i\alpha_i}^{(j)} = \mathcal{B}_{\delta_i\alpha_i}^{+}$ if $r^{(j)}$ has a zero in $j$th place,  and $\overrightarrow{T}_{\delta_i\alpha_i}^{(j)} = \mathcal{B}_{\delta_i\alpha_i}^{-}$ if $r^{(j)}$ has a one in $j$th place. Many of these subentries may be zero if the corresponding entries are signed. Now $\overrightarrow{T}_\alpha^{(j)} = P(\alpha)\prod_{i=1}^{|\delta|}\overrightarrow{T}_{\delta_i\alpha_i}^{(j)}$ and $\overrightarrow{T}_\alpha = \sum_{j=1}^{2^{|\delta|}}\overrightarrow{T}_{\alpha}^{(j)}$. The same methodology can be applied to give signed subterms $\overleftarrow{T}_{\alpha}^{(j)}$ in $\mathcal{A}[\delta|\gamma]$. From now on a ``signed subterm'' in the expansion of a determinant will mean a signed subterm not identically zero. 

If all entries $\mathcal{B}_{\delta_i\alpha_i}$ are nonzero and signed, then there is a unique signed subterm corresponding to $\alpha$, i.e. $\overrightarrow{T}_\alpha = \overrightarrow{T}_{\alpha}^{(j)}$ for some $j$. Thus we can refer to $\overrightarrow{T}_{\alpha}^{(j)}$ without further comment, to mean a signed subterm corresponding to permutation $\alpha$.

{\bf Edges in the DSR graph are associated with signed subentries.} We can associate with any nonzero signed subentry $\overleftarrow{T}_{\delta_i\beta_i}^{(j)}$ in $\mathcal{A}(\delta|\gamma)$ an S-to-R edge in $G(\delta|\gamma)$, which we term $\overleftarrow{E}_{\delta_i\beta_i}^{(j)}$. Similarly given some nonzero signed subentry $\overrightarrow{T}_{\delta_i\alpha_i}^{(k)}$ in $\mathcal{B}(\delta|\gamma)$ we get an R-to-S edge in $G(\delta|\gamma)$, termed $\overrightarrow{E}_{\delta_i\alpha_i}^{(k)}$. If for some $\alpha$, $j$, $\overrightarrow{T}_{\delta_i\alpha_i}^{(j)}$ and $\overleftarrow{T}_{\delta_i\alpha_i}^{(j)}$ are both nonzero, then by construction they have the same sign, and hence correspond to the same (undirected) edge in $G(\delta|\gamma)$, which will be termed $E_{\delta_i\alpha_i}^{(j)}$.

{\bf Signed term subgraphs.} We can associate with any signed subterm $\overrightarrow{T}_{\alpha}^{(j)}$ in $\mathcal{B}[\delta|\gamma]$ an S-to-R edge-set in $G(\delta|\gamma)$ which will be called $\overrightarrow{E}_\alpha^{(j)}$. Similarly for any signed subterm $\overleftarrow{T}_\alpha^{(j)}$ in $\mathcal{A}[\delta|\gamma]$ there is an R-to-S edge-set in $G(\delta|\gamma)$, $\overleftarrow{E}_\alpha^{(j)}$. Such edge-sets are a generalisation of term subgraphs introduced in \cite{banajicraciun} and will be called signed term subgraphs. Note that when we associate a signed term subgraph with a particular signed subterm we import this subgraph from the DSR graph including all directions and edge-labels. Further note that it only makes sense to talk about a signed term subgraph of an SR or DSR graph $G$, when $G$ is square. Given a square DSR graph $G$, we will say that a signed term subgraph $E$ in $G$ {\bf bisects} a cycle $C$ if it contains one half of the edges in $C$. This is only possible if it contains one member of a disconnecting partition of $C$.

The above notions and notation carry over to SR graphs, except that in this case we omit the arrows indicating directionality. So given a matrix-set $\mathcal{A}(\delta|\gamma)$ and associated SR graph $G_{\mathcal{A}(\delta|\gamma)}$, $T_{\alpha}$ is the subterm in $\mathcal{A}[\delta|\gamma]$ corresponding to permutation $\alpha$ of $\gamma$, $T_{\alpha}^{(j)}$ is the $j$th signed subterm in $T_{\alpha}$, and $E_{\alpha}^{(j)}$ is the corresponding signed term subgraph of $G_{\mathcal{A}(\delta|\gamma)}$. 

Many of the results to follow rely on the fact that the union of two signed term subgraphs, one with S-to-R direction and one with R-to-S direction, results in a set of cycles in the DSR graph:
\begin{lem}
\label{formaltogenuine}
Consider two oppositely directed, signed term subgraphs $\overrightarrow{E}_\alpha^{(r)}$ and $\overleftarrow{E}_\beta^{(s)}$ in some subgraph of a DSR graph. The union $\overrightarrow{E}_\alpha^{(r)} \cup \overleftarrow{E}_\beta^{(s)}$, regarded as a subgraph, consists of a set of vertex-disjoint components, each of which is either an isolated edge or a genuine cycle. 
\end{lem}
\begin{proof}
Each vertex in $\overrightarrow{E}_\alpha^{(r)} \cup \overleftarrow{E}_\beta^{(s)}$ either has incident on it two directed edges (one from $\overrightarrow{E}_\alpha^{(r)}$ and one from $\overleftarrow{E}_\beta^{(s)}$), or one undirected edge. Firstly, if two vertices in $\overrightarrow{E}_\alpha^{(r)} \cup \overleftarrow{E}_\beta^{(s)}$ are connected by an undirected edge, then there are no other edges incident on either vertex. Secondly, no two formal cycles can intersect, as then there would be a vertex with three edges incident on it. Finally, all formal cycles are genuine: consider any formal cycle consisting of edges $[e_1, e_2, \ldots, e_{2r}]$. By the definition of a signed term subgraph, it is not possible for $e_i$ and $e_{(i \mod 2r)+1}$ to belong to the same signed term subgraph, and thus the edges must have alternating S-to-R and R-to-S direction. 
\end{proof}

\subsection{How the conditions in Lemmas~\ref{genlem}~and~\ref{genlemconv} can fail}

\nopagebreak
Following the notation and terminology in Section~\ref{matsets}, consider two sets of square matrices of the same dimension $\mathcal{A} = \{\widehat{\mathcal{A}}(x)\,|\,x \in X\}$, $\mathcal{B} = \{\widehat{\mathcal{B}}(x)\,|\,x \in X\}$. We refer to $(\mathcal{A}, \mathcal{B})$ as a {\bf failed pair} if there is some $x_0 \in X$ such that $A_0 \equiv \widehat{\mathcal{A}}(x_0)$, and $B_0 \equiv \widehat{\mathcal{B}}(x_0)$ satisfy $\mathrm{det}(A_0)\mathrm{det}(B_0) < 0$. $(A_0, B_0)$ will be termed a {\bf failed instance} of $(\mathcal{A}, \mathcal{B})$. With this terminology, some pair $(\mathcal{S}, \mathcal{V})$ fail the conditions in Lemmas~\ref{genlem}~and~\ref{genlemconv} iff there is some $\delta \subset \{1, \ldots, n\}$, $\gamma \subset \{1, \ldots, m\}$ such that $(\mathcal{S}(\delta|\gamma), (-\mathcal{V})(\gamma|\delta))$ are a failed pair. 

Consider a failed pair $(\mathcal{A}, \mathcal{B})$ with failed instance $(A_0, B_0)$. There must be a nonempty set of signed subterms $\mathcal{T}_{\mathcal{A}}$ in the expansion of $\mathrm{det}(\mathcal{A})$ such that for any $\overleftarrow{T} \in \mathcal{T}_{\mathcal{A}}$, $\overleftarrow{T}\mathrm{det}(A_0) \subset \mathbb{R}_{\geq 0}$, and similarly a nonempty set of signed subterms $\mathcal{T}_{\mathcal{B}}$ in the expansion of $\mathrm{det}(\mathcal{B})$ such that for any $\overrightarrow{T} \in \mathcal{T}_{\mathcal{B}}$, $\overrightarrow{T}\mathrm{det}(B_0) \subset \mathbb{R}_{\geq 0}$. We will call the pair $(\mathcal{T}_{\mathcal{A}}, \mathcal{T}_{\mathcal{B}})$ {\bf failed subterms}.

\section{Graph-theoretic results}

\nopagebreak
All results in this section build towards a proof of Theorem~\ref{mainthm}. Two key preliminaries, Lemmas~\ref{allocycles1}~and~\ref{oneeands} below, are minor generalisations of results in \cite{banajicraciun}. The proofs are similar to those in \cite{banajicraciun}, and can be found in Appendix~\ref{app2}.

\begin{lem}
\label{allocycles1}
Consider two sets of $k \times k$ matrices $\mathcal{A}$ and $\mathcal{B}$, and the associated DSR graph $G \equiv G_{\mathcal{A},\mathcal{B}}$. Consider signed subterms $\overleftarrow{T}_{\alpha}^{(r)}$ in the expansion of $\mathrm{det}(\mathcal{A})$ and $\overrightarrow{T}_{\beta}^{(s)}$ in the expansion of $\mathrm{det}(\mathcal{B})$. If all cycles in $\overleftarrow{E}_\alpha^{(r)}\cup \overrightarrow{E}_\beta^{(s)}$ are o-cycles then $\overleftarrow{T}_\alpha^{(r)}\overrightarrow{T}_\beta^{(s)} \subset \mathbb{R}_{\geq 0}$.
\end{lem}

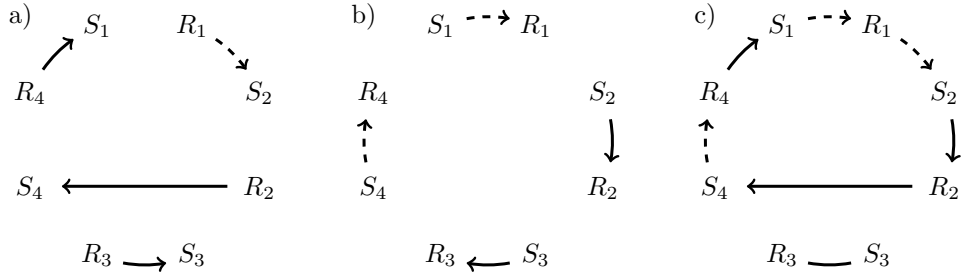
\begin{figure}[h]
\begin{minipage}{\textwidth}

\begin{minipage}{0.3\textwidth}
\begin{tikzpicture}[domain=-4:4,scale=0.55]

\node at (-3,3) {a)};

\path (67.5: 3cm) coordinate(R1);
\path (22.5: 3cm) coordinate (S2);
\path (-22.5: 3cm) coordinate (R2);
\path (-67.5: 3cm) coordinate (S3);
\path (-112.5: 3cm) coordinate (R3);
\path (-157.5: 3cm) coordinate (S4);
\path (157.5: 3cm) coordinate (R4);
\path (112.5: 3cm) coordinate (S1);
\node at (R1) {$R_1$};
\node at (R2) {$R_2$};
\node at (R3) {$R_3$};
\node at (R4) {$R_4$};
\node at (S1) {$S_1$};
\node at (S2) {$S_2$};
\node at (S3) {$S_3$};
\node at (S4) {$S_4$};

\path (67.5-12.5: 3cm) coordinate(R1end);
\path (22.5-12.5: 3cm) coordinate (S2end);
\path (-22.5-12.5: 3cm) coordinate (R2end);
\path (-67.5-12.5: 3cm) coordinate (S3end);
\path (-112.5-12.5: 3cm) coordinate (R3end);
\path (-157.5-12.5: 3cm) coordinate (S4end);
\path (157.5-12.5: 3cm) coordinate (R4end);
\path (112.5-12.5: 3cm) coordinate (S1end);

\draw[->, dashed, line width=0.04cm] (R1end)  arc (67.5-12.5:67.5-32.5:3cm);
\draw[<-, line width=0.04cm] (S3end)  arc (-67.5-12.5:-67.5-32.5:3cm);
\draw[->, line width=0.04cm] (R4end)  arc (157.5-12.5:157.5-32.5:3cm);

\draw[->,line width=0.04cm] (2, -1.1) -- (-2, -1.1);

\end{tikzpicture}
\end{minipage}
\hfill
\begin{minipage}{0.3\textwidth}
\begin{tikzpicture}[domain=-4:4,scale=0.55]

\node at (-3,3) {b)};

\path (67.5: 3cm) coordinate(R1);
\path (22.5: 3cm) coordinate (S2);
\path (-22.5: 3cm) coordinate (R2);
\path (-67.5: 3cm) coordinate (S3);
\path (-112.5: 3cm) coordinate (R3);
\path (-157.5: 3cm) coordinate (S4);
\path (157.5: 3cm) coordinate (R4);
\path (112.5: 3cm) coordinate (S1);
\node at (R1) {$R_1$};
\node at (R2) {$R_2$};
\node at (R3) {$R_3$};
\node at (R4) {$R_4$};
\node at (S1) {$S_1$};
\node at (S2) {$S_2$};
\node at (S3) {$S_3$};
\node at (S4) {$S_4$};

\path (67.5-12.5: 3cm) coordinate(R1end);
\path (22.5-12.5: 3cm) coordinate (S2end);
\path (-22.5-12.5: 3cm) coordinate (R2end);
\path (-67.5-12.5: 3cm) coordinate (S3end);
\path (-112.5-12.5: 3cm) coordinate (R3end);
\path (-157.5-12.5: 3cm) coordinate (S4end);
\path (157.5-12.5: 3cm) coordinate (R4end);
\path (112.5-12.5: 3cm) coordinate (S1end);

\draw[->, line width=0.04cm] (S2end)  arc (22.5-12.5:22.5-32.5:3cm);
\draw[->, line width=0.04cm] (S3end)  arc (-67.5-12.5:-67.5-32.5:3cm);
\draw[->, dashed, line width=0.04cm] (S4end)  arc (-157.5-12.5:-157.5-32.5:3cm);
\draw[->,dashed, line width=0.04cm] (S1end)  arc (112.5-12.5:112.5-32.5:3cm);

\end{tikzpicture}
\end{minipage}
\hfill
\begin{minipage}{0.3\textwidth}
\begin{tikzpicture}[domain=-4:4,scale=0.55]

\node at (-3,3) {c)};

\path (67.5: 3cm) coordinate(R1);
\path (22.5: 3cm) coordinate (S2);
\path (-22.5: 3cm) coordinate (R2);
\path (-67.5: 3cm) coordinate (S3);
\path (-112.5: 3cm) coordinate (R3);
\path (-157.5: 3cm) coordinate (S4);
\path (157.5: 3cm) coordinate (R4);
\path (112.5: 3cm) coordinate (S1);
\node at (R1) {$R_1$};
\node at (R2) {$R_2$};
\node at (R3) {$R_3$};
\node at (R4) {$R_4$};
\node at (S1) {$S_1$};
\node at (S2) {$S_2$};
\node at (S3) {$S_3$};
\node at (S4) {$S_4$};

\path (67.5-12.5: 3cm) coordinate(R1end);
\path (22.5-12.5: 3cm) coordinate (S2end);
\path (-22.5-12.5: 3cm) coordinate (R2end);
\path (-67.5-12.5: 3cm) coordinate (S3end);
\path (-112.5-12.5: 3cm) coordinate (R3end);
\path (-157.5-12.5: 3cm) coordinate (S4end);
\path (157.5-12.5: 3cm) coordinate (R4end);
\path (112.5-12.5: 3cm) coordinate (S1end);

\draw[->, dashed, line width=0.04cm] (R1end)  arc (67.5-12.5:67.5-32.5:3cm);
\draw[->, line width=0.04cm] (S2end)  arc (22.5-12.5:22.5-32.5:3cm);
\draw[-, line width=0.04cm] (S3end)  arc (-67.5-12.5:-67.5-32.5:3cm);
\draw[->, dashed, line width=0.04cm] (S4end)  arc (-157.5-12.5:-157.5-32.5:3cm);
\draw[->, line width=0.04cm] (R4end)  arc (157.5-12.5:157.5-32.5:3cm);
\draw[->,dashed, line width=0.04cm] (S1end)  arc (112.5-12.5:112.5-32.5:3cm);

\draw[->,line width=0.04cm] (2, -1.1) -- (-2, -1.1);

\end{tikzpicture}
\end{minipage}
\end{minipage}
\caption{\label{lemallofig}a) A signed term subgraph in the DSR graph $G_{\mathcal{A}, \mathcal{B}}$ corresponding to a signed subterm in $\mathrm{det}(\mathcal{A})$ in (\ref{alloexample}). b) A signed term subgraph in $G_{\mathcal{A}, \mathcal{B}}$ corresponding to a signed subterm in $\mathrm{det}(\mathcal{B})$ in (\ref{alloexample}). c) The union of these term subgraphs. As the terms are oppositely signed, the union of the corresponding term subgraphs contains an e-cycle. Edge-labels have been omitted.}
\end{figure}

Lemma~\ref{allocycles1} tells us that a pair of oppositely directed, signed term subgraphs in a DSR graph derived from oppositely signed subterms must contain in their union some e-cycles. For example, define the following pair of matrix-sets:
\begin{equation}
\label{alloexample}
\mathcal{A} = \left[\begin{array}{cccc}* & * & * & a_1\\-a_2& * & * & * \\ * & * & a_3 & * \\ * & a_4 & * & * \end{array}\right] \qquad \mathcal{B} = \left[\begin{array}{cccc}-b_1 & * & * & *\\ * & b_2 & * & *\\ * & * & b_3 & * \\ * & * & * & -b_4\end{array}\right]
\end{equation}
with $a_i, b_i > 0$, and $*$ indicating entries of arbitrary magnitude and sign. Consider the corresponding DSR graph $G_{\mathcal{A},\mathcal{B}}$. The entries $a_i$ in $\mathcal{A}$ define an R-to-S signed term subgraph in $G_{\mathcal{A},\mathcal{B}}$ (Figure~\ref{lemallofig}a), and the entries $b_i$ in $\mathcal{B}$ define an S-to-R signed term subgraph in $G_{\mathcal{A},\mathcal{B}}$ (Figure~\ref{lemallofig}b). As the corresponding terms in $\mathrm{det}(\mathcal{A})$ and $\mathrm{det}(\mathcal{B})$ have opposite sign, the union of these term subgraphs contains an e-cycle (Figure~\ref{lemallofig}c).

{\bf Remark.} As the DSR graph does not contain information on whether $\mathcal{A}$ and $\mathcal{B}$ are independent, Lemma~\ref{allocycles1} is actually stronger than it appears. If all cycles in $\overleftarrow{E}_\alpha^{(r)}\cup \overrightarrow{E}_\beta^{(s)}$ are o-cycles, then in fact fixing any $A \in \mathcal{A}, B \in \mathcal{B}$ and defining $\overleftarrow{T}_\alpha = P(\alpha)\prod_{i=1}^{k}A_{i\alpha_i}$ and $\overrightarrow{T}_\beta = P(\beta)\prod_{i=1}^{k}B_{i\beta_i}$, we have $\overleftarrow{T}_\alpha\overrightarrow{T}_\beta \geq 0$. The same holds true for each other DSR graph result: if a result about a pair of matrix-sets $(\mathcal{A}, \mathcal{B})$ can be inferred from the DSR graph $G_{\mathcal{A}, \mathcal{B}}$, then it is true not only for each pair $(\widehat{\mathcal{A}}(x), \widehat{\mathcal{B}}(x))$ ($x \in X$), but also for each pair $(\widehat{\mathcal{A}}(x_1), \widehat{\mathcal{B}}(x_2))$ ($x_1, x_2 \in X$).

\begin{lem}
\label{oneeands}
Consider any set of $k \times k$ matrices $\mathcal{A}$ with associated SR graph $G$. Let $\alpha$ and $\beta$ be permutations of $\{1, \ldots, k\}$ such that $T_\alpha^{(r)}$ and $T_\beta^{(s)}$ are signed subterms in the determinant expansion of $\mathcal{A}$. Assume that $E_\alpha^{(r)} \cup E_\beta^{(s)}$ contains exactly one cycle $C$, and this cycle is an es-cycle. Then $T_\alpha^{(r)} + T_\beta^{(s)} = 0$. 
\end{lem}

Lemma~\ref{oneeands} shows that having es-cycles in an SR-graph means that some terms in a determinant expansion sum to zero. The result is not affected by the possibility that some edge-labels in the graph may be $\infty$, but obviously such labels cannot occur in es-cycles themselves.

\begin{lem}
\label{pairterms}
Consider a set of square matrices $\mathcal{A}$ and the associated SR graph $G_{\mathcal{A}}$. Assume that $G_{\mathcal{A}}$ has an es-cycle $C$. Let $\mathcal{E}_C$ be the set of all signed term subgraphs in $G_{\mathcal{A}}$ which bisect $C$, and $\mathcal{T}_C$ the corresponding signed subterms in $\mathrm{det}(\mathcal{A})$. Then 
\[
\sum_{T \in \mathcal{T}_C} T = 0\,.
\]
\end{lem}
\begin{proof}
All edges in $C$ correspond to signed subentries in $\mathcal{A}$. If $\mathcal{E}_C$ is empty then we are done. Otherwise consider $E \in \mathcal{E}_C$ with corresponding signed subterm $T$. Construct the new signed term subgraph $\tilde{E} = (E\backslash C) \cup (C\backslash E)$ with corresponding signed subterm $\tilde{T}$. Clearly $\tilde{E} \in \mathcal{E}_C$. Now $E \cup \tilde{E}$ contains a single e-cycle $C$ which is an s-cycle, so by Lemma~\ref{oneeands}, $T + \tilde{T} = 0$. All signed subterms in $\mathcal{T}_C$ pair off in this way, so $\sum_{T \in \mathcal{T}_C} T = 0$. 
\end{proof}

Note that the above lemma does not imply that $\mathrm{det}(\mathcal{A}) = 0$, as not all signed term subgraphs in $G_{\mathcal{A}}$ necessarily bisect $C$.

\begin{lem}
\label{twodisjointscycles}
Consider a set of square matrices $\mathcal{A}$ and the associated SR graph $G_{\mathcal{A}}$. Let $\mathcal{C}$ be any set of {\bf edge-disjoint} es-cycles in $G$. Consider all signed subterms in the expansion of $\mathrm{det}(\mathcal{A})$ such that the corresponding signed term subgraphs bisect some es-cycle in $\mathcal{C}$. These terms all sum to zero. 
\end{lem}
\begin{proof}
Let $\mathcal{C} = \{C^{(1)}, \ldots, C^{(k)}\}$. Define $\mathcal{E}_i$ to be the set of all signed term subgraphs which bisect $C^{(i)}$, let $\mathcal{T}_i$ be the set of corresponding signed subterms, and let $\mathcal{T} = \cup_{i=1}^k \mathcal{T}_i$. By Lemma~\ref{pairterms}, all terms in $\mathcal{T}_i$ sum to zero. 

We will refer to a pair of signed term subgraphs $E$ and $\tilde{E}$ in $\mathcal{E}_i$ such that $(E \cup \tilde{E})\backslash(E \cap \tilde{E}) = C^{(i)}$ as an $i$-pair. Now consider some signed term subgraph $E \in \mathcal{E}_i \cap \mathcal{E}_j$ for some $i, j$. Because the $C^{(j)}$ are all edge-disjoint, we have the following partition of $E$:
\[
E = (E \cap C^{(i)}) \cup (E \cap C^{(j)}) \cup (E\backslash (C^{(i)} \cup C^{(j)}))\,.
\]
Let $\tilde{E} = (E \backslash C^{(i)}) \cup (C^{(i)} \backslash E)$ be the other member of the $i$-pair corresponding to $E$. We know that $(E \cap C^{(j)}) \subset (E \backslash C^{(i)})$, so 
\[
\tilde{E} \cap C^{(j)} = (E\backslash C^{(i)}) \cap C^{(j)} = E \cap C^{(j)}
\]
which makes it clear that $\tilde{E} \in \mathcal{E}_j$. So if one member of an $i$-pair is in $\mathcal{E}_j$, then so is the other. Since $\mathcal{E}_i \cap \mathcal{E}_j$ consists of $i$-pairs, it follows that $\mathcal{E}_i \backslash \mathcal{E}_j$ consists of $i$-pairs. By induction, if $i > 1$, then  $\mathcal{E}_i \backslash \cup_{j=1}^{i-1}\mathcal{E}_j$ consists of $i$-pairs. So, by Lemma~\ref{pairterms},
\[
\sum_{T \in \mathcal{T}_i\backslash \cup_{j=1}^{i-1}\mathcal{T}_j}T = 0.
\]

We have the partition $\mathcal{T} = \mathcal{T}_1 \cup (\mathcal{T}_2\backslash \mathcal{T}_1) \cup (\mathcal{T}_3\backslash (\mathcal{T}_1 \cup \mathcal{T}_2)) \cup \cdots$, and hence
\[
\sum_{T \in \mathcal{T}}T = \sum_{i=1}^k\sum_{T \in \mathcal{T}_i\backslash \cup_{j=1}^{i-1}\mathcal{T}_j}T\,\,\,= \,\,\,\sum_{i=1}^k0 \,\,\, = \,\,\,0.
\]
\end{proof}

\begin{cor}
\label{disjointscycles}
Consider a set of square matrices $\mathcal{A}$ and the associated SR graph $G_{\mathcal{A}}$. Consider any set $\mathcal{C}$ of {\bf edge-disjoint} es-cycles in $G_{\mathcal{A}}$, and assume that each signed term subgraph in $G_{\mathcal{A}}$ bisects some es-cycle from $\mathcal{C}$. Then all matrices in $\mathcal{A}$ are singular. 
\end{cor}

\begin{proof}
This follows immediately from the previous result. 
\end{proof}

\begin{figure}[h]
\begin{minipage}{0.3\textwidth}
\[
\left[\begin{array}{cccc}1& 0 & 0 & 2\\ * & * & * & 0\\0 & * & * & 0\\1 & 0 & 0 & 2 \end{array}\right]
\]

\end{minipage}
\hfill
\begin{minipage}{0.65\textwidth}
\begin{tikzpicture}[domain=0:4,scale=0.6]

\node at (1,4) {$S_1$};
\node at (4,4) {$R_1$};
\node at (7,4) {$S_2$};
\node at (10,4) {$R_3$};
\node at (1,1) {$R_4$};
\node at (4,1) {$S_4$};
\node at (7,1) {$R_2$};
\node at (10,1) {$S_3$};

\draw[-, line width=0.04cm] (1.5, 4) -- (3.5,4);
\draw[-, line width=0.04cm] (1.5, 1) -- (3.5,1);
\draw[-, line width=0.04cm] (1, 3.5) -- (1,1.5);
\draw[-, line width=0.04cm] (4, 3.5) -- (4,1.5);

\draw [thick,snake=snake,segment amplitude=.4mm,
         segment length=2mm] (4.5, 4) -- (6.5,4);
\draw [thick,snake=snake,segment amplitude=.4mm,
         segment length=2mm] (7.5, 4) -- (9.5,4);
\draw [thick,snake=snake,segment amplitude=.4mm,
         segment length=2mm] (7.5, 1) -- (9.5,1);
\draw [thick,snake=snake,segment amplitude=.4mm,
         segment length=2mm] (10, 3.5) -- (10,1.5);
\draw [thick,snake=snake,segment amplitude=.4mm,
         segment length=2mm] (7, 3.5) -- (7,1.5);

\node at (2.5, 4.3) {$\scriptstyle{1}$};
\node at (2.5, 1.3) {$\scriptstyle{2}$};

\node at (1.3, 2.5) {$\scriptstyle{2}$};
\node at (4.3, 2.5) {$\scriptstyle{1}$};
\end{tikzpicture}
\end{minipage}
\caption{\label{exoneeands} A matrix-set and the corresponding SR graph. The entries $*$ are of unknown magnitude and sign (and may be unsigned), and thus wavy edges in the SR graph may correspond to single edges or pairs of edges, and have unknown edge-labels. However, it is easy to see that every signed term subgraph in the DSR graph bisects the es-cycle $S_1\!-\!R_1\!-\!S_4\!-\!R_4$, and thus the matrices are singular.}
\end{figure}
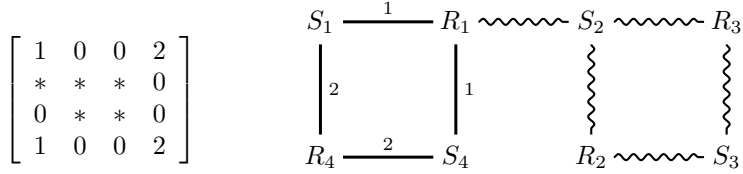

The matrix-set and corresponding SR graph in Figure~\ref{exoneeands} provide an illustration of the previous results.

\begin{lem}
\label{disjointscycles2}
Consider a pair $(\mathcal{A},\mathcal{B})$ of sets of $k \times k$ matrices, and the associated DSR graph $G = G_{\mathcal{A},\mathcal{B}}$. Assume that $(\mathcal{A}, \mathcal{B})$ are a failed pair with failed instance $(A_0, B_0)$ and failed subterms $(\mathcal{T}_{\mathcal{A}},\mathcal{T}_{\mathcal{B}})$ corresponding to sets of signed term subgraphs $(\mathcal{E}_{\mathcal{A}},\mathcal{E}_{\mathcal{B}})$. Choose some $\overrightarrow{T}_* \in \mathcal{T}_{\mathcal{B}}$ with corresponding signed term subgraph $\overrightarrow{E}_* \in \mathcal{E}_{\mathcal{B}}$. Let $\mathcal{C}$ be the set of all e-cycles in $G$ each of which lies in $\overrightarrow{E}_* \cup \overleftarrow{E}$ for some $\overleftarrow{E} \in \mathcal{E}_{\mathcal{A}}$. Then
\begin{enumerate}
\item $\mathcal{C}$ contains an e-cycle which fails to be an s-cycle, or
\item $\mathcal{C}$ contains two e-cycles which fail to be edge-disjoint.
\end{enumerate}
\end{lem}
\begin{proof}
$\mathcal{T}_{\mathcal{A}}$ and $\mathcal{T}_{\mathcal{B}}$ are nonempty (by the definition of a failed pair), and so, by Lemma~\ref{allocycles1}, $\mathcal{C}$ is not empty. Assume the result is false, i.e. $\mathcal{C}$ consists of edge-disjoint es-cycles $C^{(1)}, \ldots, C^{(r)}$. All edges in all cycles in $\mathcal{C}$ correspond to edges with the same edge-labels in $G_{\mathcal{A}}$ (the SR graph corresponding to $\mathcal{A}$): otherwise some edge must be imported only from $G_{\mathcal{B}}$, and hence carry an edge-label of $\infty$, causing some cycle in $\mathcal{C}$ to fail to be an s-cycle. Let $\mathcal{E}_i$ be the set of all signed term subgraphs in $G_{\mathcal{A}}$ which bisect $C^{(i)}$, with corresponding signed subterms $\mathcal{T}_i$, and let $\mathcal{T}^{'} = \cup \mathcal{T}_i$. By construction, $\mathcal{T}^{'} \supset \mathcal{T}_{\mathcal{A}}$. By Lemma~\ref{twodisjointscycles}, the sum of terms in $\mathcal{T}^{'}$ is zero. Since $A_0$ is nonsingular, there must be some signed subterm $\overleftarrow{T} \not \in \mathcal{T}^{'}$ such that $\overleftarrow{T}\mathrm{det}(A_0)\subset \mathbb{R}_{\geq 0}$. By definition, $\overleftarrow{T} \in \mathcal{T}_{\mathcal{A}}$, contradicting the fact that $\mathcal{T}^{'} \supset \mathcal{T}_{\mathcal{A}}$.
\end{proof}

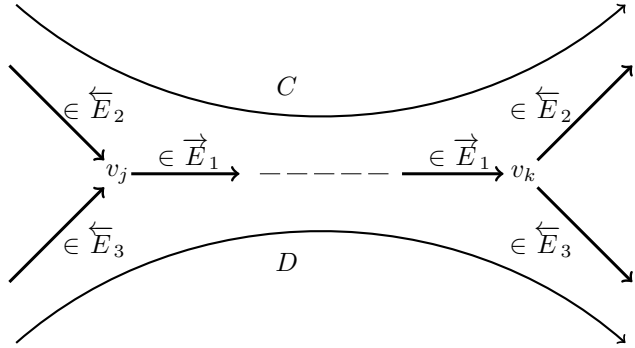
\begin{figure}[h]
\begin{center}
\begin{tikzpicture}[domain=0:4,scale=0.9]

\node at (2,3) {$v_j$};
\node at (8,3) {$v_k$};
\node at (1.6,4) {${}\in \overleftarrow{E}_2$};
\node at (1.6,2) {${}\in \overleftarrow{E}_3$};
\node at (3,3.3) {${}\in \overrightarrow{E}_1$};
\node at (7,3.3) {${}\in \overrightarrow{E}_1$};
\node at (8.2,4) {${}\in \overleftarrow{E}_2$};
\node at (8.2,2) {${}\in \overleftarrow{E}_3$};

\node at (4.5,4.3) {$C$};
\node at (4.5,1.7) {$D$};
\draw[->, thick] (0.5,5.5) .. controls (3,3.3) and (7,3.3) .. (9.5,5.5);
\draw[->, thick] (0.5,0.5) .. controls (3,2.7) and (7,2.7) .. (9.5,0.5);

\draw[->, line width=0.04cm] (0.4,4.6) -- (1.8,3.2);
\draw[->, line width=0.04cm] (0.4,1.4) -- (1.8,2.8);
\draw[->, line width=0.04cm] (2.2,3) -- (3.8,3);
\draw[very thin] (4.1,3) -- (4.4,3);
\draw[very thin] (4.5,3) -- (4.8,3);
\draw[very thin] (4.9,3) -- (5.2,3);
\draw[very thin] (5.3,3) -- (5.6,3);
\draw[very thin] (5.7,3) -- (6.0,3);
\draw[->, line width=0.04cm] (6.2,3) -- (7.7,3);
\draw[->, line width=0.04cm] (8.2,3.2) -- (9.6,4.6);
\draw[->, line width=0.04cm] (8.2,2.8) -- (9.6,1.4);

\end{tikzpicture}
\end{center}
\caption{\label{StoRfig} A portion of the DSR graph corresponding to the situation in Theorem~\ref{StoR2}. Because any vertex with three edges incident on it must have one from {$E_1$}, one from $E_2$ and one from $E_3$, this forces any component of the intersection between $C$ and $D$ to be of odd length.}
\end{figure}

The next result is about the geometry of subgraphs of DSR graphs constructed as the union of exactly three signed term subgraphs:

\begin{lem}
\label{StoR2}
Consider a square DSR graph $G$ containing signed term subgraphs, $\overrightarrow{E}_1$, $\overleftarrow{E}_2$ and $\overleftarrow{E}_3$. Assume that there is a cycle $C$ in $\overrightarrow{E}_1 \cup \overleftarrow{E}_2$ and another, distinct, cycle $D$ in $\overrightarrow{E}_1 \cup \overleftarrow{E}_3$. $C$ and $D$ must either be edge and vertex-disjoint, or must have S-to-R intersection. 
\end{lem}
\begin{proof}
Each edge in $C$ lies in exactly one of $\overrightarrow{E}_1$ or $\overleftarrow{E}_2$, and similarly each edge in $D$ lies in exactly one of $\overrightarrow{E}_1$ or $\overleftarrow{E}_3$. Consequently, any edge in $C \cap D$ lies either in $\overrightarrow{E}_1$ or in $\overrightarrow{E}_2 \cap \overleftarrow{E}_3$, but not in both. Thus $C$ and $D$ each have a natural orientation and any edge in $C \cap D$ has the same $C$-orientation and $D$-orientation (i.e. $C$ and $D$ are compatibly oriented). Define $G_{123} = \overrightarrow{E}_1 \cup \overleftarrow{E}_2 \cup \overleftarrow{E}_3$, and ignore signs and labels on edges. No vertex from $G_{123}$ can have more than three edges from $G_{123}$ incident on it, so it is impossible for $C$ and $D$ to have a vertex but no edges in their intersection as then this vertex would have four edges incident on it. If $C$ and $D$ are vertex-disjoint then we are done. 

So assume that $C$ and $D$ share a vertex and choose any such vertex $v_0$. Since $C$ and $D$ are distinct, following $C$ backwards from $v_0$ we must come to a first vertex $v_j$ (possibly $v_0$) with three edges from $G_{123}$ incident on it, two incoming edges (in $C \backslash D$ and $D \backslash C$), and one outgoing edge (in $C \cap D$). Since $\overrightarrow{E}_1$, $\overleftarrow{E}_2$ and $\overleftarrow{E}_3$ are term subgraphs, the outgoing edge must lie in $\overrightarrow{E}_1$, while the incoming edges must lie in $\overrightarrow{E}_2$ and $\overleftarrow{E}_3$. Similarly following $C$ forwards from $v_0$ we must come to a vertex $v_k$ (possibly $v_0$) with three edges from $G_{123}$ incident on it, one incoming edge (in $C \cap D$) and two outgoing edges (in $C \backslash D$ and $D \backslash C$). This time the incoming edge must lie in $\overrightarrow{E}_1$, while the outgoing edges must lie in $\overrightarrow{E}_2$ and $\overleftarrow{E}_3$. (see Figure~\ref{StoRfig}). Since it starts and ends with an edge from $\overrightarrow{E}_1$, the path in $C \cap D$ from $v_j$ to $v_k$ has an odd number of edges. As $v_0$ was an arbitrary vertex in $C \cap D$, $C$ and $D$ have S-to-R intersection.
\end{proof}

As an immediate corollary of Lemma~\ref{StoR2}:
\begin{cor}
\label{StoR3}
Consider a DSR graph $G$ with a signed, directed term subgraph $\overrightarrow{E}_1$. Let $\overleftarrow{E}_2, \ldots, \overleftarrow{E}_k$ ($k \geq 2$) be a set of oppositely directed term subgraphs in $G$. Define the set of e-cycles $\mathcal{C}$ in $G$ as follows: a cycle $C$ is in $\mathcal{C}$ iff $C \subset (\overrightarrow{E}_1 \cup \overleftarrow{E}_j)$ for some $j \in \{2, \ldots, k\}$. Then any two cycles in $\mathcal{C}$ are either edge and vertex-disjoint, or have S-to-R intersection. 
\end{cor}
\begin{proof}
Consider two cycles, $C$ and $D$ in $\mathcal{C}$. Let $\overleftarrow{E}_i$ be some signed term subgraph such that $\overleftarrow{E}_i \cup \overrightarrow{E}_1 \supset C$ and $\overleftarrow{E}_j$ be some signed term subgraph such that $\overleftarrow{E}_j \cup \overrightarrow{E}_1 \supset D$. If $i=j$, then clearly $C$ and $D$ must be edge and vertex-disjoint since no vertex in $\overleftarrow{E}_i \cup \overrightarrow{E}_1$ has more than two edges incident on it. Otherwise, apply Lemma~\ref{StoR2} to $\overrightarrow{E}_1, \overleftarrow{E}_i$ and $\overleftarrow{E}_j$ to get that $C$ and $D$ must either be edge and vertex-disjoint or must have S-to-R intersection. 
\end{proof}

\begin{figure}[h]
\begin{minipage}{0.48\textwidth}
\begin{tikzpicture}[domain=0:4,scale=0.6]

\node at (-1,2.5) {$\overrightarrow{E}_1 \cup\overleftarrow{E}_2 =$};
\node at (1,4) {$S_1$};
\node at (4,4) {$R_3$};
\node at (7,4) {$S_2$};
\node at (1,1) {$R_1$};
\node at (4,1) {$S_3$};
\node at (7,1) {$R_2$};

\draw[->, line width=0.04cm] (4.5, 4) -- (6.5,4);
\draw[<-, line width=0.04cm] (4.5, 1) -- (6.5,1);

\draw[-, line width=0.04cm] (1, 3.5) -- (1,1.5);
\draw[-, line width=0.04cm] (4, 3.5) -- (4,1.5);
\draw[-, line width=0.04cm] (7, 3.5) -- (7,1.5);

\node at (1.3, 2.5) {$\scriptstyle{1}$};
\node at (4.3, 2.5) {$\scriptstyle{1}$};

\node at (5.5, 4.3) {$\scriptstyle{1}$};
\node at (5.5, 1.3) {$\scriptstyle{1}$};

\node at (7.3, 2.5) {$\scriptstyle{1}$};
\end{tikzpicture}
\end{minipage}
\hfill
\begin{minipage}{0.48\textwidth}
\begin{tikzpicture}[domain=0:4,scale=0.6]

\node at (-1,2.5) {$\overrightarrow{E}_1 \cup\overleftarrow{E}_3 =$};

\node at (1,4) {$S_1$};
\node at (4,4) {$R_3$};
\node at (7,4) {$S_2$};
\node at (1,1) {$R_1$};
\node at (4,1) {$S_3$};
\node at (7,1) {$R_2$};

\draw[<-, line width=0.04cm] (1.5, 4) -- (3.5,4);
\draw[->, line width=0.04cm] (1.5, 1) -- (3.5,1);

\draw[-, line width=0.04cm] (1, 3.5) -- (1,1.5);
\draw[-, line width=0.04cm] (4, 3.5) -- (4,1.5);
\draw[-, line width=0.04cm] (7, 3.5) -- (7,1.5);

\node at (2.5, 4.3) {$\scriptstyle{1}$};
\node at (2.5, 1.3) {$\scriptstyle{1}$};

\node at (1.3, 2.5) {$\scriptstyle{1}$};
\node at (4.3, 2.5) {$\scriptstyle{1}$};

\node at (7.3, 2.5) {$\scriptstyle{1}$};
\end{tikzpicture}
\end{minipage}
\begin{minipage}{\textwidth}
\begin{center}
\begin{tikzpicture}[domain=0:4,scale=0.6]

\node at (-2,2.5) {$\overrightarrow{E}_1 \cup\overleftarrow{E}_2\cup\overleftarrow{E}_3\,\,=\,\,$};

\node at (1,4) {$S_1$};
\node at (4,4) {$R_3$};
\node at (7,4) {$S_2$};
\node at (1,1) {$R_1$};
\node at (4,1) {$S_3$};
\node at (7,1) {$R_2$};

\draw[<-, line width=0.04cm] (1.5, 4) -- (3.5,4);
\draw[->, line width=0.04cm] (1.5, 1) -- (3.5,1);
\draw[->, line width=0.04cm] (4.5, 4) -- (6.5,4);
\draw[<-, line width=0.04cm] (4.5, 1) -- (6.5,1);

\draw[-, line width=0.04cm] (1, 3.5) -- (1,1.5);
\draw[-, line width=0.04cm] (4, 3.5) -- (4,1.5);
\draw[-, line width=0.04cm] (7, 3.5) -- (7,1.5);

\node at (2.5, 4.3) {$\scriptstyle{1}$};
\node at (2.5, 1.3) {$\scriptstyle{1}$};

\node at (1.3, 2.5) {$\scriptstyle{1}$};
\node at (4.3, 2.5) {$\scriptstyle{1}$};

\node at (5.5, 4.3) {$\scriptstyle{1}$};
\node at (5.5, 1.3) {$\scriptstyle{1}$};

\node at (7.3, 2.5) {$\scriptstyle{1}$};
\end{tikzpicture}
\end{center}

\end{minipage}
\caption{\label{SRprob} Subgraphs of the DSR graph for the pair $A, \mathcal{B}$ in (\ref{StoReq}) above. By Lemma~\ref{StoR2} the distinct es-cycles are forced to have S-to-R intersection.}
\end{figure}
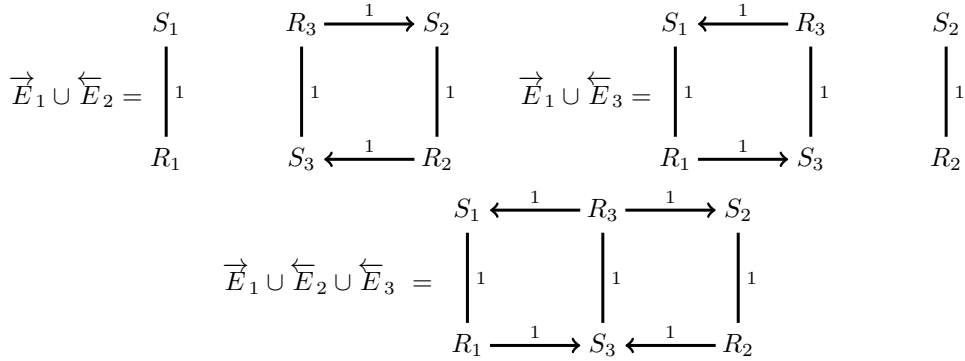

As an illustration of Lemma~\ref{StoR2} and Corollary~\ref{StoR3} consider the following matrix-sets (one of which happens to consist of a single matrix):
\begin{equation}
\label{StoReq}
A = \left[\begin{array}{ccc}1 & 0 & 1\\0 & 1 & 1\\1 & 1 & 1\end{array}\right],\quad \mathcal{B} = \left[\begin{array}{ccc}a & 0 & 0\\0 & b & 0\\0 & 0 & c\end{array}\right]\,\,.
\end{equation}
Assume that $a, b, c > 0$. $A$ and $\mathcal{B}$ are a failed pair since $\mathrm{det}(A) \mathrm{det}(\mathcal{B}) = -abc$. Define $\overrightarrow{T}_1 = \mathcal{B}_{11}\mathcal{B}_{22}\mathcal{B}_{33} = abc$, $\overleftarrow{T}_2 = -A_{11}A_{23}A_{32} = -1$, and $\overleftarrow{T}_3= -A_{13}A_{22}A_{31} = -1$, with corresponding signed term subgraphs $\overrightarrow{E}_1, \overleftarrow{E}_2$ and $\overleftarrow{E}_3$. The subgraphs $\overrightarrow{E}_1 \cup\overleftarrow{E}_2$, $\overrightarrow{E}_1 \cup\overleftarrow{E}_3$ and $\overrightarrow{E}_1 \cup\overleftarrow{E}_2 \cup \overleftarrow{E}_3 $ of the DSR graph are shown in Figure~\ref{SRprob}. Each of $\overrightarrow{E}_1 \cup\overleftarrow{E}_2$, $\overrightarrow{E}_1 \cup\overleftarrow{E}_3$ contains an es-cycle. As these are not edge and vertex-disjoint, they have S-to-R intersection, consisting of the edge $S_3\!-\!R_3$.

We now come to the proof of Theorem~\ref{mainthm} which states that if an interaction network defined by matrix-sets $\mathcal{S}$ and $\mathcal{V}$ fails the conditions in Lemma~\ref{genlem}, then the associated DSR graph fails Condition~($*$).

{\it Proof of Theorem~\ref{mainthm}.} If $\mathcal{S}$ and $\mathcal{V}$ fail the conditions in Lemma~\ref{genlem}, there are some nonempty $\delta \subset \{1, \ldots, n\}$ and $\gamma \subset \{1, \ldots, m\}$ such that $\mathcal{A} \equiv \mathcal{S}(\delta|\gamma)$ and $\mathcal{B} \equiv -\mathcal{V}(\gamma|\delta)$ are a failed pair, with nonempty sets of failed subterms $(\mathcal{T}_{\mathcal{A}},\mathcal{T}_{\mathcal{B}})$. Recall that a term from $\mathcal{T}_{\mathcal{A}}$ and one from $\mathcal{T}_{\mathcal{B}}$ correspond to distinct signed term subgraphs in $G(\delta|\gamma)$. Choose some $\overrightarrow{T}_* \in \mathcal{T}_{\mathcal{B}}$ with corresponding signed term subgraph $\overrightarrow{E}_*$ in $G(\delta|\gamma)$. Let $\mathcal{E}_{\mathcal{A}}$ be the signed term subgraphs corresponding to terms in $\mathcal{T}_{\mathcal{A}}$. 

Let $\mathcal{C}$ be the set of e-cycles defined as follows: an e-cycle $C$ is in $\mathcal{C}$ iff it lies in $\overrightarrow{E}_* \cup \overleftarrow{E}$ for some $\overleftarrow{E} \in \mathcal{E}_{\mathcal{A}}$. By Lemma~\ref{disjointscycles2}, $\mathcal{C}$ must either contain an e-cycle which fails to be an s-cycle or two es-cycles which fail to be edge-disjoint. In the latter case, by Corollary~\ref{StoR3}, these must have S-to-R intersection. In either case Condition~($*$) is failed. \endproof

\newpage

\section{Examples}

\subsection{A system with three variables}
\label{exnotinteraction}

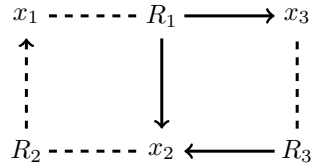
\begin{figure}[h]
\begin{center}
\begin{tikzpicture}[domain=0:4,scale=0.6]

\node at (1,4) {$x_1$};
\node at (4,4) {$R_1$};
\node at (7,4) {$x_3$};

\node at (1,1) {$R_2$};
\node at (4,1) {$x_2$};
\node at (7,1) {$R_3$};

\draw[-, dashed, line width=0.04cm] (1.5, 1) -- (3.5,1);
\draw[-, dashed, line width=0.04cm] (1.5, 4) -- (3.5,4);

\draw[<-, line width=0.04cm] (4.5, 1) -- (6.5,1);
\draw[->, line width=0.04cm] (4.5, 4) -- (6.5,4);

\draw[->, dashed, line width=0.04cm] (1, 1.5) -- (1,3.5);
\draw[<-, line width=0.04cm] (4, 1.5) -- (4,3.5);
\draw[-, dashed, line width=0.04cm] (7, 1.5) -- (7,3.5);

\end{tikzpicture}
\end{center}
\caption{\label{fignotinter}The DSR graph corresponding to the trivial decomposition $J = JI$ for System (\ref{exnotinter}). All edge-labels are $\infty$ and have been omitted. By inspection, the graph contains no e-cycles and hence satisfies Condition~($*$). As System (\ref{exnotinter}) has full outflows, it is thus injective.}
\end{figure}

\nopagebreak
Consider the following dynamical system on $\mathbb{R}^3$:
\begin{equation}
\label{exnotinter}
\begin{array}{ccl}
\dot x_1 & = & f_1(x_1, x_2) - q_1(x_1)\\
\dot x_2 & = & f_2(x_1, x_2, x_3) - q_2(x_2)\\
\dot x_3 & = & f_3(x_1, x_3) - q_3(x_3)\\
\end{array}
\end{equation}
where $\frac{\partial q_i}{\partial x_i} > 0$, and moreover $\frac{\partial f_1}{\partial x_1}, \frac{\partial f_1}{\partial x_2}, \frac{\partial f_2}{\partial x_2}, \frac{\partial f_3}{\partial x_3} \leq 0$, and $\frac{\partial f_2}{\partial x_1}, \frac{\partial f_2}{\partial x_3}, \frac{\partial f_3}{\partial x_1} \geq 0$. Without further information, there is no obvious decomposition of $[f_1,f_2, f_3]^T$ except for trivial decompositions. Defining $v_{ij} = |\partial f_i/\partial x_j|$ we have the Jacobian
\[
J \equiv \left[\begin{array}{rrr}-v_{11} & -v_{12} & 0\\v_{21} & -v_{22} & v_{23}\\ v_{31} & 0 & -v_{33}\end{array}\right] = \left[\begin{array}{rrr}-v_{11} & -v_{12} & 0\\v_{21} & -v_{22} & v_{23}\\ v_{31} & 0 & -v_{33}\end{array}\right]\left[\begin{array}{rrr}1 & 0 & 0\\0 & 1 & 0\\ 0 & 0 & 1\end{array}\right].
\]
Using this decomposition we can construct the DSR graph shown in Figure~\ref{fignotinter} which contains no e-cycles. From this, we deduce that System~\ref{exnotinter} is injective on $\mathbb{R}^3$ (or any rectangular subset of $\mathbb{R}^3$). This example shows that the trivial decomposition can be useful. The relationship between DSR graphs corresponding to trivial decompositions, and the usual interaction graphs \cite{gouze98,soule} is discussed further in the concluding section and forthcoming work.

{\bf Remark.} In this example, if $\frac{\partial f_1}{\partial x_1}, \frac{\partial f_2}{\partial x_2}, \frac{\partial f_3}{\partial x_3} < 0$, then since the DSR graph contains no e-cycles at all, this means that even without outflows the system has $P^{(-)}$ Jacobian forbidding multiple equilibria. We do not pursue this here.

\subsection{A single very simple reaction}
\label{ex1section}

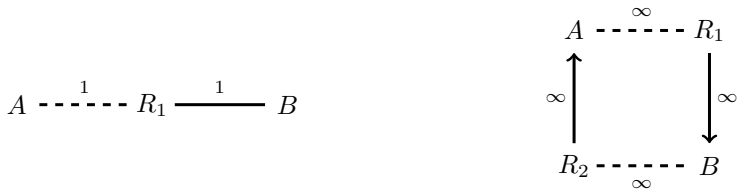
\begin{figure}[h]
\begin{minipage}{\textwidth}
\begin{minipage}{0.45\textwidth}
\hspace{1cm}
\begin{tikzpicture}[domain=0:4,scale=0.6]

\node at (1,1) {$A$};
\node at (4,1) {$R_1$};
\node at (7,1) {$B$};

\draw[-, dashed, line width=0.04cm] (1.5, 1) -- (3.5,1);
\draw[-, line width=0.04cm] (4.5, 1) -- (6.5,1);

\node at (2.5,1.4) {$\scriptstyle{1}$};
\node at (5.5,1.4) {$\scriptstyle{1}$};

\end{tikzpicture}
\end{minipage}
\hfill
\begin{minipage}{0.45\textwidth}
\hspace{1cm}
\begin{tikzpicture}[domain=0:4,scale=0.6]

\node at (1,4) {$A$};
\node at (4,4) {$R_1$};
\node at (1,1) {$R_2$};
\node at (4,1) {$B$};

\draw[-, dashed, line width=0.04cm] (1.5, 1) -- (3.5,1);
\draw[-, dashed, line width=0.04cm] (1.5, 4) -- (3.5,4);
\draw[->, line width=0.04cm] (1, 1.5) -- (1,3.5);
\draw[<-, line width=0.04cm] (4, 1.5) -- (4,3.5);

\node at (2.5,4.4) {$\scriptstyle{\infty}$};
\node at (2.5,0.6) {$\scriptstyle{\infty}$};
\node at (0.6,2.5) {$\scriptstyle{\infty}$};
\node at (4.4,2.5) {$\scriptstyle{\infty}$};

\end{tikzpicture}
\end{minipage}

\end{minipage}
\caption{\label{SRverybasic}{\em Left:} The natural DSR graph for the reaction $A \rightleftharpoons B$ constructed from the decomposition $F = Sv$. The graph contains no cycles. {\em Right:} The DSR graph graph constructed from the decomposition $F = F \circ \mathrm{id}$. The graph contains an e-cycle which fails to be an s-cycle.}
\end{figure}

\nopagebreak
Consider a single chemical reaction involving two substrates $A \rightleftharpoons B$. Denoting the concentration of $A$ by $a$, that of $B$ by $b$, the reaction rate by $v(a, b)$, and ignoring any inflows and outflows, we get the dynamical system
\[
\left[\begin{array}{c}\dot a\\\dot b\end{array}\right] = F(a, b) \equiv \left[\begin{array}{r}-v(a, b)\\v(a, b)\end{array}\right]\,.
\]
Defining $S = [-1, 1]^T$ to be the stoichiometric matrix of the system, gives the natural decomposition $F(a, b) = Sv(a, b)$. However we could also choose to write $F= F\circ \mathrm{id}$, where $\mathrm{id}(\cdot)$ is the identity on $\mathbb{R}^2$. These two choices give us the following two product forms for the Jacobian $DF(a, b)$: 
\[
DF(a, b) = \left[\begin{array}{r}-1\\1\end{array}\right]\begin{array}{c}\left[\begin{array}{cc}v_a & v_b\end{array}\right] \\{}\end{array} \quad \mbox{and}\quad DF(a, b) = \left[\begin{array}{rr}-v_a & -v_b\\v_a & v_b\end{array}\right]\left[\begin{array}{cc}1\frac{}{} & 0 \\0\frac{}{} & 1\end{array}\right]\,,
\]
where $v_a = \frac{\partial v}{\partial a}$ and $v_b = \frac{\partial v}{\partial b}$. With the natural assumption that $v_a \geq 0$ and $v_b \leq 0$, corresponding to the two products are the DSR graphs shown in Figure~\ref{SRverybasic}. 

The natural DSR graph has no cycles and it is immediate that the system is injective provided both substrates are subject to degradation/outflow. On the other hand the DSR graph obtained from the decomposition $F= F\circ \mathrm{id}$ contains an e-cycle which is not an s-cycle, and cannot be used to draw such conclusions. This example illustrates the importance of the choice of decomposition, discussed further in the conclusions.\footnote{It happens, in this case, that if for each fixed value of $(a, b)$ we generate a DSR graph based on the decomposition $F = F \circ \mathrm{id}$, then these all satisfy Condition~($*$) because the e-cycle is in each case an s-cycle. However it is often the case that for one decomposition the DSR graph satisfies Condition~($*$) globally, while for another it fails Condition~($*$), both pointwise and globally.} 

\subsection{Corollary~\ref{tworeac}} 
\label{tworeacexample}

\begin{figure}[h]
\begin{minipage}{0.6\textwidth}
\[
\left[
\begin{array}{cc}
\mathcal{S}_{ji} & 0\\
\mathcal{S}_{li} & \mathcal{S}_{lk}
\end{array}
\right], \qquad
\left[
\begin{array}{cc}
\mathcal{V}_{ij} & \mathcal{V}_{il} \\
x & \mathcal{V}_{kl}
\end{array}
\right],
\]

\end{minipage}
\hfill
\begin{minipage}{0.3\textwidth}
\begin{tikzpicture}[domain=0:4,scale=0.6]

\node at (1,4) {$S_j$};
\node at (4,4) {$R_k$};
\node at (1,1) {$R_i$};
\node at (4,1) {$S_l$};

\draw[->, thick] (1.6,4.2) .. controls (2.2,4.5) and (2.8,4.5) .. (3.4,4.2);
\draw[->, dashed, thick] (1.6,3.8) .. controls (2.2,3.5) and (2.8,3.5) .. (3.4,3.8);

\draw [<-, thick,snake=snake,segment amplitude=.4mm,
         segment length=2mm,line before snake=1mm] (1, 3.5) -- (1,1.5);

\draw [->, thick,snake=snake,segment amplitude=.4mm,
         segment length=2mm,line after snake=1mm] (4, 3.5) -- (4,1.5);

\draw [<-, thick,snake=snake,segment amplitude=.4mm,
         segment length=2mm,line before snake=1mm] (1.5, 1) -- (3.5,1);

\node at (2.5,4.7) {$\scriptstyle{\infty}$};
\node at (2.5,3.3) {$\scriptstyle{\infty}$};

\end{tikzpicture}
\end{minipage}
\caption{\label{SRpart} Submatrices and subgraph of the DSR graph corresponding to the situation described in Corollary~\ref{tworeac}. Wavy edges correspond to edges of unknown sign and label or possibly edge-pairs. The subgraph of the DSR graph has been drawn for $\mathcal{S}_{li}=\mathcal{V}_{ij}=\mathcal{V}_{kl}=0$. Allowing these to take nonzero values simply makes some directed edges undirected. Clearly at least one of the cycles connecting $R_i$, $S_j$, $R_k$ and $S_l$ must be an e-cycle which fails to be an s-cycle as it includes an edge with edge-label $\infty$. }
\end{figure}
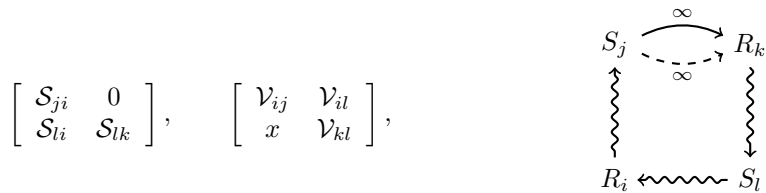

\nopagebreak
Consider the result in Corollary~\ref{tworeac} which gave submatrices and subgraph shown in Figure~\ref{SRpart}. The subgraph has an e-cycle which fails to be an s-cycle. 

The fact that the DSR graph contains an e-cycle that fails to be an s-cycle does not in itself prove that the system must fail the conditions in Lemmas~\ref{genlem}~and~\ref{genlemconv} -- for this the matrix-theoretic formulation is necessary. It is an easy matter to find counterexamples \cite{banajicraciun} illustrating that Condition~($*$) is not necessary to ensure that a system is $P_0^{(-)}$. This example does illustrate however that with certain assumptions, certain motifs in a DSR graph may be {\em sufficient} to ensure that the conditions are failed, a fact which will be explored in future work.

\subsection{A famous interaction network: the ``repressilator''}

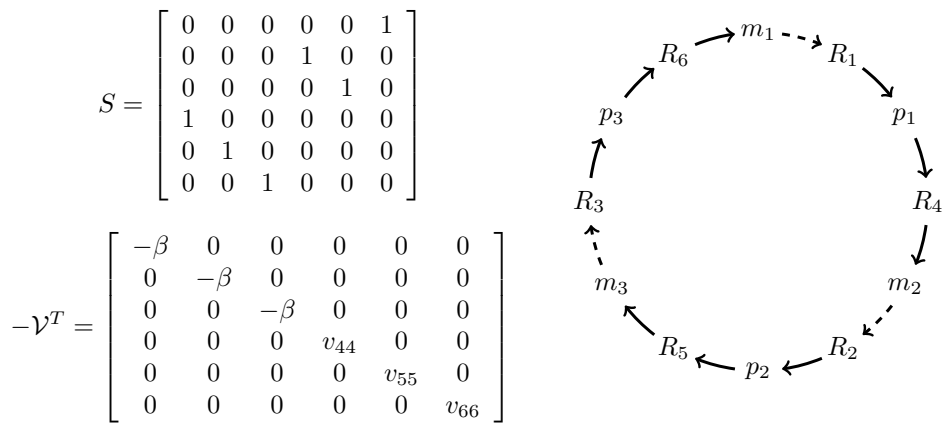
\begin{figure}[h]
\begin{minipage}{\textwidth}
\begin{minipage}{0.58\textwidth}
\[
\begin{array}{c}
S = \left[\begin{array}{cccccc}0&0&0&0&0&1\\ 0&0&0&1&0&0\\ 0&0&0&0&1&0 \\ 1 & 0 & 0 & 0 & 0 & 0\\ 0 & 1 & 0 & 0 & 0 & 0\\ 0 & 0 & 1 & 0 & 0 & 0\end{array}\right]\\\\
-\mathcal{V}^T = \left[\begin{array}{cccccc}-\beta&0&0&0&0&0\\ 0&-\beta&0&0&0&0\\ 0&0&-\beta&0&0&0 \\ 0 & 0 & 0 & v_{44} & 0 & 0\\ 0 & 0 & 0 & 0 & v_{55} & 0\\ 0 & 0 & 0 & 0 & 0 & v_{66}\end{array}\right]
\end{array}
\]

\end{minipage}
\hfill
\begin{minipage}{0.4\textwidth}
\begin{tikzpicture}[domain=-4:4,scale=0.75]

\path (90: 3cm) coordinate(S1);
\path (60: 3cm) coordinate(R1);
\path (30: 3cm) coordinate(S4);
\path (0: 3cm) coordinate(R4);
\path (-30: 3cm) coordinate(S2);
\path (-60: 3cm) coordinate(R2);
\path (-90: 3cm) coordinate(S5);
\path (-120: 3cm) coordinate(R5);
\path (-150: 3cm) coordinate(S3);
\path (180: 3cm) coordinate(R3);
\path (150: 3cm) coordinate(S6);
\path (120: 3cm) coordinate(R6);

\node at (R1) {$R_1$};
\node at (R2) {$R_2$};
\node at (R3) {$R_3$};
\node at (R4) {$R_4$};
\node at (R5) {$R_5$};
\node at (R6) {$R_6$};
\node at (S1) {$m_1$};
\node at (S2) {$m_2$};
\node at (S3) {$m_3$};
\node at (S4) {$p_1$};
\node at (S5) {$p_2$};
\node at (S6) {$p_3$};

\path (90-8: 3cm) coordinate(S1end);
\path (60-8: 3cm) coordinate(R1end);
\path (30-8: 3cm) coordinate(S4end);
\path (-8: 3cm) coordinate(R4end);
\path (-30-8: 3cm) coordinate(S2end);
\path (-60-8: 3cm) coordinate(R2end);
\path (-90-8: 3cm) coordinate(S5end);
\path (-120-8: 3cm) coordinate(R5end);
\path (-150-8: 3cm) coordinate(S3end);
\path (180-8: 3cm) coordinate(R3end);
\path (150-8: 3cm) coordinate(S6end);
\path (120-8: 3cm) coordinate(R6end);

\draw[->, dashed, line width=0.04cm] (S1end)  arc (90-8:90-22:3cm);
\draw[->, line width=0.04cm] (R1end)  arc (60-8:60-22:3cm);
\draw[->, line width=0.04cm] (S4end)  arc (30-8:30-22:3cm);
\draw[->, line width=0.04cm] (R4end)  arc (-8:-22:3cm);
\draw[->, dashed, line width=0.04cm] (S2end)  arc (-30-8:-30-22:3cm);
\draw[->, line width=0.04cm] (R2end)  arc (-60-8:-60-22:3cm);
\draw[->, line width=0.04cm] (S5end)  arc (-90-8:-90-22:3cm);
\draw[->, line width=0.04cm] (R5end)  arc (-120-8:-120-22:3cm);
\draw[->, dashed, line width=0.04cm] (S3end)  arc (-150-8:-150-22:3cm);
\draw[->, line width=0.04cm] (R3end)  arc (180-8:180-22:3cm);
\draw[->, line width=0.04cm] (S6end)  arc (150-8:150-22:3cm);
\draw[->, line width=0.04cm] (R6end)  arc (120-8:120-22:3cm);

\end{tikzpicture}

\end{minipage}
\end{minipage}
\caption{\label{figrepres}The matrices $S$ and $-\mathcal{V}^T$, and the resulting DSR graph for the repressilator. The quantities $v_{ii}$ are all positive. The DSR graph consists of a single o-cycle. Edge-labels have been omitted.}
\end{figure}

\nopagebreak
A famous example of a synthetic oscillator is the {\em repressilator} described in \cite{elowitz}. The system is described by the differential equations:
\[
\begin{array}{c}\dot m_1\\ \dot m_2 \\ \dot m_3 \\ \dot p_1\\ \dot p_2 \\ \dot p_3\end{array} \begin{array}{c}=\\=\\=\\=\\=\\=\end{array} \begin{array}{c}\alpha_0\\ \alpha_0 \\ \alpha_0 \\ {}\\ {}\\ {}\end{array} \begin{array}{c}+\\+ \\+ \\ {}\\ {}\\ {}\end{array} \begin{array}{c}f(p_3)\\ f(p_1) \\ f(p_2)\\\beta m_1\\ \beta m_2 \\ \beta m_3\end{array} \begin{array}{c}-\\-\\-\\-\\-\\-\end{array} \begin{array}{c}m_1\\ m_2 \\ m_3 \\ \beta p_1\\ \beta p_2\\ \beta p_3\end{array}
\] 
where $\alpha_0, \beta > 0$ and $f$ is a decreasing function. The matrices $S$, $-\mathcal{V}^T$ and the DSR graph $G \equiv G_{S, -\mathcal{V}^T}$ for the system are shown in Figure~\ref{figrepres}. With the given outflow conditions the system is injective on any rectangular domain as $G$ consists of an o-cycle. Injectivity is also easily derived using the trivial decomposition, or using results from \cite{gouze98,soule}.

\subsection{Reaction systems with unknown influences} 

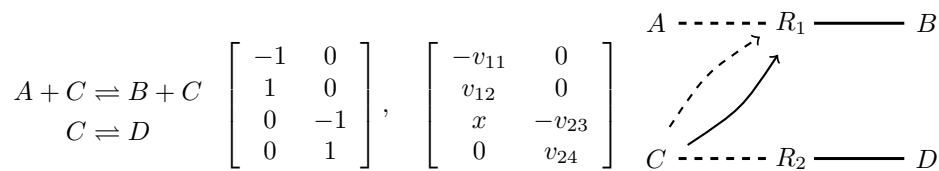
\begin{figure}[h]
\begin{minipage}{\textwidth}
\begin{minipage}{0.2\textwidth}
\begin{eqnarray}
A + C & \rightleftharpoons & B + C \nonumber\\
C  & \rightleftharpoons & D \nonumber
\end{eqnarray}
\end{minipage}
\hfill
\begin{minipage}{0.42\textwidth}
\[
\left[
\begin{array}{cc}
-1 & 0\\
1 & 0\\
0 & -1\\
0 & 1
\end{array}
\right], \quad \left[
\begin{array}{cc}
-v_{11} & 0\\
v_{12} & 0\\
x & -v_{23}\\
0 & v_{24}
\end{array}
\right]
\]
\end{minipage}
\hfill
\begin{minipage}{0.36\textwidth}
\begin{tikzpicture}[domain=0:4,scale=0.6]

\node at (1,4) {$A$};
\node at (4,4) {$R_1$};
\node at (7,4) {$B$};
\node at (1,1) {$C$};
\node at (4,1) {$R_2$};
\node at (7,1) {$D$};

\draw[-, dashed, line width=0.04cm] (1.5, 4) -- (3.5,4);
\draw[-, dashed, line width=0.04cm] (1.5, 1) -- (3.5,1);
\draw[-, line width=0.04cm] (4.5, 4) -- (6.5,4);
\draw[-, line width=0.04cm] (4.5, 1) -- (6.5,1);

\draw[->, dashed, thick] (1.3,1.7) .. controls (2.1,2.9) and (2.1,2.9) .. (3.3,3.7);
\draw[->, thick] (1.7,1.3) .. controls (2.9,2.1) and (2.9,2.1) .. (3.7,3.3);

\end{tikzpicture}
\end{minipage}
\end{minipage}
\caption{\label{SRone}A reaction system with matrices $S$ and $-\mathcal{V}^T$, and the DSR graph. $C$ modulates the first reaction with influence of unknown sign. The quantities $v_{ij}$ are all nonnegative (as expected for any reasonable kinetics), while $x$ is unsigned. Edge-labels have been omitted. There are no cycles in the graph, and so the system is a $P_0^{(-)}$ system.}
\end{figure}

\nopagebreak
We present a reaction network involving two reactions and four substrates which is an example of a system which (with outflow) is injective for all kinetics despite the fact that we don't know whether an influence is activatory or inhibitory. Intuitively, this is because it has hierarchichal structure, decomposing into two subsystems, each well-behaved, and with only one-way influence between the two. The system, with corresponding matrices $S$ and $-\mathcal{V}^T$, and resulting DSR graph, is shown in Figure~\ref{SRone}. As the graph contains no cycles, injectivity of the system with outflows is immediate. 

\subsection{The importance of reversibility}

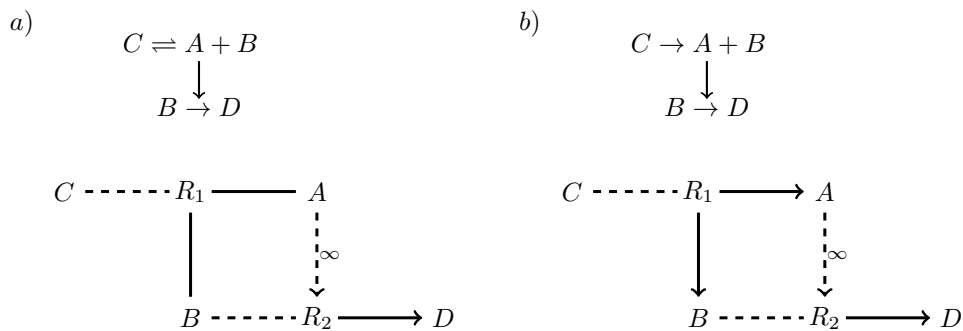
\begin{figure}[h]
\begin{minipage}{\textwidth}
\begin{minipage}{0.48\textwidth}
\begin{tikzpicture}[domain=0:4,scale=0.56]

\node at (0,8) {$a)$};

\node at (4,7.5) {$C  \rightleftharpoons A + B$};
\draw[->, thick] (4.2, 7.1) -- (4.2,6.2);
\node at (4.2,6) {$B \rightarrow D$};

\node at (1,4) {$C$};
\node at (4,4) {$R_1$};
\node at (7,4) {$A$};
\node at (4,1) {$B$};
\node at (7,1) {$R_2$};
\node at (10,1) {$D$};

\draw[-, dashed, line width=0.04cm] (1.5, 4) -- (3.5,4);
\draw[-, line width=0.04cm] (4.5, 4) -- (6.5,4);
\draw[-, dashed, line width=0.04cm] (4.5, 1) -- (6.5,1);
\draw[->, line width=0.04cm] (7.5, 1) -- (9.5,1);

\draw[-, line width=0.04cm] (4, 3.5) -- (4,1.5);
\draw[->, dashed, line width=0.04cm] (7, 3.5) -- (7,1.5);

\node at (7.3, 2.5) {$\scriptstyle{\infty}$};
\end{tikzpicture}
\end{minipage}
\hfill
\begin{minipage}{0.48\textwidth}
\begin{tikzpicture}[domain=0:4,scale=0.56]

\node at (0,8) {$b)$};

\node at (4,7.5) {$C  \rightarrow A + B$};
\draw[->, thick] (4.2, 7.1) -- (4.2,6.2);
\node at (4.2,6) {$B \rightarrow D$};

\node at (1,4) {$C$};
\node at (4,4) {$R_1$};
\node at (7,4) {$A$};
\node at (4,1) {$B$};
\node at (7,1) {$R_2$};
\node at (10,1) {$D$};

\draw[-, dashed, line width=0.04cm] (1.5, 4) -- (3.5,4);
\draw[->, line width=0.04cm] (4.5, 4) -- (6.5,4);
\draw[-, dashed, line width=0.04cm] (4.5, 1) -- (6.5,1);
\draw[->, line width=0.04cm] (7.5, 1) -- (9.5,1);

\draw[->, line width=0.04cm] (4, 3.5) -- (4,1.5);
\draw[->, dashed, line width=0.04cm] (7, 3.5) -- (7,1.5);

\node at (7.3, 2.5) {$\scriptstyle{\infty}$};
\end{tikzpicture}
\end{minipage}
\end{minipage}
\caption{\label{figreversible} Two reaction systems and their DSR graphs. All edges except for those labelled $\infty$ have edge-label $1$. If the first reaction is irreversible, then the system is a $P_0^{(-)}$ system.}
\end{figure}

\nopagebreak
Consider the reaction systems and their respective DSR graphs shown in Figure~\ref{figreversible}, drawn for any kinetics satisfying the N1C condition \cite{banajicraciun}. In each case $A$, a product of the first reaction, activates the second reaction. In Figure~\ref{figreversible}a the first reaction is assumed to be reversible. The associated DSR graph has an e-cycle which fails to be an s-cycle and so cannot be used to make injectivity claims. In Figure~\ref{figreversible}b the first reaction is assumed to be irreversible (in the strong sense that the products of the reaction cannot influence the reaction rate). In this case the DSR graph contains no cycles at all, and so the system with outflows is injective. Intuitively, this is because in the first case $A$ and $B$ may be able to activate their own production, while in the second they cannot.

\subsection{Back to the TCA cycle}

\nopagebreak
We return to the DSR graphs of the TCA cycle model in \cite{cortassa} shown in Figure~\ref{SRTCA}. A couple of points are noteworthy. First, several of the reactions appear as reversible in the DSR graph (in particular the production of fumarate and oxaloacetate) because although the functional forms given in \cite{cortassa} may only permit the reaction to proceed in one direction, the reaction rates are inhibited by the reaction products. Secondly, given a conserved pair such as NAD/NADH, one can write $[\mathrm{NAD}] = \mathrm{NAD_{tot}} - [\mathrm{NADH}]$ for some constant $\mathrm{NAD_{tot}}$ to eliminate $\mathrm{NAD}$ from the system and simplify the DSR graph. Subsequently all claims refer to the system with some fixed value of $\mathrm{NAD_{tot}}$, i.e. on some fixed stoichiometric class. This procedure works well with a conserved pair, but is not necessarily helpful with more general conserved quantities, possibly introducing a number of additional edges into the DSR graph.

At the first stage of construction, shown in Figure~\ref{SRTCA}a, the DSR graph has a single es-cycle and fulfils Condition~($*$). Adding in the inhibition by oxaloacetate of succinate dehydrogenase (Figure~\ref{SRTCA}b) introduces a single o-cycle with undefined stoichiometry. It does {\em not} introduce an e-cycle with undefined stoichiometry, precisely because some of the reactions are irreversible (in the strong sense), indicating again the importance of considerations of reversibility which would be missed by analysis using the SR graph alone. Introducing $\mathrm{NADH}$ adds a single S-vertex and three extra positive edges (with edge-labels $1$) to the DSR graph (Figure~\ref{SRTCA}c). All these are undirected whether or not the reactions are reversible in the chemical sense, because an increase in the level of NADH is equivalent to a drop in the level of NAD, and so affects the reaction rate. This adds a great number of cycles to the graph, but they are all o-cycles. Finally, adding the AAT-catalysed reaction (Figure~\ref{SRTCA}d), leads to the creation of yet more cycles. Now the graph violates Condition~($*$), both because it contains e-cycles which fail to be s-cycles (e.g. OAA-$R_6$-FUM-$R_7$-MAL-$R_8$-NADH-$R_3$-$\alpha\mathrm{KG}$-$R_9$-OAA) and because it contains es-cycles with S-to-R intersection. 

Note that even at this stage some model quantities, in particular ADP, calcium, and membrane potential are omitted from the analysis. This example illustrates that even in models where we cannot positively rule out multiple equilibria, by identifying cycles which cause Condition~($*$) to be violated, we can speculate on the mechanisms by which multistationarity may arise.

\section{Discussion and Conclusions}
\label{secfinal}

\nopagebreak
We have shown that the DSR graph, an object closely related to the pathway diagrams and interaction diagrams drawn by applied scientists in various fields, can be used to rule out the possibility of multiple equilibria or multiple nondegenerate equilibria. Given a systems written in one of the forms 
\[
\dot x = f(v(x)) - Q(x) \quad \mbox{or} \quad \dot x = f(v(x)) - Q_\theta(x)\,,
\]
we construct the DSR graph for the associated systems $\dot x = f(v(x))$, and check whether Condition~($*$) holds for this graph. If it does then we know that all of $f(v(x)) - Q(x)$ are injective and that all of $\dot x = f(v(x)) - Q_\theta(x)$ forbid multiple nondegenerate equilibria on the relative interior of any invariant affine set. 

Algorithmic development is clearly important at this stage. Simple, widely available, computational tools to test whether Condition~($*$) holds for a DSR graph will open up the possibility of routine preliminary model analysis: modellers working with qualitative models will be able to check whether their models, as a result of structure alone, forbid multiple equilibria, before parametrisation and simulation. We are currently working on creating user-friendly software that performs this test.

It is worth remembering that matrix-theoretic approaches give sharper results than the graph-theoretic ones. As mentioned in the text, Lemmas~\ref{genlem}~and~\ref{genlemconv} simplify in a number of special cases where we have additional knowledge of a system (e.g. reversibility of all processes). This specialisation is a task for future work. 

A number of special cases can be treated via minor variants on this analysis. Some systems on the nonnegative orthant, for example population models, generally have boundary equilibria, and one question is whether they admit multiple positive equilibria. The general theory developed here can be applied in this context, for example to a Lotka-Volterra type model for the evolution of $n$ species of the form $\dot x_i = x_i G_i(x)$ \cite{HRThieme,XQZhao} where $x = [x_1, \ldots, x_n]^T$ is the nonnegative vector of species concentrations, and the functions $G_i(x)$ encode information on the interactions between species. Positive equilibria must satisfy $G_i(x) = 0$, and hence we are interested in injectivity of $G(x) = [G_1(x), \ldots, G_n(x)]^T$ in the interior of the positive orthant, which can be treated using the techniques developed here. 

The example in Section~\ref{exnotinteraction} illustrated that even trivial decompositions can be useful. In fact, in forthcoming work, we will show how given an arbitrary dynamical system $\dot x = F(x)$ on a rectangular domain, and using the trivial decomposition $F = F\circ \mathrm{id}$, it is possible to obtain stronger results than the well-known interaction graph results on injectivity (results in \cite{soule} for example). However, the example in Section~\ref{ex1section} illustrated the important and subtle point that in some cases the strongest results are obtained by choosing the ``correct'' decomposition of the functions we hope to show are injective. This is true for both the matrix-theoretic and the graph-theoretic results. The following difficult question remains open: given an arbitrary set of functions, how do we choose the decomposition which will give us the sharpest results on injectivity (assuming that such a choice exists)? 

Finally, the importance of cycle structure in SR and DSR graphs goes beyond questions of injectivity. Cycle structure in SR graphs is, for example, linked closely to monotonicity \cite{angelileenheersontag} (several results in \cite{banajidynsys} also have immediate graph-theoretic interpretations). The close relationship between cycle structure and the possibility of complex behaviour has been understood for some time for models whose structure can be represented via an interaction graph \cite{kaufman}, but the corresponding results for general interaction networks represented via SR/DSR graphs (or possibly other variants of these ideas) are less complete. Such study has important implications for model caricature and model simplification.

\section*{Acknowledgements}
GC acknowledges support from NIH grant R01GM86881. MB acknowledges support from EPSRC grant EP/D060982/1. We would like to thank Stefan M\"uller for helpful comments on our treatment of the repressilator, and the reviewers of this paper for a number of useful suggestions.

\appendix

\section{Definitions and notation}
\label{app1}

\nopagebreak
In each of the the definitions below, $M$ is an $n \times m$ matrix, $\delta \subset \{1, \ldots, n\}$, and $\gamma \subset \{1, \ldots, m\}$ with $|\delta| = |\gamma|$.

$M(\delta|\gamma)$ is the submatrix of $M$ with rows indexed by $\delta$ and columns indexed by $\gamma$. A {\bf principal submatrix} of $M$ is a submatrix of the form $M(\delta|\delta)$. Determinants of submatrices of $M$ are termed {\bf minors} of $M$. {\bf Principal minors} are determinants of principal submatrices. 

{\bf $P$ matrices} are square matrices all of whose principal minors are positive. They are by definition nonsingular. If $\,-M$ is a $P$ matrix, then $M$ is a {\bf $P^{(-)}$ matrix}. Each $k \times k$ principal minor of a $P^{(-)}$ matrix has sign $(-1)^k$.

{\bf $P_0$ matrices} (so termed in \cite{hershkowitz}) are matrices in the closure of the set of $P$ matrices. These are matrices all of whose principal minors are nonnegative. We will term $M$ a {\bf $P_0^{(-)}$ matrix} if $-M$ is a $P_0$ matrix. A $k \times k$ minor of a $P_0^{(-)}$ matrix is either zero or has sign $(-1)^k$. The zero matrix is a $P_0$ and $P_0^{(-)}$ matrix.

{\bf Qualitative classes.} A matrix $M$ determines the qualitative class $\mathcal{Q}(M)$ \cite{brualdi} consisting of all matrices with the same sign pattern as $M$. Explicitly, $\mathcal{Q}(M)$ consists of all matrices $X$ with the same dimensions as $M$, and satisfying $M_{ij} > 0 \Rightarrow X_{ij} > 0$, $M_{ij} < 0 \Rightarrow X_{ij} < 0$ and $M_{ij} = 0 \Rightarrow X_{ij} = 0$. 

$\mathcal{Q}_0(M)$ is the closure of $\mathcal{Q}(M)$. Explicitly, $\mathcal{Q}_0(M)$ consists of all matrices $X$ with the same dimensions as $M$ such that $M_{ij} > 0 \Rightarrow X_{ij} \geq 0$, $M_{ij} < 0 \Rightarrow X_{ij} \leq 0$ and $M_{ij} = 0 \Rightarrow X_{ij} = 0$. 

A square matrix $M$ is {\bf sign nonsingular} if all matrices in $\mathcal{Q}(M)$ are nonsingular. A matrix-set $\mathcal{M}$ has {\bf signed determinant} if either $\mathrm{det}(\mathcal{M}) = 0$, or $\mathrm{det}(\mathcal{M}) \subset (-\infty, 0)$ or $\mathrm{det}(\mathcal{M}) \subset (0, \infty)$.

\section{Multiple nondegenerate equilibria in systems with conserved quantities}
\label{apppersistence}

\nopagebreak
The discussion in this appendix is closely related both to Theorem~2 in \cite{craciun2} and to Proposition~1 in \cite{soule}. 

Consider a dynamical system $\dot x = f(x)$ defined on $X \subset \mathbb{R}^n$. An equilibrium $p$ of the system is {\bf nondegenerate} if $Df(p)$ (the Jacobian at $p$) has no zero eigenvalues. When $X$ is foliated by invariant sets we are generally interested in behaviour on one of these invariant sets. Assume that $\dot x = f(x)$ preserves some $C^1$ function $E: \mathbb{R}^n \to \mathbb{R}^l$ ($1 \leq l < n$), so that $\dot E(x) = 0$ along trajectories. Since any level set of the form $E_C \equiv \{x\,|\, E(x) = C\}$ in invariant, the pertinent question is whether the system restricted to $E_C \cap X$ admits multiple nondegenerate equilibria, i.e. whether the system $\dot x = f(x)$ can have two equilibria $p, q \in (E_C \cap X)$ such that $Df(p)$ and $Df(q)$ have no eigenvectors corresponding to zero eigenvalues, tangent to $E_C$. 

Assume that $X$ is a rectangular subset of $\mathbb{R}^n$. We already know that if (\ref{genprodform1}) is a $P_0^{(-)}$ system, then (\ref{genprodformout1}) is injective on $X$. We now show that if (\ref{genprodformout1}) is injective on $X$, then (\ref{genprodformconserved}) is incapable of multiple nondegenerate equilibria in the following sense: let $E_C$ be any invariant affine subspace of $\mathbb{R}^n$. Then the system restricted to $E_C \cap X$ can have no multiple nondegenerate equilibria in $\mathrm{ri}(E_C \cap X)$, the relative interior of $E_C \cap X$.

In order to show this, we need the following basic persistence property of nondegenerate equilibria: if a $C^1$ vector field $f$ on some smooth manifold $\mathcal{S}$ (possibly with boundary) has a nondegenerate equilibrium $p \in \mathrm{int}(\mathcal{S})$ then any vector field $g$ on $\mathcal{S}$, close to $f$ in the $C^1$ topology, has an equilibrium $p^{'}$, close to $p$, in $\mathrm{int}(\mathcal{S})$. This follows, for example, by isolating $p$ in a sufficiently small closed neighbourhood $U_p \subset \mathrm{int}(\mathcal{S})$ and using the invariance of the Brouwer degree in $\mathrm{int}(U_p)$ under small perturbations of $f$.

\begin{lem}
\label{lemconserved}
Consider a dynamical system $\dot x = f(x)$ on $X \subset \mathbb{R}^n$, where $X$ is a convex, forward invariant set. Assume that
\begin{enumerate}
\item Given any $\lambda > 0$, all functions of the form $f(x) - \lambda x$ are injective on $X$.
\item There is an affine subspace $E_C \subset \mathbb{R}^n$ such that $E_C \cap X$ is invariant for $\dot x = f(x)$. 
\end{enumerate}
Then the system $\dot x = f(x)$ restricted to $E_C \cap X$ cannot have more than one nondegenerate equilibrium in $\mathrm{ri}(E_C \cap X)$. 
\end{lem}
\begin{proof}
The result is trivial if $E_C \cap X$ is empty or if $E_C$ is a point. Otherwise, either $E_C = \mathbb{R}^n$, in which case $E_C \cap X = X$, or $E_C$ has dimension $n - l$ for some $1 \leq l \leq n-1$, and there is a surjective linear function $E: \mathbb{R}^n \to \mathbb{R}^l$, and a vector $C \in \mathbb{R}^l$ such that $E_C \equiv \{x\,|\,Ex = C\}$, and for $x \in E_C$, $Ef(x)=0$. In either case, choose any vector $k \in E_C \cap X$. Define the systems $\dot x = f(x) + \lambda(k-x) \equiv F_\lambda(x)$ with $\lambda > 0$. By convexity of $X$ and the fact that $X$ is forward invariant for $\dot x = f(x)$, we have that $X$ is forward invariant for $\dot x = F_\lambda(x)$. Further, $\dot x = F_\lambda(x)$ leaves $E_C \cap X$ invariant: this is obvious when $E_C = \mathbb{R}^n$; when $E_C \not = \mathbb{R}^n$, we have $Ef(x) + \lambda(Ek-Ex) = 0$ when $x \in E_C$. By assumption $f(x) - \lambda x$, and hence $F_\lambda(x)$, are injective on $X$. Now suppose that $\dot x = f(x)$ contains two nondegenerate equilibria in $\mathrm{ri}(E_C \cap X)$. Since, for small $\lambda$, $F_\lambda(x)$ is $C^1$ close to $f(x)$, by the above persistence arguments, the system $\dot x = F_\lambda(x)$ must contain two (nondegenerate) equilibria on $\mathrm{ri}(E_C \cap X)$ contradicting the fact that $F_\lambda$ are injective. 
\end{proof}

{\bf Remark.} In Lemma~\ref{lemconserved}, rather than choosing equilibria, i.e. values of $x$ satisfying $f(x)=0$, we can consider values of $x$ satisfying $f(x) = c$ for any $c$. Theorem~4w in \cite{gale} implies that if $X$ is an open rectangular subset of $\mathbb{R}^n$, and, for some $\theta$, (\ref{genprodformconserved}) is a $P_0^{(-)}$ system with nonsingular Jacobian everywhere in $X$, then it is injective on the interior of $X$. This result can be seen as a special case of Lemma~\ref{lemconserved}.

We remark that extensions of Lemma~\ref{lemconserved} to allow nonlinear integrals $E$ are also possible -- the key question is whether given some fixed level set $E_C$ we can construct a family of $C^1$ systems $\dot x = F_\lambda(x)$ which preserve $E_C$, are injective for $\lambda > 0$, and such that $F_0(x) = f(x)$. We do not pursue this here.

\section{Further discussion of DSR graphs}
\label{appDSR}

\begin{figure}[h]
\begin{center}
\begin{tikzpicture}[domain=0:4,scale=0.6]

\node at (0,5.5) {$a)$};
\node at (1,5) {$S$};
\node at (5,5) {$R$};
\draw[->, line width=0.04cm] (1.7,5) -- (4.3,5);

\node at (7,5.5) {$b)$};
\node at (8,5) {$S$};
\node at (12,5) {$R$};
\draw[<-, line width=0.04cm] (8.7,5) -- (11.3,5);

\node at (14,5.5) {$c)$};
\node at (15,5) {$S$};
\node at (19,5) {$R$};
\draw[-, line width=0.04cm] (15.7,5) -- (18.3,5);

\node at (0,3.5) {$d)$};
\node at (1,3) {$S$};
\node at (5,3) {$R$};
\draw[dashed, ->, line width=0.04cm] (1.7,3) -- (4.3,3);

\node at (7,3.5) {$e)$};
\node at (8,3) {$S$};
\node at (12,3) {$R$};
\draw[dashed, <-, line width=0.04cm] (8.7,3) -- (11.3,3);

\node at (14,3.5) {$f)$};
\node at (15,3) {$S$};
\node at (19,3) {$R$};
\draw[dashed, -, line width=0.04cm] (15.7,3) -- (18.3,3);

\end{tikzpicture}
\end{center}
\caption{\label{motifs1} All single edge connections between an S-vertex and an R-vertex in a DSR graph. a) Increases in $S$ decrease $R$ ($S$ inhibits the interaction), but $S$ is unaffected by the interaction. b) The interaction increases $S$, but its rate is not affected by $S$. c) Increases in $S$ decrease $R$ and the interaction increases $S$. d) Increases in $S$ increase $R$ ($S$ activates the interaction), but $S$ is unaffected by the interaction. e) The interaction decreases $S$, but its rate is not affected by $S$. f) Increases in $S$ increase $R$ and the interaction decreases $S$.}
\end{figure}
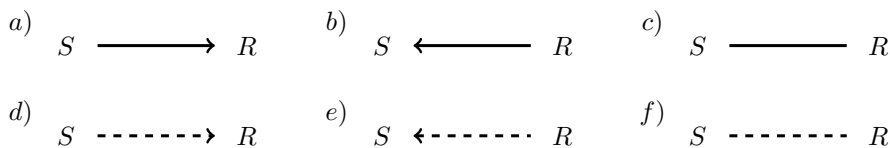

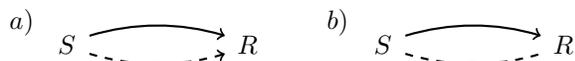
\begin{figure}[h]
\begin{center}
\begin{tikzpicture}[domain=0:4,scale=0.6]

\node at (0,1.5) {$a)$};
\node at (1,1) {$S$};
\node at (5,1) {$R$};
\draw[->, thick] (1.5,1.2) .. controls (2.5,1.5) and (3.5,1.5) .. (4.5,1.2);
\draw[->, dashed, thick] (1.5,0.8) .. controls (2.5,0.5) and (3.5,0.5) .. (4.5,0.8);

\node at (7,1.5) {$b)$};
\node at (8,1) {$S$};
\node at (12,1) {$R$};
\draw[->, thick] (8.5,1.2) .. controls (9.5,1.5) and (10.5,1.5) .. (11.5,1.2);
\draw[-, dashed, thick] (8.5,0.8) .. controls (9.5,0.5) and (10.5,0.5) .. (11.5,0.8);

\end{tikzpicture}
\end{center}
\caption{\label{motifs2} Two common double edge connections between vertices in a DSR graph. a) Increases in $S$ may cause either an increase or a decrease in $R$, but $S$ itself is unaffected by the interaction. b) Increases in $S$ may cause either an increase or a decrease in $R$, and $S$ is decreased by the interaction.}
\end{figure}

\nopagebreak
DSR graphs can be understood via possible {\em motifs}. Consider some species with concentration $S$ which participates in some interaction with rate $R$. In Figure~\ref{motifs1} all six possible single-edge connections between an S-vertex and an R-vertex are illustrated. There are nine possible double edge connections between two vertices: These are created by taking any of the motifs in Figure~\ref{motifs1}a,~b,~or~c, and combining with any motif from Figure~\ref{motifs1}d,~e,~or~f. Two common ones are shown in Figure~\ref{motifs2}, and their meanings are described in the caption.

The DSR graph is a natural amalgamation of directed versions of two SR graphs $\overleftarrow{G}_{\mathcal{S}}$ and $\overrightarrow{G}_{\mathcal{-V^T}}$, associated with two sets of matrices, $\mathcal{S}$ and $\mathcal{V}$. All the matrix-theoretic results treat $\mathcal{S}$ and $\mathcal{V}$ symmetrically, and ideally this would be reflected in the construction of the DSR graph. However, this would involve introducing two sets of edge-labels onto the graph. In our construction, for simplicity, only one set of edge-labels is imported from $G_{\mathcal{S}}$, based on the practical fact that we often expect $\mathcal{S}$, rather than $\mathcal{V}$, to have constant entries and hence nontrivial edge-labels. However, it is straightforward to extend the treatment, allowing two sets of edge-labels, if required in applications. 

For a CRN, the DSR graph encodes information about irreversibility of reactions, and modulation by quantities which do not formally participate in the reaction. Formally, it is the usual SR graph with the following modifications: some edges have become directed, some short cycles have been replaced with single directed edges of defined sign, but edge-label $\infty$, and some edge-labels on short cycles have been changed. If an o-cycle (resp. e-cycle, resp. s-cycle) survives, then it remains an o-cycle (resp. e-cycle, resp. s-cycle). It is clear that results for CRNs obtained using DSR graphs are sharper than those using SR graphs.

\section{Generalisations of results from \cite{banajicraciun}}
\label{app2}

\nopagebreak
The results here are self-contained, but \cite{banajicraciun} contains a more detailed discussion of the relationship between permutations of ordered sets and cycles in SR graphs.

\begin{lem}
\label{permsigns}
Consider a permutation $\alpha$ of some ordered set. Let $\alpha$ be written as the product of disjoint cycles from some set $\mathcal{C}$. Let $\theta = \cup_{c \in \mathcal{C}}c$. Then $P(\alpha)$, the parity of $\alpha$, is given by 
\[
P(\alpha) = (-1)^{\left|\theta\right| - \left|\mathcal{C}\right|}.
\]
\end{lem}
\begin{proof}
This follows by writing any permutation as the product of disjoint cycles and noting that a $k$-cycle is an even permutation if $k$ is odd and vice versa.
\end{proof}

The following two lemmas are key results on the relationship between the signs of two signed subterms and the cycle structure of the union of the corresponding signed term subgraphs, first for SR graphs, and then for DSR graphs. The proof of the first result is presented in full for completeness. 
\begin{lem}
\label{prodformula}
Consider any set of $k \times k$ matrices $\mathcal{A}$ with corresponding SR graph $G_\mathcal{A}$. Consider any two signed subterms $T_{\alpha}^{(r)}$ and $T_{\beta}^{(s)}$ in the determinant expansion of $\mathcal{A}$, corresponding to permutations $\alpha$ and $\beta$ of $\{1, \ldots, k\}$ and signed term subgraphs $E_\alpha^{(r)}$ and $E_\beta^{(s)}$ in $G_\mathcal{A}$. Then
\begin{equation}
\label{prodform}
T_\alpha^{(r)} T_\beta^{(s)} \subset (-1)^{|\mathcal{C}_e|}\mathbb{R}_{\geq 0}
\end{equation}
where $|\mathcal{C}_e|$ is the number of e-cycles in $E_\alpha^{(r)} \cup E_\beta^{(s)}$.
\end{lem}

\begin{proof}
By definition
\[
T_\alpha^{(r)} T_\beta^{(s)} = P(\alpha)P(\beta)\prod_{i =1}^k T_{i \alpha_i}^{(r)}T_{i \beta_i}^{(s)}\,.
\]
Define 
\[
Z \equiv P(\alpha)P(\beta)\prod_{i =1}^k \mathrm{sign}(E_{i \alpha_i}^{(r)})\mathrm{sign}(E_{i \beta_i}^{(s)})\,.
\]
Proving the lemma is equivalent to proving that $Z = (-1)^{|\mathcal{C}_e|}$. 

Let $\theta$ be the set of indices for which $E_{i \alpha_i}^{(r)}$ and $E_{i \beta_i}^{(s)}$ are distinct edges in $G_{\mathcal{A}}$. When $i \in \{1, \ldots, k\}\backslash\theta$, $\mathrm{sign}(E_{i \alpha_i}^{(r)})\mathrm{sign}(E_{i \beta_i}^{(s)}) = 1$. So
\[
Z = P(\alpha)P(\beta)\,\,\prod_{i \in \theta} \mathrm{sign}(E_{i \alpha_i}^{(r)})\mathrm{sign}(E_{i \beta_i}^{(s)}).
\]

Now the edge set 
\[
\bigcup_{i \in \theta}\left(E_{i \alpha_i}^{(r)} \cup E_{i \beta_i}^{(s)}\right)
\]
consists precisely of the set of cycles in $E_\alpha^{(r)} \cup E_\beta^{(s)}$. These cycles are vertex-disjoint since no more than two edges from this set can be incident on any vertex. Let the set of o-cycles in this set be $\mathcal{C}_o$ and the set of e-cycles be $\mathcal{C}_e$, with  $\mathcal{C} =  \mathcal{C}_o \cup \mathcal{C}_e$. Associate with each cycle $c \in \mathcal{C}_o \cup \mathcal{C}_e$ the corresponding index set $\tilde{c}$, i.e., $i \in \tilde{c} \Leftrightarrow E_{i\alpha_i}^{(r)}, E_{i\beta_i}^{(s)} \in c$. Thus corresponding to the sets $ \mathcal{C}_o$ and $ \mathcal{C}_e$ are the sets of index sets $\tilde{\mathcal{C}_o}$ and $\tilde{\mathcal{C}_e}$. Since any two cycles are edge-disjoint, $\tilde{\mathcal{C}_o} \cup \tilde{\mathcal{C}_e}$ is a partition of $\theta$, and we can define
\[
\theta_o \equiv \bigcup_{\tilde c \in \tilde{\mathcal{C}_o}}\tilde c,\quad\theta_e \equiv \bigcup_{\tilde c \in \tilde{\mathcal{C}_e}}\tilde c\quad\mbox{with}\quad |\theta_o| = \sum_{\tilde{c} \in \tilde{\mathcal{C}_o}}|\tilde{c}|,\quad |\theta_e| = \sum_{\tilde{c} \in \tilde{\mathcal{C}_e}}|\tilde{c}|\,.
\] 
Clearly $\theta = \theta_o \cup \theta_e$. We can write
\begin{eqnarray}
\prod_{i \in \theta} T_{i \alpha_i}^{(r)}T_{i \beta_i}^{(s)} & = & \left(\prod_{i \in \theta_o} T_{i \alpha_i}^{(r)}T_{i \beta_i}^{(s)}\right)\left(\prod_{i \in \theta_e} T_{i \alpha_i}^{(r)}T_{i \beta_i}^{(s)}\right) \nonumber\\
&  = &\left(\prod_{\tilde{c} \in \tilde{\mathcal{C}_o}}\prod_{i \in \tilde{c}} T_{i \alpha_i}^{(r)}T_{i \beta_i}^{(s)}\right)\left(\prod_{\tilde{c} \in \tilde{\mathcal{C}_e}}\prod_{i \in \tilde{c}} T_{i \alpha_i}^{(r)}T_{i \beta_i}^{(s)}\right) \,.\nonumber
\end{eqnarray}

So
\begin{eqnarray}
\hspace{-0.7cm}Z & = & P(\alpha)P(\beta)\left(\prod_{\tilde{c} \in \tilde{\mathcal{C}_o}}\prod_{i \in \tilde{c}} \mathrm{sign}(E_{i \alpha_i}^{(r)})\mathrm{sign}(E_{i \beta_i}^{(s)})\right)\left(\prod_{\tilde{c} \in \tilde{\mathcal{C}_e}}\prod_{i \in \tilde{c}} \mathrm{sign}(E_{i \alpha_i}^{(r)})\mathrm{sign}(E_{i \beta_i}^{(s)})\right) \nonumber \\
& = & P(\alpha)P(\beta)\left(\prod_{c \in \mathcal{C}_o}(-1)^{|c|-1}\right)\left(\prod_{c \in \mathcal{C}_e}(-1)^{|c|}\right) \nonumber\\
& = & P(\alpha)P(\beta)(-1)^{|\theta_o| + |\theta_e|- |\mathcal{C}_o|} \nonumber\\
& = & P(\alpha)P(\beta)(-1)^{|\theta|- |\mathcal{C}_o|} \,. \nonumber
\end{eqnarray}
Applying Lemma~\ref{permsigns} to $\beta\circ\alpha^{-1}$ gives us that 
\[
P(\alpha)P(\beta) = P(\beta\circ\alpha^{-1}) = (-1)^{p-q}\,,
\]
where $q$ is the number of cycles in $\beta\circ\alpha^{-1}$ and $p$ is half the number of elements in these cycles. Now there is a one-to-one correspondence between nontrivial cycles in $\beta\circ\alpha^{-1}$ and long cycles in $E_\alpha^{(r)} \cup E_\beta^{(s)}$, with a cycle of length $l$ in  $\beta\circ\alpha^{-1}$ corresponding to a cycle of length $2l$ in $E_\alpha^{(r)} \cup E_\beta^{(s)}$. Let $q^{'}$ be the number of short cycles in $E_\alpha^{(r)} \cup E_\beta^{(s)}$, so that $q + q^{'} = |\mathcal{C}|$. Since there are precisely two edges in a short cycle, $q^{'}$ is also half the number of edges in short cycles in $E_\alpha^{(r)} \cup E_\beta^{(s)}$, so that $p + q^{'} = |\theta|$. This gives us that $p-q = \left|\theta\right| - \left|\mathcal{C}\right|$, so that:
\[
P(\alpha)P(\beta) = (-1)^{\left|\theta\right| - \left|\mathcal{C}\right|}\,.
\]
Completing the argument:
\[
Z = (-1)^{\left|\theta\right| - \left|\mathcal{C}\right|}(-1)^{|\theta|- |\mathcal{C}_o|} = (-1)^{2|\theta| - |\mathcal{C}| - |\mathcal{C}_o|} = (-1)^{|\mathcal{C}| + |\mathcal{C}_o|} = (-1)^{|\mathcal{C}_e|}\,.
\]
This proves the result. 
\end{proof}

In the above result, it is perfectly possible to have $\alpha = \beta$ and/or $r=s$. Further the result is independent of edge-labels, so some of these may be $\infty$.

For DSR graphs this result becomes:
\begin{lem}
\label{prodformula1}
Consider two sets of $k \times k$ matrices $\mathcal{A}$ and $\mathcal{B}$, and the associated DSR graph $G$. Consider any two signed subterms $\overleftarrow{T}_{\alpha}^{(r)}$ in $\mathrm{det}(\mathcal{A})$ and $\overrightarrow{T}_{\beta}^{(s)}$ in $\mathrm{det}(\mathcal{B})$, corresponding to permutations $\alpha$ and $\beta$ of $\{1, \ldots, k\}$ and corresponding to oppositely directed, signed term subgraphs $\overleftarrow{E}_\alpha^{(r)}$ and $\overrightarrow{E}_\beta^{(s)}$ in $G$. Then
\begin{equation}
\label{prodform1}
\overleftarrow{T}_\alpha^{(r)} \overrightarrow{T}_\beta^{(s)} \subset (-1)^{|\mathcal{C}_e|}\mathbb{R}_{\geq 0}
\end{equation}
where $|\mathcal{C}_e|$ is the number of e-cycles in $\overleftarrow{E}_\alpha^{(r)} \cup \overrightarrow{E}_\beta^{(s)}$.
\end{lem}

\begin{proof}
Note first of all that by Lemma~\ref{formaltogenuine}, all formal cycles in $\overleftarrow{E}_\alpha^{(r)} \cup \overrightarrow{E}_\beta^{(s)}$ are genuine cycles and are vertex-disjoint. From here the proof proceeds identically to the result for SR graphs. 
\end{proof}

We now prove that if two signed term subgraphs have only o-cycles in their intersection, then the corresponding signed subterms have the same sign. First for SR graphs:
\begin{lem}
\label{allocycles}
Consider any set of $n \times n$ matrices $\mathcal{A}$. Let $T_\alpha^{(r)}$ and $T_\beta^{(s)}$ be signed subterms in the expansion of $\mathrm{det}(\mathcal{A})$. If all cycles in $E_\alpha^{(r)}\cup E_\beta^{(s)}$ are o-cycles then $T_\alpha^{(r)}T_\beta^{(s)} \subset \mathbb{R}_{\geq 0}$.
\end{lem}
\begin{proof}
If all cycles in $E_\alpha^{(r)}\cup E_\beta^{(s)}$ are o-cycles, then $|\mathcal{C}_e| = 0$, and applying Eq.~(\ref{prodform}) immediately gives $T_\alpha^{(r)} T_\beta^{(s)}  \subset \mathbb{R}_{\geq 0}$.
\end{proof}

The result for DSR graphs was stated as Lemma~\ref{allocycles1}.

{\it Proof of Lemma~\ref{allocycles1}.} If all cycles in $\overleftarrow{E}_\alpha^{(r)}\cup \overrightarrow{E}_\beta^{(s)}$ are o-cycles, then $|\mathcal{C}_e| = 0$, and applying Eq.~(\ref{prodform1}) immediately gives $\overleftarrow{T}_\alpha^{(r)} \overrightarrow{T}_\beta^{(s)}  \subset \mathbb{R}_{\geq 0}$. \endproof

Finally, we prove Lemma~\ref{oneeands} which stated that when two signed term subgraphs of an SR graph contain in their union a single cycle, and this is an es-cycle, then the corresponding terms sum to zero.

{\it Proof of Lemma~\ref{oneeands}.} Note that the fact that $C$ is an es-cycle implies that none of the edges in $C$ have edge-label $\infty$, even if some of the other edges in $G$ may have edge-label $\infty$. Thus $C$ is a long cycle. By definition
\[
T_\alpha^{(r)} + T_\beta^{(s)} = P(\alpha)\prod_{i =1}^k T_{i \alpha_i}^{(r)} + P(\beta)\prod_{i =1}^k T_{i \beta_i}^{(s)}\,.
\]
As usual, let $\theta$ be the set of indices for which $E_{i \alpha_i}^{(r)}$ and $E_{i \beta_i}^{(s)}$ are distinct edges in $G$, so that $\{E_{i \alpha_i}^{(r)}\}_{i \in \theta}$ and $\{E_{i \beta_i}^{(s)}\}_{i \in \theta}$ are precisely the edges in the unique es-cycle $C$. Defining $C_1 = \{E_{i \alpha_i}^{(r)}\}_{i \in \theta}$ and $C_2 = \{E_{i \beta_i}^{(s)}\}_{i \in \theta}$ gives us a disconnecting partition of $C$. Since $C$ is an s-cycle, $\mathrm{val}(C_1)$ and $\mathrm{val}(C_2)$ are defined and equal. Define $Z \equiv \prod_{i \in \{1, \ldots, k\}\backslash\theta} T_{i \alpha_i}^{(r)}$. We can write
\begin{eqnarray}
T_\alpha^{(r)} + T_\beta^{(s)} &= & Z\left(P(\alpha)\prod_{i \in \theta} T_{i \alpha_i}^{(r)} + P(\beta)\prod_{i \in \theta} T_{i \beta_i}^{(s)}\right) \nonumber \\
& = & P(\alpha)Z\left(\prod_{i \in \theta} T_{i \alpha_i}^{(r)} + P(\beta\circ\alpha^{-1})\prod_{i \in \theta} T_{i \beta_i}^{(s)}\right) \nonumber \\
& = & P(\alpha)Z\left(\mathrm{sign}(C_1)\mathrm{val}(C_1) + P(\beta\circ\alpha^{-1})\mathrm{sign}(C_2)\mathrm{val}(C_2)\right)\,. \nonumber 
\end{eqnarray}
$\beta \circ \alpha^{-1}$ can be written as a single cycle of length $|\theta|$, and so from Lemma~\ref{permsigns}, $P(\beta \circ \alpha^{-1}) = (-1)^{|\theta|-1}$. I.e.,
\[
T_\alpha^{(r)} + T_\beta^{(s)} = P(\alpha)Z\left(\mathrm{sign}(C_1)\mathrm{val}(C_1) + (-1)^{|\theta|-1}\mathrm{sign}(C_2)\mathrm{val}(C_2)\right)\,.
\]
Since $C$ is an e-cycle:
\[
\mathrm{sign}(C_2)/\mathrm{sign}(C_1) = \mathrm{sign}(C_1)\mathrm{sign}(C_2) = \mathrm{sign}(C) = (-1)^{|\theta|}.
\] 
Substituting into the expression for $T_\alpha^{(r)} + T_\beta^{(s)}$ gives:
\[
T_\alpha^{(r)} + T_\beta^{(s)} = P(\alpha)\,Z\,\mathrm{sign}(C_1)\left(\mathrm{val}(C_1) - \mathrm{val}(C_2)\right)\,.
\]
Since $C$ is an s-cycle, $\mathrm{val}(C_1) - \mathrm{val}(C_2)=0$, giving $T_\alpha^{(r)} + T_\beta^{(s)} = 0$. \endproof

\bibliographystyle{unsrt}

\end{document}